\documentclass[a4paper,11pt,english]{amsart} 
\usepackage{amssymb,amsmath,amsfonts,latexsym,color,graphicx,stmaryrd,dsfont}
\usepackage{enumerate}
\usepackage{mathtools}
\usepackage{graphicx}
\usepackage{float}

\usepackage{hyperref} 
\hypersetup{pdfborder=0 0 0, 
	    colorlinks=true,
	    citecolor=blue,
	    linkcolor=blue,
	    urlcolor=blue,
	    pdfauthor={Johannes Sjoestrand, Martin Vogel}
	   }
\usepackage[font=small,labelfont=bf]{caption}
\usepackage[margin=2.5cm]{geometry}
\numberwithin{equation}{section}
\newcommand{\defeq}{\stackrel{\mathrm{def}}{=}}
\newcommand{\C}{\mathds{C}}
\newcommand{\R}{\mathds{R}}
\newcommand{\N}{\mathds{N}}
\newcommand{\Z}{\mathds{Z}}
\newcommand{\spec}{\mathrm{Spec}}
\newcommand{\supp}{\mathop{\rm supp}} 
\newcommand{\tr}{\mathrm{tr}}
\newcommand{\dist}{\mathrm{dist}}
\newcommand{\coker}{\mathrm{coker}}
\newcommand{\e}{\mathrm{e}}
\newcommand{\HS}{\mathrm{HS}}

\newcommand{\mO}{\mathcal{O}}

\newtheorem{dref}{Definition}[section] 
\newtheorem{lemma}[dref]{Lemma}
\newtheorem{theo}[dref]{Theorem} 
\newtheorem{prop}[dref]{Proposition}
\newtheorem{remark}[dref]{Remark} 

\newtheorem{cor}[dref]{Corollary}
\title[Toeplitz band matrices with small random perturbations]
{Toeplitz band matrices with small random perturbations}
\author{Johannes Sj\"ostrand}
\address[Johannes Sj\"ostrand]{IMB, 
  Universit\'e de Bourgogne Franche-Comt\'e, 
  UMR 5584 du CNRS, 
  9, avenue Alain Savary - BP 47870 FR-21078 Dijon Cedex, France.}
\email{johannes.sjostrand@u-bourgogne.fr}
\author{Martin Vogel}
\address[Martin Vogel]{Institut de Recherche Math{\'e}matique Avanc{\'e}e - UMR 7501 CNRS, Universit{\'e} de Strasbourg, 7 rue René-Descartes, 67084 Strasbourg Cedex, France.}
\email{vogel@math.unistra.fr}
 \date{}
\begin{document}
\begin{abstract}
We study the spectra of $N\times N$ Toeplitz band matrices perturbed by small complex Gaussian 
random matrices, in the regime $N\gg 1$. We prove a probabilistic Weyl law, which provides an precise asymptotic formula for the number of eigenvalues in certain domains, which may depend on $N$, 
 with probability sub-exponentially (in $N$) close to $1$. We show that most eigenvalues of the 
 perturbed Toeplitz matrix are at a distance of at most 
 $\mO(N^{-1+\varepsilon})$, for all $\varepsilon >0$, to the curve in the complex plane given by the symbol of the unperturbed Toeplitz 
 matrix. 
\end{abstract}
\maketitle
\setcounter{tocdepth}{1}
\tableofcontents
\section{Introduction} 
Let $N_{\pm} \geq 0$ be in $\N$, such that either $N_+\neq 0$ or $N_-\neq 0$, 
and consider the operator
\begin{equation}\label{int.0}
p(\tau) \defeq \sum_{j=-N_-}^{N_+} a_j \tau^j, \quad a_{-N_-},
a_{-N_-+1},\dots, a_{N_+} \in \C,\ a_{\pm N_\pm}\ne 0, 
\end{equation}
acting on $\ell^2(\Z)$, or more generally on functions $\psi: \Z\to \C$, where 
\begin{equation}\label{int.1}
 (\tau u)(k) \stackrel{\mathrm{def}}{=} u(k-1),
\end{equation}
defines the translation to the right by one unit. We shall work on
$\Z$, on an interval in $\Z$ and on $\Z / M\Z$, for some $ \N \ni M
\geq 1$. The symbol of $\tau =\exp( -iD_x)$ is $1/\zeta $, with $\zeta
=e^{i\xi }$. Therefore, the symbol of the operator \eqref{int.0} is
given by the meromorphic function
\begin{equation}\label{int.0.1}
\C \ni \zeta \mapsto p(1/\zeta) = \sum_{j=-N_-}^{N_+} a_j \zeta^{-j}.
\end{equation}
\par
We obtain a \emph{Toeplitz band matrix} from the operator $p(\tau)$ by restricting it 
to the finite dimensional space $\C^N$. Indeed, we let $N\geq 1$ and identify $\C^N$ 
with $\ell^2([1,N])$, $[1,N]=\{ 1,2,..,N\}$, and also with $\ell^2_{[1,N]}({\Z})$ (the space 
of all $u\in \ell^2({\Z})$ with support in $[1,N]$). Then, we consider the $N\times N$ 
\emph{Toeplitz band matrix} 
\begin{equation}\label{int.2}
P_N\defeq 1_{[1,N]}\,p(\tau)\,1_{[1,N]},
\end{equation}
acting on $\C^N \simeq \ell^2_{[1,N]}({\Z})$.  
\\
\par
\begin{figure}[ht]
  \centering
   \begin{minipage}[b]{0.49\linewidth}
  \centering
   \includegraphics[width=\textwidth]{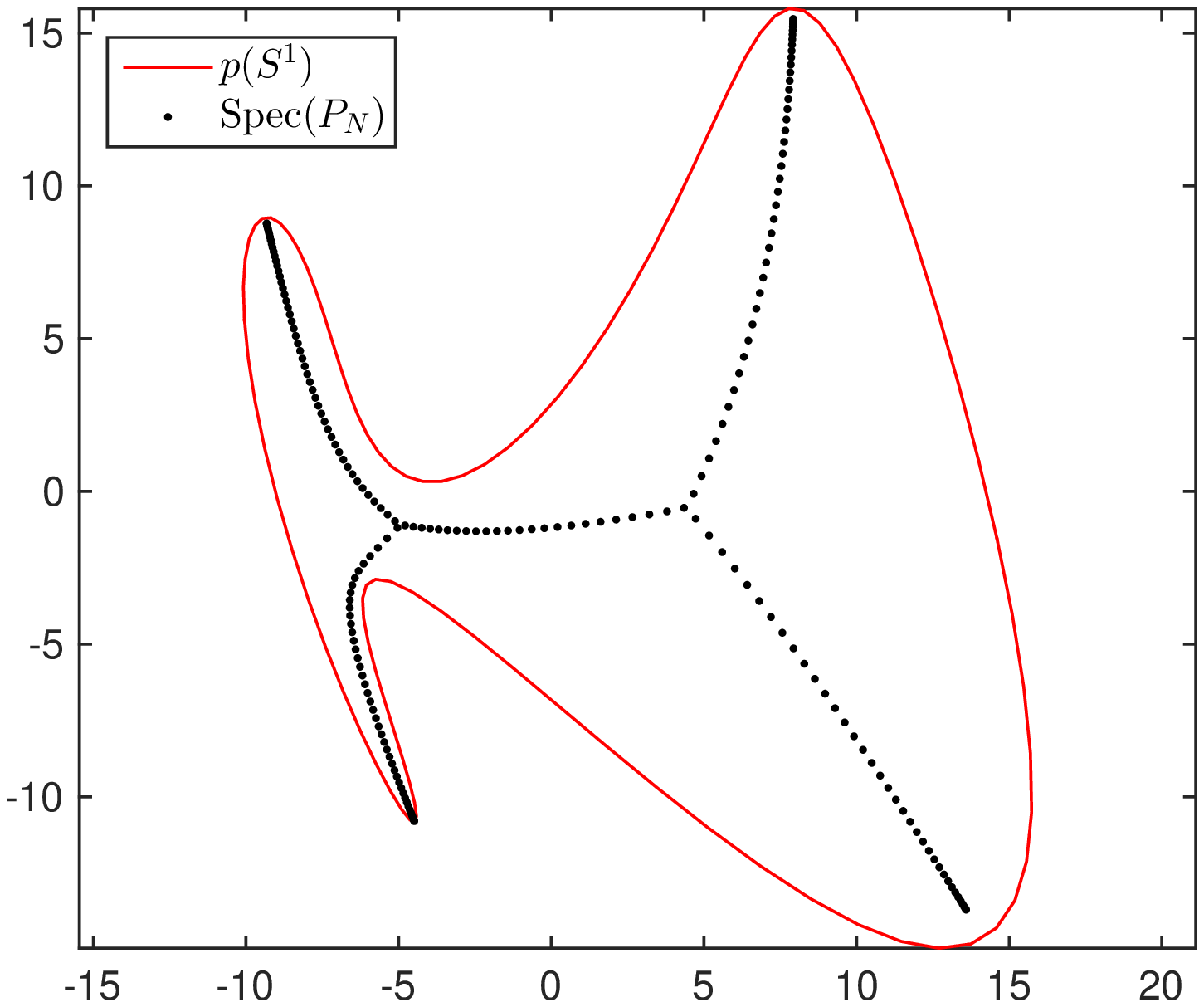}
 \end{minipage}
 \hspace{0cm}
 \begin{minipage}[b]{0.49\linewidth}
  \includegraphics[width=\textwidth]{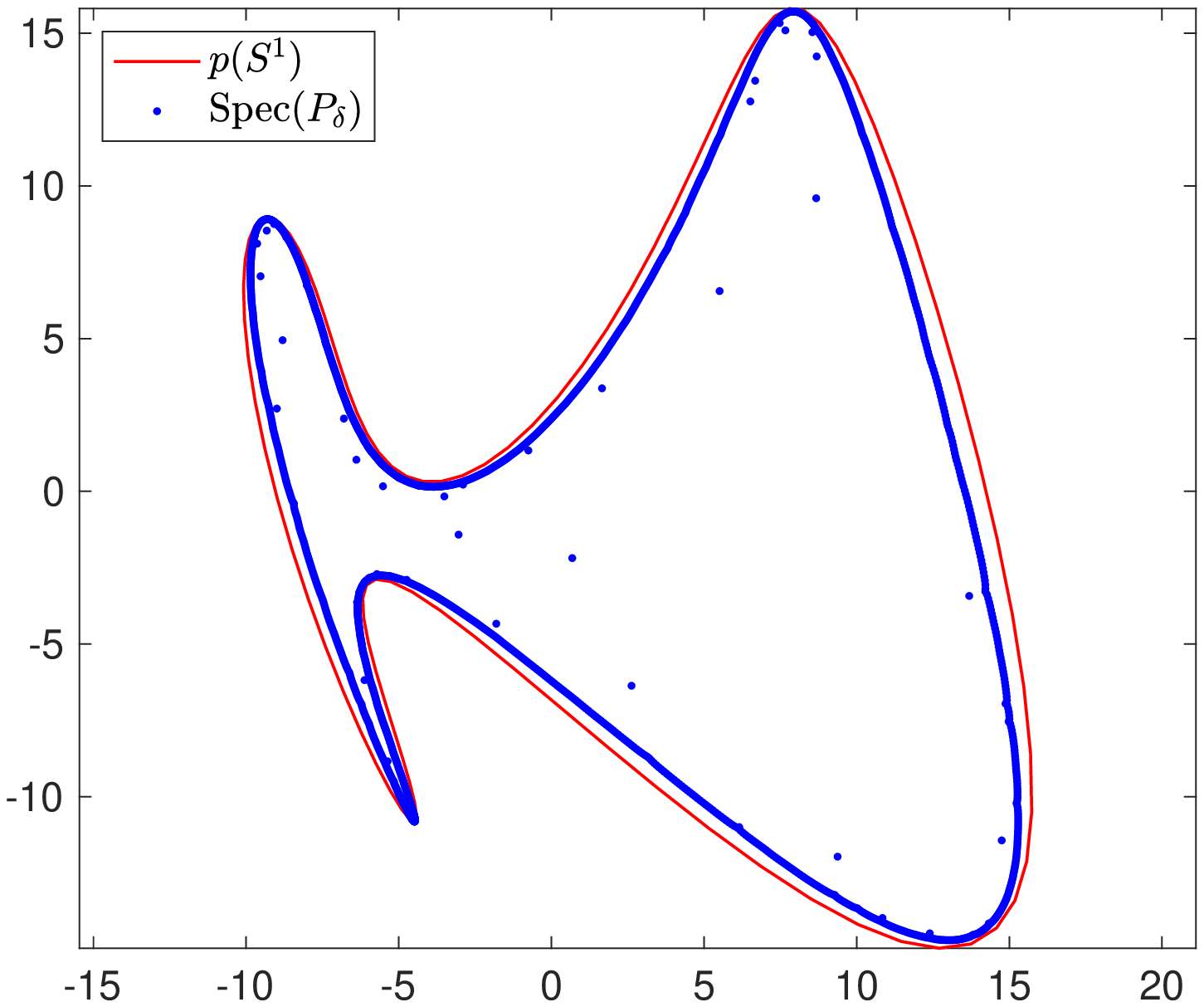}
 \end{minipage}
 \caption{The pictures on the left hand side shows the spectrum of the Toeplitz matrix 
 	$P_N$, with $N=100$, 
	given by the symbol $p(1/\zeta)=2i\zeta^{-1}+\zeta^2+\frac{7}{10}\zeta^3$ 
	and the right hand side shows the spectrum of a random perturbation $P_{\delta}$, 
	as in \eqref{i0.1} below, with coupling constant $\delta=10^{-14}$ and dimension $N=1000$. 
	The red line shows the symbol curve $p(S^1)$.}
  \label{fig1}
\end{figure}

The translation operator $\tau$ on $\ell^2(\Z)$ is unitary, i.e. $\tau^*=\tau^{-1}$, 
so one can easily see that $p(\tau)$ is a normal operator, meaning that it 
commutes with its adjoint. The Fourier transform shows that the spectrum of 
$p(\tau)$ \eqref{int.0} acting on $\ell^2(\Z)$ is purely absolutely continuous 
and given by 
\begin{equation}\label{int.2.2}
	\spec(p(\tau)) = p(S^1).
\end{equation}
The restriction $P_{\N}=p(\tau)|_{\ell^2(\N)}$ of $p(\tau)$ to $\ell^2(\N)$, is in general no longer 
normal, except for specific choices of $N_+,N_-$ and the coefficients 
$a_j$. The essential spectrum of the Toeplitz operator $P_{\N}$ \eqref{int.0} is still 
given by $p(S^1)$. However, we gain additional pointspectrum in all loops of $p(S^1)$ 
with non-zero winding number, i.e.  
\begin{equation}\label{int.2.3}
	\spec(P_{\N}) = p(S^1) \cup \{ z\in \C; \mathrm{ind}_{p(S^1)}(z)\neq 0 \}.
\end{equation}
Here, by a result of Krein \cite[Theorem 1.15]{BoSi99} (see also Proposition \ref{g_prop4.1} below) 
the winding number of $p(S^1)$ around the point $z\not\in p(S^1)$ is related to the Fredholm 
index of $P_{\N}-z$: 
\begin{equation}\label{int.2.4}
	\mathrm{Ind}(P_{\N}-z) = - \mathrm{ind}_{p(S^1)}(z).
\end{equation}
For every $\epsilon >0$, the spectrum of the finite Toeplitz matrix $P_N$ \eqref{int.2} satisfies 
\begin{equation}\label{int.2.5}
	\spec(P_{N}) \subset \spec(P_{\N})+D(0,\epsilon )
\end{equation}
for $N>0$ sufficiently large, where $D(z,r)$ denotes the open disc of
radius $r$, centered at $z$. The limit of $\spec(P_{N})$ as $N\to
\infty$ is contained in a union of analytic 
arcs inside $ \spec(P_{\N})$, see \cite[Theorem 5.28]{BoSi99}. This phenomenon can also be 
observed in the numerical simulations presented in Figures \ref{fig1}, \ref{fig2}. 
\par
However, we will show that after a small random perturbation of $P_N$, most of the eigenvalues of the
perturbed operator will be very close to the curve $p(S^1)$, see Figures \ref{fig1}, \ref{fig2}. 
\subsection{Adding a small random perturbation}
Let $(\mathcal{M},\mathcal{A},\mathds{P})$ denote a probability space 
and let $\mathcal{H}_N(\C^{N\times N}, \| \cdot \|_{\HS})$ denote the space 
of $N\times N$ complex valued matrices equipped with the Hilbert-Schmidt norm. 
Consider the random matrix 
$$
\mathcal{M} \ni \omega \mapsto Q_{\omega}\defeq Q_{\omega}(N) \defeq  (q_{j,k}(\omega))_{1\leq j,k\leq N} 
\in \mathcal{H}_N
$$
with Gaussian law 
\begin{equation*}
	(Q_{\omega})_*(d\mathds{P}) = \pi^{-N^2} \e^{-\|Q\|_{\HS}^2} L(dQ),
\end{equation*}
where $L$ denotes the Lebesgue measure on $\C^{N\times N}$. 
We are interested in the spectrum of the random perturbations of 
the matrix $P_N^0=P_N$: 
\begin{equation}\label{i0.1}
 P_N^{\delta} \defeq P_N^0 + \delta Q_{\omega}, 
 \quad 0 \leq \delta \ll 1.
\end{equation}
Notice that the entries $q_{j,k}(\omega)$ of $Q_{\omega}$ are independent and 
identically distributed complex Gaussian random variables with 
expectation $0$, and variance $1$. 
\par 
We recall that the probability distribution of a complex Gaussian 
random variable $\alpha \sim \mathcal{N}_{\C}(0,1)$, defined on the  
probability space $(\mathcal{M},\mathcal{A},\mathds{P})$, is given by 
\begin{equation*}
	\alpha_*(d\mathds{P}) = \pi^{-1} \e^{-|\alpha|^2} L(d\alpha),
\end{equation*}
where $L(d\alpha)$ denotes the Lebesgue measure on $\C$. 
If $\mathds{E}$ denotes the expectation with respect to the probability 
measure $\mathds{P}$, then 
\begin{equation*}
	\mathds{E}[\alpha] = 0, \quad \mathds{E}[|\alpha|^2] = 1.
\end{equation*}
In this paper we consider the Gaussian case for the sake of simplicity. However, we believe that 
our method can be adapted to the case of more general complex valued random 
matrices. The main difficulty lies in showing that the logarithm of the determinant of 
a certain matrix valued stochastic process is not too small with probability close to $1$ (see Proposition 
\ref{prop:LB} below). 
\section{Main results}\label{int}
\setcounter{equation}{0}
We will provide precise eigenvalue asympotics for the eigenvalues of $P_N^{\delta}$ 
 in certain domains which show that most eigenvalues of $P_N^{\delta}$ 
 are close to the curve $p(S^1)$ with probability sub-exponentially (in $N$) close to $1$, 
 see Theorem \ref{thm:t1} below. We also prove eigenvalue asymptotics in \emph{thin} 
 $N$-dependent domains in scales up to order $N^{-1+\varepsilon}$, for every $\varepsilon >0$. 
 This shows in particular that for every $\varepsilon >0$, with probability sub-exponentially 
 (in $N$) close to $1$, 
 most eigenvalues can be found at a distance $\leq \mO(N^{-1+\varepsilon})$ 
 from $p(S^1)$, see Theorem \ref{thm:t2} for the precise statement. 
 \par
Our results also provide an upper bound on the number of eigenvalues of $P_N^{\delta}$ 
which remain far from the curve $p(S^1)$. Finally, we will show that our results 
on the eigenvalue asymptotics of $P_N^{\delta}$ imply the almost sure weak 
convergence of the empirical measure of eigenvalues of $P_N^{\delta}$ to the 
uniform measure on $p(S^1)$, see Corollary \ref{thm:t0}. This corresponds to 
the leading term of our asymptotic result. 
\subsection{Eigenvalue asymptotics in fixed smooth domains}\label{sec:2.R2}
Let $\Omega \Subset \C$ be an open 
simply connected set with smooth boundary $\partial\Omega$ which is independent 
of $N$. We suppose that 
\begin{enumerate}[($\Omega$1)]
	\item \label{O1} $\partial\Omega$ intersects  $p(S^1)$ in at most finitely many 
		points;
	\item \label{O2} the points of intersection are non-degenerate, i.e. 
		\begin{equation}\label{fd.1}
			\partial_{\zeta}p \neq 0 \hbox{ on } p^{-1}(\partial\Omega \cap p(S^1) );
		\end{equation}
              \item \label{O3} $\partial\Omega$ intersects $p(S^1)$
                transversally, in the following precise sense :
		for each $z_0\in\partial\Omega \cap p(S^1)$ let $\gamma_k\subset p(S^1)$, 
		$k=1,\dots, n$ denote the mutually distinct segments of $p(S^1)$ passing 
		through $z_0$, i.e. 
		each $\gamma_k$ is given by the image of a small neighborhood in $S^1$ 
		of a point in $p^{-1}(z_0)\cap S^1$. Then $\gamma _k$
                and $\partial \Omega $ intersect transversally at $z_0$. 
\end{enumerate}
We then have the following result:
\begin{theo}\label{thm:t1}
  Let $p$ be as in \eqref{int.0}, set $M = N_+ + N_-$ and let
  $P_N^{\delta}$ be as in \eqref{i0.1}. Let $\Omega$ be as above,
  satisfying conditions
  \textnormal{($\Omega$\ref{O1})--($\Omega$\ref{O3})} and pick a
  $\varepsilon_0\in]0,1[$. There exists a constant $C>0$, such that, for
  $N>1$ sufficiently large, if 
\begin{equation}\label{int.3}
	C\e^{- N^{\varepsilon_0}/(2M)}\leq \delta \leq \frac{N^{-4}}{C}, 
\end{equation}
then we have that 
 \begin{equation}\label{t1.1}
	\bigg| 
	\#(\spec(P^{\delta}_N)\cap \Omega )
	-  \frac{N}{2\pi} \int_{p^{-1}(\Omega\cap p(S^1))}L_{S^1}(d\theta) \bigg| 
	\leq 
	\mO( N^{\varepsilon_0}\log N ) .
\end{equation}
with probability 
 \begin{equation}\label{t1.2}
	 \geq 1 - 
	 \mO(\log N)
	 \left(\e^{-N^2} +\delta^{-M} \e^{-\frac{1}{2}N^{\varepsilon_0}}\right).
\end{equation}
\end{theo}
Let us give some remarks on this result. The $\e^{-N^2}$ term in the estimate \eqref{t1.2} is 
an artifact from the proof where we restrict to the event that $ \| Q_{\omega}\|_{\HS} \leq C N$ 
which occurs with probability $\geq 1 - \e^{-N^2}$, see \eqref{mark1}. In fact, in the proof 
we can reduce this restriction to $ \| Q_{\omega}\|_{\HS} \leq C \sqrt{N}$ which results in \eqref{t1.1} 
holding with probability \eqref{t1.2} with $\e^{-N^2}$ exchanged by $\e^{-N}$. Moreover, the 
Theorem holds for $C\e^{- N^{\varepsilon_0}/(2M)}\leq \delta \leq \frac{N^{-5/2}}{C}$.
\par
The factor $N^{\varepsilon_0}$ in the error estimate in \eqref{t1.1} is a consequence 
of our aim to show that \eqref{t1.1} holds with probability which is sub-exponentially 
close to $1$. However, it is clear from 
the proof, see Proposition \ref{prop:LB}, that if we were to settle for a probability 
$\geq 1 - N^{-\kappa}$, for every $\kappa >0$, 
then we can ameliorate the error estimate in \eqref{t1.1} to $\mO((\log N)^2)$. 
\begin{figure}[ht]
 \begin{minipage}[b]{0.49\linewidth}
  \centering
   \includegraphics[width=\textwidth]{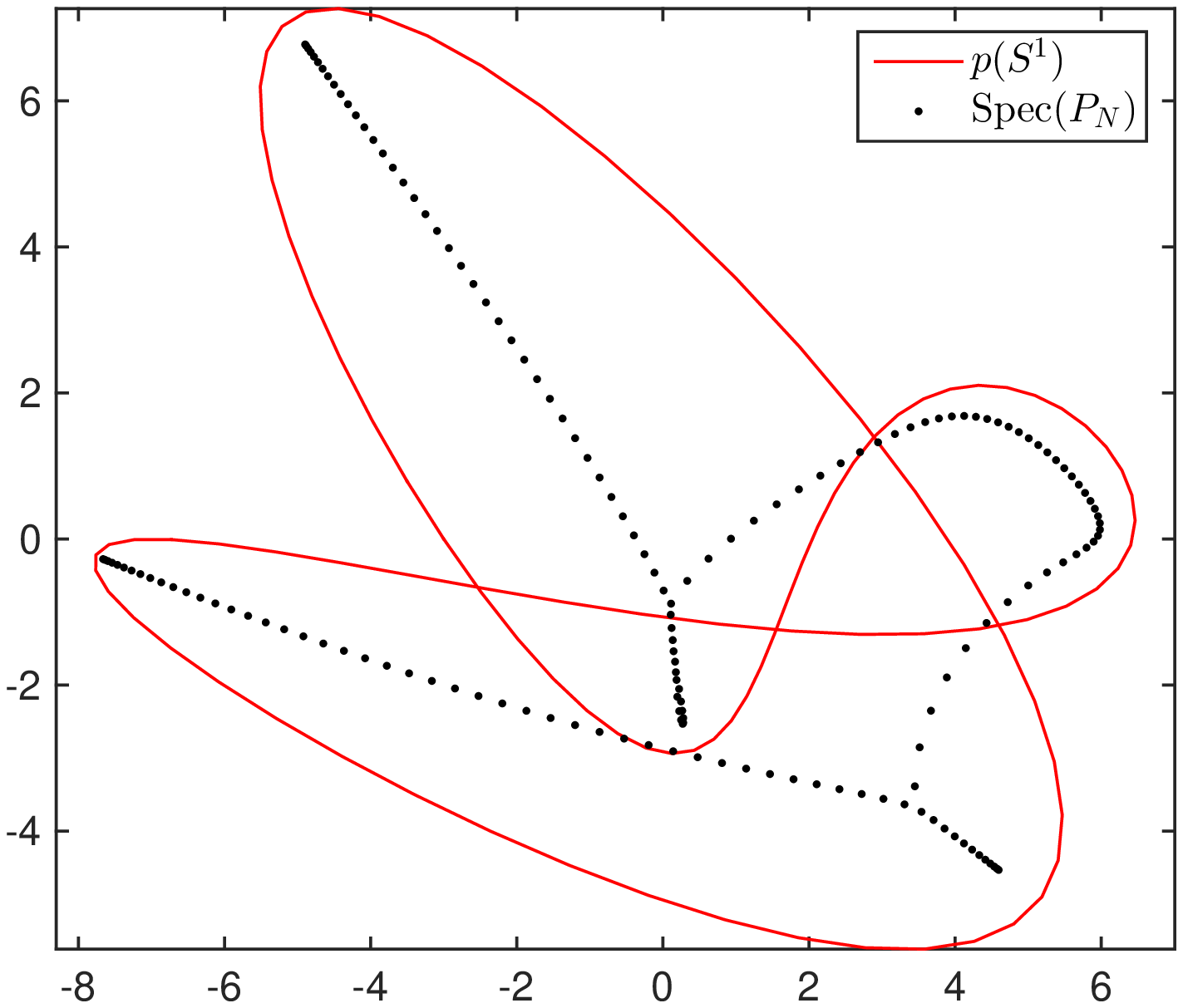}
 \end{minipage}
 \hspace{0cm}
 \begin{minipage}[b]{0.49\linewidth}
  \includegraphics[width=\textwidth]{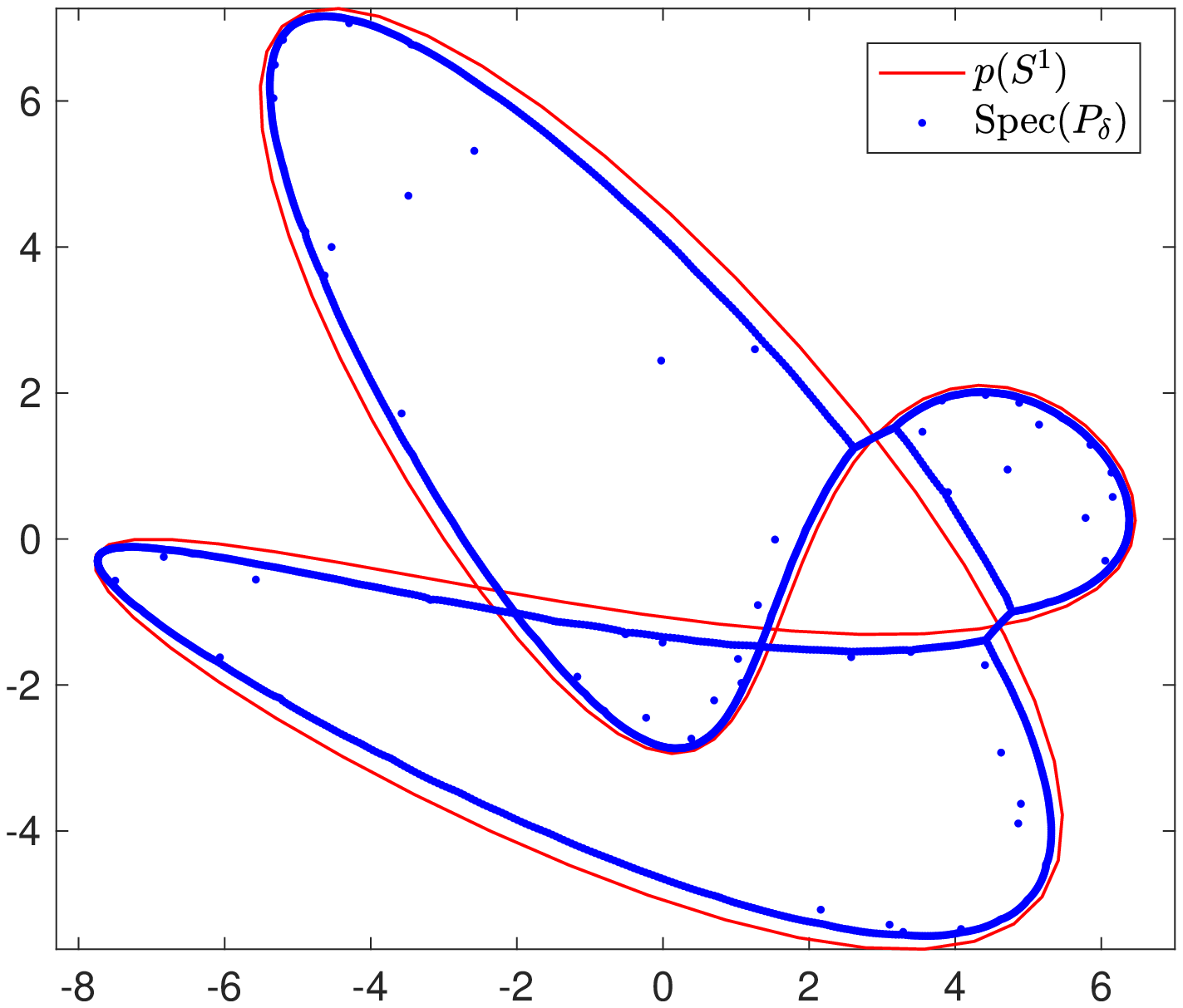}
 \end{minipage}
 \caption{The pictures on the left hand side shows the spectrum of the Toeplitz matrix $P_N$, 
 	with $N=100$,  	given by the symbol 
	$p(1/\zeta)= 2\zeta^{-3}-\zeta^{-2}+2i\zeta^{-1}-4\zeta^2-2i\zeta^{3}$ and 
	the right hand side shows the spectrum of a random 	perturbation $P_{\delta}$, 
	as in \eqref{i0.1}, with coupling 
	constant $\delta=10^{-14}$ and $N=1000$. The red line shows the symbol curve $p(S^1)$.}
  \label{fig2}
\end{figure}
\par
We provide a more detailed version of this result in Theorem \ref{thm:t2} below. There, 
we present a Weyl law in probability for the eigenvalues of $P_N^{\delta}$ in \emph{thin} 
$N$-dependent domains $\Omega_N$ with, roughly speaking, a width $\geq CN^{-1+}$, 
and whose boundary is uniformly Lipschitz. See Theorem \ref{thm:t2} below for more details. 
\subsection{Convergence of the empirical measure and related results}
Another way to see the limiting behavior of the spectrum of $P_N^{\delta}$ \eqref{i0.1}  is 
to study the limits of the \emph{empirical measure} of the eigenvalues of $ P_N^{\delta}$, 
defined by 
\begin{equation}\label{eq:cm1}
	\xi_N \defeq \frac{1}{N}\sum_{\lambda \in \spec( P_N^{\delta})} \delta_{\lambda}
\end{equation}
where the eigenvalues are counted including multiplicity and $\delta_{\lambda}$ denotes 
the Dirac measure at $\lambda \in\C$. The Markov inequality implies that 
\begin{equation}\label{mark1}
	\mathds{P} [ \| Q_{\omega}\|_{\HS} \leq C N ] \geq 1 - \e^{-N^2},
\end{equation}
for $C>0$ large enough. The operator norm of $P_N$ \eqref{int.2} satisfies
\begin{equation*}
	\| P_N\| \leq \| p \|_{L^{\infty}(S^1)}. 
\end{equation*}
If $\delta \leq N^{-1}$, then the Borel-Cantelli Theorem shows that, almost surely, 
$\xi_N$ has compact support for $N>0$ sufficiently large. 
\par
From Theorem \ref{thm:t1} we will deduce that, almost surely, $\xi_N$ converges weakly to the 
uniform distribution on $p(S^1)$.  
\begin{cor}\label{thm:t0}
	Let $\varepsilon_0 \in ]0,1[$, let $p$ be as in \eqref{int.0} and write $M = N_+ + N_-$. Then, 
	there exists a constant $C>0$ such that if (\ref{int.3}) holds,
	\begin{equation*}
	C\e^{-N^{\varepsilon_0}/(2M)} \leq \delta \leq \frac{N^{-4}}{C},
	\end{equation*}
	then, almost surely, 
	\begin{equation}\label{int.4}
		\xi_N \rightharpoonup  p_*\left(\frac{1}{2\pi} L_{S^1}\right), \quad N\to \infty, 
	\end{equation}
	weakly, where $L_{S^1}$ denotes the Lebesgue measure 
	on $S^1$.
\end{cor}
Our strategy to prove the precise eigenvalue asymptotics presented in Theorems \ref{thm:t1} 
and \ref{thm:t2} also provides an alternative proof of the above result via the convergence of the 
associated logarithmic potentials, see Section \ref{sec:LP}.
\\
\par
Similar results to Corollary \ref{thm:t0} have been proven in various settings. 
In the recent work \cite{BPZ18}, the authors consider a general sequence of deterministic 
complex $N\times N$ matrices $M_N$ perturbed by complex Gaussian random 
matrices $Q_{\omega} = Q_{\omega}(N)$, as in \eqref{i0.1}. They study the 
empirical measure $\xi_N$ of the eigenvalues of $\mathcal{M}_N := M_N + N^{-\gamma} Q_{\omega}$, 
$\gamma >1/2$, defined as in \eqref{eq:cm1}. The authors show 
that the \emph{Logarithmic potential} $L_{\xi_N}(z)$, $z\in\C$, (see Section \ref{sec:LP} below for a definition) associated with $\xi_N$, asymptotically coincides with a deterministic function 
$g_N(z)$ in probability at each point $z$, for which the number of singular values of $(M_N -z\mathrm{Id})$ smaller than $N^{-\gamma +1/2 + \delta_N}$, $0<\delta_N = o(1)$ as $N\to \infty$, is of order $o(N(\log N)^{-1})$  as $N\to \infty$. Since the weak convergence of the random measure $\xi_N$ can be deduced from the point wise convergence of the 
Logarithmic potential $L_{\xi_N}(z)$ (see Section \ref{sec:LP} below for details and references), 
this result shows that studying the weak convergence of the empirical measure $\xi_N$ 
can be reduced to deterministic calculation involving only the unperturbed matrix $M_N$. 
\par
Moreover, in \cite{BPZ18,BPZ18b}, the authors consider the special case of $M_N$ being given by 
a band Toeplitz matrix, i.e. $M_N = P_N$ with $p$ as in \eqref{int.0}. 
In this case they show that the convergence \eqref{int.4} holds weakly in probability for 
a coupling constant $\delta = N^{-\gamma}$, with $\gamma  >1/2$. Furthermore, 
they prove a version of this theorem for Toeplitz matrices with non-constant coefficients in 
the bands, see \cite[Theorem 1.3, Theorem 4.1]{BPZ18}. Their methods are quite different 
from ours. They compute directly the $\log |\det \mathcal{M}_N -z|$ by relating it to 
$\log |\det M_N(z)|$, where $M_N(z)$ is a truncation of $M_N -z$, where the smallest singular 
values of $M_N-z$ have been excluded. The level of truncation however depends on the strength 
of the coupling constant and it necessitates a very detailed analysis of the small singular values 
of $M_N -z$.
\\
\par
In the earlier work \cite{GuMaZe14}, the authors prove that the convergence \eqref{int.4} holds weakly in probability for the Jordan bloc matrix $P_N$ with $p(\tau) = \tau^{-1}$ \eqref{int.0} and 
 a perturbation given by a complex Gaussian random matrix whose entries 
are independent complex Gaussian random variables whose variances vanishes (not necessarily at the same speed) polynomially fast, with minimal decay of order $N^{-1/2+}$. 
\par
In \cite{Wo16}, using a replacement principle developed in \cite{TVK10}, it was shown that the result of \cite{GuMaZe14} holds for perturbations given by 
complex random matrices whose entries are independent and identically distributed random 
complex random variables with expectation $0$ and variance $1$ and a coupling constant 
$\delta = N^{-\gamma}$, with $\gamma >2 $. 
\\
\par
In \cite{DaHa09}, the authors showed that in the case of large Jordan block matrix $p(\tau) = \tau^{-1}$, most eigenvalues of the perturbed matrix $P_N^{\delta}$ lie in the annulus 
$$\{z\in \C; (\delta N)^{1/N} \e^{-\sigma}  \leq |z| \leq (\delta N)^{1/N}\},$$
for any fixed $\sigma >0$, 
with probability $\geq 1 - \mO(N^{-2})$. Moreover, the authors show that there are at most 
$\mO(\sigma^{-1}\log N )$ eigenvalues of $P_N^{\delta}$  outside this annulus, 
with probability $\geq 1 - \mO(N^{-2})$.
\par
A version of Theorem \ref{thm:t1}, concerning the special cases of large Jordan block matrices 
$p(\tau) = \tau^{-1}$ and large bi-diagonal matrices $p(\tau) = a\tau + b\tau^{-1}$, $a,b\in \C$, 
have been proven in \cite{Sj19,SjVo15b}. 
\subsection{Spectral instability}\label{sec:2.D1}
In general, the spectra of non-selfadjoint operators can be highly unstable under 
small perturbations due to the lack of good control over the norm 
of the resolvent. This phenomenon, sometimes referred to as \emph{pseudospectral effect} or 
\emph{spectral instability}, can be observed in the case of non-normal Toeplitz matrices 
$P_N$ \eqref{int.2}, as illustrated in Figures \ref{fig1} and \ref{fig2}. To quantify the 
zone of spectral instability in the complex Plane, one defines the 
$\varepsilon$-\emph{pseudospectrum} of a linear operator $P$ acting on some complex
 Hilbert space $\mathcal{H}$ as follows: for $\varepsilon>0$ set 
\begin{equation}
	\spec_{\varepsilon}(P) \defeq \{ z\in \C; \|(P-z)^{-1}\| > \varepsilon^{-1}\}. 
\end{equation}
The points $z\in\C$ in the $\varepsilon$-{pseudospectrum} of $P$ are precisely the points $z\in\C$ for whom there exists a bounded linear operator $Q$ acting on $\mathcal{H}$ with 
$\|Q\|\leq 1$, such that $z\in \spec(P+\varepsilon Q)$, see \cite{TrEm05,Da07} for a detailed 
exposition. 
\par
For the Toeplitz band matrices $P_N$, we have that any fixed point in $\C\backslash p(S^1)$ with
\begin{equation}\label{sc1}
	z\notin\{0,+\infty\} \quad \text{and} \quad z\neq a_0, \text{ when } 
	N_+ \text{ or } N_- =0,
\end{equation}
which is contained in the pointspectrum of $P_{\N}$ \eqref{int.2.3} is contained in the 
$C\e^{-N/C}$-pseudospectrum of $P_N$. Recall from \eqref{int.2.3} that the pointspectrum of $P_{\N}$ 
in $\C\backslash p(S^1)$ is given by the points $z$ around which the curve $p(S^1)$ has a non-zero 
winding number $\mathrm{ind}_{p(S^1)}(z)\neq 0$. In fact, provided that we avoid the special cases \eqref{sc1}, we have that 
\begin{itemize}
\item if $\mathrm{ind}_{p(S^1)}(z)<0 $, then the Fredholm index of $P_{\N}-z$ satisfies 
	\begin{equation*}
	 	\mathrm{Ind}(P_{\N}-z) = \dim \mathrm{ker} (P_{\N}-z) = - \mathrm{ind}_{p(S^1)}(z);
	\end{equation*}
\item if $\mathrm{ind}_{p(S^1)}(z)>0 $,
	\begin{equation*}
	 	\mathrm{Ind}(P_{\N}-z) = - \dim \mathrm{ker} (P_{\N}-z)^* 
		= - \mathrm{ind}_{p(S^1)}(z),
	\end{equation*}
\end{itemize}
see Propositions \ref{g_prop3} and \ref{g_prop4.1}. 
Moreover, these kernels are spanned by exponentially decaying functions, see the discussion in Section \ref{gendisc_s1}. In the first case, restricting such a function $u \in \mathrm{ker} (P_{\N}-z)$ to the interval $[1,N]$ yields 
an approximate solution to the equation $(P_{N}-z)u =0$, sometimes called a \emph{quasimode}. More precisely, setting 
$e_+ = \|\mathbf{1}_{[1,N]}  u \|^{-1}\mathbf{1}_{[1,N]} u $, we get that 
	\begin{equation*}
	 	(P_{N}-z)e_+ = \mO(\e^{-N/C}). 
	\end{equation*}
Similarly, we get in the second case an $e_-\in \ell^2([1,N])$, $\|e_-\| =1$, with 
	\begin{equation*}
	 	(P_{N}-z)^*e_- = \mO(\e^{-N/C}). 
	\end{equation*}
These exponentially precise quasimodes show that any fixed $z$ with 
$\mathrm{ind}_{p(S^1)}(z)\neq 0$ satisfying \eqref{sc1}, is contained in the 
$C\e^{-N/C}$-pseudospectrum of $P_N$. 
\\
\par
On the other hand, for any compact set $\Omega\Subset \C\backslash p(S^1)$, with 
$z\in\Omega$ satisfying \eqref{sc1} and 
	\begin{equation*}
	 	\mathrm{ind}_{p(S^1)}(z) = 0,
	\end{equation*}
 we have that that for $N>0$ sufficiently large $\| (P_N-z)^{-1}\| = \mO(1)$ uniformly 
 for $z\in\Omega$, see Proposition \ref{resUB}. 
 Hence, outside the spectrum of $P_{\N}$ \eqref{int.2.3}
is a zone of spectral stability for $P_N$. This explains why the eigenvalues 
of $P_N^{\delta}$ can (with high probability) only be found in a small neighborhood 
of  the spectrum of $P_{\N}$. 
\par
However, only analysing the pseudospectrum does not 
yield any information on \emph{where} the eigenvalues of $P_N^{\delta}$ can be found. 
Theorem \ref{thm:t1}, shows that with probability very close to one, all but 
$\mO(N^{\varepsilon_0}\log N)$ many eigenvalues of $P_N^{\delta}$ can be found close to 
the curve $p(S^1)$. Theorem \ref{thm:t2} below shows that still probability very close to one, 
most eigenvalues of $P_N^{\delta}$ are at a distance of $\leq N^{-1+\varepsilon}$, for 
every $\varepsilon>0$, from $p(S^1)$, see \eqref{t2.1} for the precise error estimate. 
\par
It would be interesting to perform a precise analysis of the boundary 
of the $\varepsilon$-pseudospectrum of $P_N$ to see whether the eigenvalues of $P_N^{\delta}$  accumulate there, as in the case of small random perturbations of semiclassical differential 
operators in \cite{Vo14}.
\subsection{Outline of the proof}
The overall strategy of the proof is based on a \emph{Grushin} reduction. In Section \ref{grush} 
we review the basic idea of such a reduction and we set up a Grushin problem $\mathcal{P}_N$ by considering the operator $p(\tau)$ \eqref{int.0} on 
the discrete circle $\Z/\widetilde{N}\Z$, $\widetilde{N}=N+N_-+N_+$,
$$
	\mathcal{P}_N = p(\tau): \ell^2(\Z/\widetilde{N}\Z) \to \ell^2(\Z/\widetilde{N}\Z) ,
$$
which can be used to describe the eigenvalues of the unperturbed operator $P_N$. In Section \ref{gendisc} we provide a general discussion of band Toeplitz matrices 
and their Fredholm properties. However, for this paper only Sections \ref{gendisc_s2} and \ref{sec:SpG} are of immediate importance as we discuss properties of $p(\tau)$ on $\Z/\widetilde{N}\Z$.
\par
In Section \ref{sec:GP}, we will use the Grushin problem for the unperturbed operator $P_N$ 
to set up a \emph{Grushin Problem} $\mathcal{P}_N^{\delta}$ for the perturbed 
$P_N^{\delta}$, resulting in an effective description of its eigenvalues 
$$
	\log \det (P_N^{\delta}-z) = \log \det \mathcal{P}_N^{\delta}(z) + \log \det E_{-+}^{\delta}(z), 
$$
with probability $\geq 1 - \e^{-N^2}$. Here, $E_{-+}(z)^{\delta}$ is an $(N_++N_-)\times(N_++N_-)$ complex valued matrix. Furthermore, the Grushin problem shows that we have a trivial upper bound 
on the quantity $\log \det E_{-+}^{\delta}(z)$. In Section \ref{sec:LB}, we show that with probability 
very close to $1$ we have a quantitative lower bound on $\log \det E_{-+}^{\delta}(z)$. 
\par
To obtain our main results on eigenvalue asympotics from this description we apply  
a general estimate \cite{Sj09b} on the number of zeros of a holomorphic function 
$u(z;N)$ of exponential growth. We will recall 
this result in Section \ref{sec:c1} below, see Theorem \ref{thm:Count}. Roughly speaking, 
if the available information is
\begin{enumerate}[(i)]
\item an \emph{upper bound} $\log|u(z;N)| \leq N\phi(z)$, for $z$ near the boundary 
	$\partial \Omega$ and $\phi$ a subharmonic continuous function and 
\item a \emph{lower bound} $\log|u(z;N)| \geq N(\phi(z)-\varepsilon_j)$, with 
	$\varepsilon_j \geq 0$, for finitely many points $z_j = z_j(N)$, $j=1,\dots, M(N)$, 
	which are situated near the boundary of $\partial\Omega$,
\end{enumerate}
then the number of zeros of $u$ in $\Omega$ is given by 
$$
\#(u^{-1}(0)\cap \Omega) \sim \frac{N}{2\pi} \int_{\Omega} \Delta\phi L(dz),
$$
asymptotically as $N\to +\infty$. In Section \ref{sec:c2} we check that our effective 
description for $\log \det (P_N^{\delta}-z)$ satisfies the required upper bound (i), 
and in Section \ref{sec:c3}, using Section \ref{sec:LB}, we check the lower bound (ii). 
\par
 In Section \ref{sec:c4} we then use these bounds in combination with 
Theorem \ref{thm:Count} to prove Theorem \ref{thm:t1}.
\par
In Section \ref{sec:c5} we provide a more general version of Theorem \ref{thm:t1} for $N$-dependent domains. 
Finally, in Section \ref{sec:LP} we give two proofs of Corollary \ref{thm:t0} via the method of logarithmic 
potentials. 
\\
\\
\paragraph{\textbf{Acknowledgments}} The first author was supported by PRC CNRS/RFBR 2017-2019 No.1556 ``Multi-dimensional semi-classical problems of Condensed Matter Physics and Quantum Dynamics''. 
The second author was supported by the Erwin Schr{\"o}dinger Fellowship J4039-N35, by the National Science Foundation grant DMS-1500852 and by CNRS Momentum. We are grateful to the Institut Mittag-Leffler for a stimulating environment. 
\section{A general discussion of Toeplitz band matrices}\label{gendisc}
Let $z\in\C$ and recall \eqref{int.1}. The exponential function $u:\Z \ni \nu \mapsto \zeta^{\nu}$, 
for $\zeta\in \C\backslash\{0\}$, is a solution to 
\begin{equation}\label{g1}
	(p(\tau)-z)u=0,
\end{equation}
if and only if 
\begin{equation}\label{g2}
	p(1/\zeta ) - z = 0. 
\end{equation}
Here, we assume that 
\begin{equation}\label{g2.0}
	z\notin \{0, \infty\}
\end{equation}
Suppose furthermore that 
\begin{equation}\label{g2.1}
	z\neq a_0, \hbox{ when } N_-=0.
\end{equation}
Then \eqref{g2} is equivalent to the following polynomial equation
\begin{equation}\label{g3}
	\sum_{j=0}^{N_++N_-} a_{N_+-j} \zeta^{j} - z\zeta^{N_+} = 0. 
\end{equation}
This is a polynomial equation of degree $N_++N_-$ (when $N_-=0$ 
we have $a_0-z\neq 0$ by \eqref{g2.1}). It has $N_-+ N_+$ roots, counted with their multiplicity.
\par
If $z\notin p(S^1)$, no root is in $S^1$, and we let 
\begin{equation}\label{g4}
	\zeta^+_1, \dots, \zeta^+_{m_+} \quad \text{be the roots in } D(0,1)
\end{equation}
and 
\begin{equation}\label{g5}
	\zeta^-_1, \dots, \zeta^-_{m_-} \quad \text{be the roots in } \C\backslash D(0,1),
\end{equation}
repeated according to their multiplicity. Notice that 
\begin{equation}\label{g6}
	m_+ + m_- = N_+ + N_-.
\end{equation}
\subsection{Remark on exponential solutions}
Let $z\in \C\backslash(\{0\}\cup p(S^1))$. We strengthen assumption \eqref{g2.1} and assume that 
\begin{equation}\label{g9}
	\hbox{if }N_+ ~\text{or} ~ N_- = 0, ~\text{then}~ a_0 \neq  z. 
\end{equation}
Let $\zeta _1,\zeta _2,...,\zeta _m\in \C\setminus \{0 \}$ be the
distinct roots of the characteristic equation (\ref{g2}):
$$
p(1/\zeta )-z=0.
$$
Let $\mathrm{mult\,}(\zeta _j)\ge 1$ be the corresponding multiplicity
so that 
\begin{equation}\label{g6.1}
\sum_1^m\mathrm{mult\,}(\zeta _j)=N_++N_- .
\end{equation}
Similarly to \eqref{g4}, \eqref{g5}, we let 
\begin{equation}\label{g6.2}
	\zeta^+_1, \dots, \zeta^+_{\widetilde{m}_+} \quad \text{be the \emph{distinct} roots in } D(0,1) 
	\text{ with multiplicities } 1\leq\mathrm{mult}(\zeta_j^+) < +\infty,
\end{equation}
and 
\begin{equation}\label{g6.3}
    \zeta^-_1, \dots, \zeta^-_{\widetilde{m}_-} \quad \text{be the  \emph{distinct}  roots in } \C\backslash D(0,1) \text{ with multiplicities } 1\leq\mathrm{mult}(\zeta_j^-) < +\infty, 
\end{equation}
so that $\widetilde{m}_- + \widetilde{m}_+ = m $ in \eqref{g6.1}. Notice also that 
\begin{equation}\label{g6.4}
 	\sum_1^{\widetilde{m}_{\pm}} \mathrm{mult\,}(\zeta _j^{\pm})= m_{\pm}.
\end{equation}
The functions 
$$
\Z\ni \nu \mapsto f_{\zeta ,k}(\nu ):=(\zeta \partial _\zeta )^k(\zeta ^\nu ),\ 0\le k\le
\mathrm{mult\,}(\zeta )-1
$$
are solutions to 
\begin{equation}\label{exp.1}
(p(\tau) -z)f_{\zeta ,k}=0,
\end{equation}
for $\zeta =\zeta _1,...,\zeta _m$. In fact, if $\zeta $ is such a
root, then for $\omega$ close to $\zeta$ 
$$
(p(\tau) -z)(\omega ^\nu )=(p(1/\omega )-z)\omega ^\nu =\mathcal{
  O}((\omega -\zeta )^{\mathrm{mult\,}(\zeta )})
$$
and applying $(\omega\partial _\omega )^k$ with $0\le k\le
\mathrm{mult\,}(\zeta )-1$, and then putting $\omega $ equal to $\zeta
$, we get (\ref{exp.1}).

\par More generally, let $\zeta _1,...,\zeta _m\in \C\setminus \{0 \}$
be distinct numbers and let $1\le m_j<\infty $, $1\le j \le m$.
\begin{prop}\label{exp1}
The functions $f_{\zeta _j,k}\,: \Z\to \C$, $1\le j\le m$, $0\le
k\le m_j-1$ are linearly independent. More precisely, if $K\subset \Z$
is an interval with $\# K=m_1+m_2+...+m_m$, then ${{f_{\zeta
      _j,k}}_\vert}_{K}$ form a basis in $\ell^2(K)$. 
\end{prop}
\begin{proof}
We first prove the linear independence of $f_{\zeta_j,k}$ as 
functions on $\Z$. 
\begin{lemma}\label{lem:exp2}
	Let $\zeta_j$, $j=1,\dots, J$, be finitely many distinct elements of $S^1$. 
	If $a_j\in \C$, $j=1,\dots, J$, and 
	$\lim\limits_{\nu\to +\infty}\sum a_j \zeta_j^{\nu} = 0$, then $a_j=0$. 
\end{lemma}
\begin{proof}[Proof of Lemma \ref{lem:exp2}]
	Write $\zeta_j = \e^{i\sigma_j}$, $\sigma_j\in \R$ and let 
	$\delta_{\sigma_j}\in \mathcal{D}'(S^1)$ be the delta function centered 
	at $\sigma_j$. Then we have 
	\begin{equation*}
		\lim\limits_{\nu\to +\infty} \mathcal{F}^{-1}\left( \sum a_j \delta_{\sigma_j}\right)
		(\nu) = 0
	\end{equation*}
	where $\mathcal{F}^{-1}(u)(\nu) = \frac{1}{2\pi} \int_{S^1} u(x) \e^{ix\nu} dx$. 
	Let $\chi\in \mathcal{C}^{\infty}(S^1)$, $\chi(\sigma_{j_0}) =1$, 
	$\chi(\sigma_j)=0$, $j \neq j_0$. Then 
	\begin{equation*}
	\begin{split}
		\lim\limits_{\nu\to +\infty} \mathcal{F}^{-1}\left(a_{j_0} \delta_{\sigma_{j_0}}\right)
		\!(\nu) &= 
		\lim\limits_{\nu\to +\infty} \mathcal{F}^{-1}\left(\chi \sum
		a_{j} \delta_{\sigma_{j}}\right)
		\!(\nu)\\
		 &= 
		\lim\limits_{\nu\to +\infty} \mathcal{F}^{-1}(\chi ) *\mathcal{F}^{-1}\left(\sum
		a_{j} \delta_{\sigma_{j}}\right)
		\!(\nu)=0, 
	\end{split}
	\end{equation*}
	where $* $ indicates the standard convolution on $\ell^p(\Z)$. 
	Hence, $a_{j_0} =0$ for any $j_0=1,\dots, J$.
\end{proof}
Now consider 
\begin{equation*}
	\sum_{j=1}^m \sum_{k=0}^{m_j-1} a_{j,k}f_{\zeta_j,k} = 0 
	\hbox{ on } \Z 
\end{equation*}
and notice that 
\begin{equation*}
	 f_{\zeta_j,k} = (\zeta \partial_{\zeta})^k(\zeta^{\nu})_{\zeta =\zeta_j} 
	 		   = \nu^k \zeta_j^{\nu}.
\end{equation*}
Let $S=\{j; |\zeta_j| = \max\limits_{\widetilde{j}} |\zeta_{\widetilde{j}}|\}$, $M= \max\limits_{j\in S} m_j$ and write $\zeta_j = \e^{s+i\sigma_j}$. Then we get 
\begin{equation*}
	\lim\limits_{\nu \to \infty} \sum_{\substack{j\in S, \\ m_j = M}} 
	a_{j,M-1}\,\e^{i\nu \sigma_j} = 0.
\end{equation*}
Lemma \ref{lem:exp2} then implies that $a_{j,k}=0$ when 
$|\zeta_j | = \max_{\widetilde{j}} |\zeta_{\widetilde{j}}|$ and 
$k= m_j-1$ is maximal. Repeating this 
procedure we get $a_{j,k} =0$, $0\leq k \leq m_j-1$, $j\in S$. 
Repeating the procedure we finally get $a_{j,k} = 0$ for all $j,k$ and 
we have shown that $f_{\zeta_j,k}$ are independent as functions 
on $\Z$. 
\\
\par
Let 
\begin{equation*}
	Q_{\infty} = \prod_{1}^m (\tau - 1/\zeta_j)^{m_j} = 
	\tau^{\widetilde{m}} + a_1 \tau^{\widetilde{m}-1} + \dots + 
	a_{\widetilde{m}}, \quad 
	\widetilde{m} = m_1 + \dots +m_m.
\end{equation*}
Then as in the case of $p(\tau) -z$, the functions $f_{\zeta_j,k}$, 
$1\leq j \leq m$, $0\leq k\leq m_j -1 $ satisfy 
\begin{equation*}
	Q_{\infty} f_{\zeta_j,k} =0 .
\end{equation*}
Assume that a linear combination $u$ of these functions vanishes 
on the interval $K$ of length $m_1+\dots+ m_m = \widetilde{m}$. Then 
$Q_{\infty}u=0$ on $\Z$, $u=0$ on $K$, and we conclude that $u=0$ on 
$\Z$. Hence $f_{\zeta_j,k}|_K$, $1\leq j\leq m$, $0\leq k \leq m_j-1$ are 
linearly independent. 
\end{proof}
\subsection{Operators on the line and circulant matrices}\label{gendisc_s2} 
Let $S_N\defeq \Z/NZ$, for $N\in \N\backslash \{0\}$. In applications we will replace $N$ 
by $N_++N_-+N$. By convention we set $S_{\infty}=\Z$. 
\\
\par
Recall \eqref{int.0}. We are interested in 
\begin{equation}\label{g17}
	(p(\tau) - z)u=v, \quad u,v\in \ell^2(\Z). 
\end{equation}
Let $\mathcal{F}u(\xi ) = \sum_{k=-\infty}^{\infty} u(k)\e^{-ik\xi }$, so that 
$\mathcal{F}:\ell^2(\Z)\to L^2(S^1,\frac{d\xi }{2\pi})$ is unitary. We have 
\begin{equation}\label{g18}
	(\mathcal{F}\tau u)(\xi ) = \sum_{k=-\infty}^{\infty} u(k-1)\e^{-ik\xi }
						=\e^{-i\xi }(\mathcal{F} u)(\xi ),
\end{equation}
explaining why $e^{-i\xi }=1/\zeta $ is the symbol of $\tau $. 
Hence, application of $\mathcal{F}$ to \eqref{g17} gives the equivalent equation 
\begin{equation}\label{g19}
	(p(\e^{-i\xi }) - z)\widehat{u}=\widehat{v}, \quad 
	\widehat{u}=\mathcal{F}u,\ \ \widehat{v}=\mathcal{F}v.
\end{equation}
Thus, $\spec(p(\tau) ) = p(S^1)$ and if $z\notin p(S^1)$, we can invert \eqref{g19}
\begin{equation}\label{g19.1}
	\widehat{u}(\xi )=\frac{1}{p(\e^{-i\xi }) - z}\widehat{v}(\xi )
\end{equation}
Applying $(\mathcal{F}^{-1}\widehat{u})(k) = \frac{1}{2\pi }
\int_{S^1}\e^{ik\xi }\widehat{u}(\xi )d\xi $, 
we get 
\begin{equation}\label{g20}
	(p(\tau) -z)^{-1}v = K_{\infty}* v,
\end{equation}
where
\begin{equation}\label{g21}
	 K_{\infty}(z;k) = \frac{1}{2\pi }
         \int_{S^1}\frac{1}{p(\e^{-i\xi }) - z}\e^{ik\xi }d\xi .
\end{equation}
In this formula, $S^1$ is identified with $\R/ 2\pi \Z$. Introduce $\zeta = \e^{i\xi}$ as the 
new integration variable, so that $dx = \frac{d\zeta}{i\zeta}$. Then \eqref{g21} becomes 
\begin{equation}\label{g21.1}
	 K_{\infty}(z;k) = \frac{1}{2\pi i} \int_{S^1}\frac{1}{p(1/\zeta) - z}\zeta^k\frac{d\zeta}{\zeta},
\end{equation}
where now $S^1$ is the boundary of the unit disk $D(0,1)\subset \C$. Recall 
\eqref{g6.2}, \eqref{g6.3} and write $m_j^{\pm}=\mathrm{mult}(\zeta_j^{\pm})$. 
If $k \gg 1$, we shrink the contour to $0$ and get by the residue theorem 
\begin{equation}\label{g22.1}
\begin{split}
	 K_{\infty}(z;k) &= 
	 \sum_{j=1}^{\widetilde{m}_+} \lim_{\zeta\to\zeta_j^+} 
	 \frac{1}{(m_j^+-1)!}
	 \frac{d^{m_j^+-1}}{d\zeta^{m_j^+-1}}
	 \frac{\zeta^{k-1}(\zeta-\zeta_j^+)^{m_j^+}}{p(1/\zeta)-z} \\ 
	 & = 
	 \sum_{j=1}^{\widetilde{m}_+} 
	 \sum_{l=1}^{m_j^+}\binom{k-1}{l} \,b^+_{j,l}\,( \zeta_j^+)^{k-l-1}, 
	 \quad b^+_{j,l}\in\C.
\end{split}
\end{equation}
If $k\ll -1$, we use \eqref{g21.1}, enlarge the contour to $|\zeta|=R$, $R\to\infty$, and get
\begin{equation}\label{g23}
\begin{split}
	 K_{\infty}(z;k) &= 
	 -\sum_{j=1}^{\widetilde{m}_-} \lim_{\zeta\to\zeta_j^-} 
	 \frac{1}{(m_j^--1)!}
	 \frac{d^{m_j^--1}}{d\zeta^{m_j^--1}}
	 \frac{\zeta^{k-1}(\zeta-\zeta_j^-)^{m_j^-}}{p(1/\zeta)-z} \\ 
	 & = 
	 -\sum_{j=1}^{\widetilde{m}_-} 
	 \sum_{l=1}^{m_j^-}\binom{k-1}{l} \,b^-_{j,l}\,( \zeta_j^-)^{k-l-1}, 
	 \quad b^-_{j,l}\in\C.
\end{split}
\end{equation}
\begin{remark}\label{remSR}
	When all roots of the polynomial \eqref{g3} are simple, then we have by 
	\eqref{g4},\eqref{g5}, \eqref{g6.2}, \eqref{g6.2} as well as \eqref{g22.1}, \eqref{g23} 
	that 
	\begin{equation}\label{g23a}
	 K_{\infty}(z;k) = 
	 \begin{cases}
	 \sum_{j=1}^{m_+}\frac{1}{\partial _\zeta ( p(1/\zeta))_{\zeta=\zeta^+_j}}(\zeta^+_j)^{k-1}, 
	 ~ \text{ if } k\geq 1, \\
	 - \sum_{j=1}^{m_-}\frac{1}{\partial _\zeta ( p(1/\zeta))_{\zeta=\zeta^-_j}}(\zeta^-_j)^{k-1}, 
	  ~ \text{ if } k\leq -1.
	 \end{cases}
\end{equation}
\end{remark}
Notice that $K_{\infty}(z;k)$ decays exponentially as $|k|\to\infty$. Hence, we can solve \eqref{g17} 
for $u,v\in \ell^{\infty}$. 
\par
If $v\in \ell^2(S_N)$, then we can view $v$ as an $N$-periodic function 
on $\Z$ and the solution $u$ is $N$-periodic and given by \eqref{g20}. 
\par
Let $\Omega\subset\Z$ be a finite set of cardinal $\# \Omega = N$ such that 
\begin{equation*}
	\begin{cases}
		(\Omega + jN) \cap (\Omega + kN) = \emptyset \text{ for } j\neq k \\
		\bigcup\limits_{j\in\Z}\Omega+jN = \Z.
 	\end{cases}
\end{equation*}
Let $N\geq N_+ + N_- +1$. Still when $u,v$ are $N$-periodic we make \eqref{g20} more explicit 
\begin{equation}\label{g24}
	\begin{split}
	 u(\nu )& = \sum_{\mu \in \Z} K_{\infty}(z;\nu-\mu)v(\mu) 
	 	    = \sum_{j \in \Z}\sum_{\mu \in \Omega + jN} K_{\infty}(z;\nu-\mu)v(\mu) \\
		   & = \sum_{j \in \Z}\sum_{\mu \in \Omega} K_{\infty}(z;\nu-\mu-jN)v(\mu)
		    = \sum_{\mu \in \Omega} K_{N}(z;\nu-\mu)v(\mu),
        \end{split}
\end{equation}
where 
\begin{equation}\label{g25}
	K_{N}(z;\nu-\mu) = \sum_{j \in \Z}K_{\infty}(z;\nu-\mu-jN)
\end{equation}
and the series converges geometrically. We check that $K_{N}(z;\nu+N) = K_{N}(z;\nu)$. 
Identifying $\Omega \simeq S_N$, and defining 
\begin{equation}\label{g25.1}
	P_{S_N} \defeq p(\tau) : {\ell^2(S_N) \to\ell^2(S_N)},
\end{equation}
we get 
\begin{prop}\label{g_prop4} 
  If $z\notin p(S^1)$, then $z\notin \spec(P_{S_N})$ and the resolvent
  $(P_{S_N}-z)^{-1}$ is given by
	\begin{equation}\label{g26}
          (P_{S_N}-z)^{-1}v(\nu) = (K_N(z) * v)(\nu) = \sum_{\mu\in S_N}K_{N}(z;\nu-\mu)v(\mu)
	\end{equation}
	with
	\begin{equation}\label{g27}
	K_{N}(z;\nu) = \sum_{\widetilde{\nu}\in \pi^{-1}(\nu)}K_{\infty}(z;\widetilde{\nu})
	\end{equation}
	where $\pi: \Z \to S_N$ is the natural projection. 
\end{prop}
A consequence of \eqref{g25} is the following: Choose $\Omega=[-\frac{N}{2},\frac{N}{2}[$ when 
$N$ is even and $\Omega=[-\frac{N-1}{2},\frac{N+1}{2}[$ when $N$ is odd. Then, 
\begin{equation}\label{g28}
	K_{N}(z;\nu) = K_{\infty}(z;\nu) + \mO\left( \e^{-\frac{N}{C}} \right), \quad \nu \in \Omega. 
\end{equation}
\subsection{The spectrum of $P_{S_N}$}\label{sec:SpG}
Using the finite Fourier transform $\ell^2(S_N) \to\ell^2(\widehat{S}_N)$, 
with $\widehat{S}_N=\{\e^{\frac{2\pi ik}{N}};k=0,\dots, N-1\}$, it is easy to 
prove that 
\begin{equation}\label{g27.5}
\spec(P_{S_N})=p(\widehat{S}_N).
\end{equation}
In this section we study the spectrum of the normal operator $P_{S_N}$, see \eqref{g25.1} and 
\eqref{gu7} below, in 
\begin{equation}\label{sp3}
	\gamma = p(S^1)\cap \Omega 
\end{equation}
with $\Omega$ as in Section \ref{sec:2.R2}. 
\subsubsection{A Weyl law for $P_{S_N}$}
We present a Weyl law for the eigenvalues of $P_{S_N}$, 
which we shall use later on to count the eigenvalues of small perturbations 
of the operator $P_N$ \eqref{int.2}.
\par
Let $\gamma$ be as in \eqref{sp3}. First notice that by
 \eqref{g27.5} 
\begin{equation}\label{sp4}
	\#\{ \spec(P_{S_N}) \cap \gamma \} = 
	\#\{\widehat{S}_{N} \cap p^{-1}(\gamma) \}.
\end{equation}
Since two consecutive points of $\widehat{S}_{N}$ differ by an angle 
of $2\pi/N$, we get that
 \begin{equation}\label{sp5}
	\#\{\widehat{S}_{N} \cap p^{-1}(\gamma) \} 
	=\frac{N}{2\pi} \int_{p^{-1}(\gamma)}L_{S^1}(d\theta) +\mO(1),
\end{equation}
where the measure $L_{S^1}(d\theta)$ in the integral denotes the Lebesgue 
measure on $S^{1}$. Combining \eqref{sp4}, \eqref{sp5}, we get 
\begin{equation}\label{sp6}
	\#\{ \spec(P_{S_N}) \cap \gamma \}
	=\frac{N}{2\pi} \int_{p^{-1}(\gamma)}L_{S^1}(d\theta)+\mO(1).
\end{equation}
\subsubsection{Local eigenvalue spacing for $P_{S_N}$}
Let $z_0 \in p(S^1)$ be such that 
\begin{equation}\label{sp2.0}
	dp \neq 0  \text{ on } p^{-1}(z_0).
\end{equation}
\begin{prop}\label{prop:EigDist}
Let $p$ be as in \eqref{int.0.1} and let $z_0\in p(S^1)$ be such that \eqref{sp2.0} holds.
Then, there exist a constant $C>0$ and an open neighborhood $U\subset \C$ of $z_0$, 
such that $p^{-1}(U)$ is the union of finitely many disjoint open sets $V_i\subset \C$, 
$i=1,\dots, M$. Moreover, on each non-empty segment $\gamma_i = p(V_i\cap S^1)$ we have that 
\begin{equation}\label{spn0}
	\min\limits_{\substack{z,w\in p(\widehat{S}_{\widetilde{N}})\cap \gamma_i  \\ w\neq z}}
	|z -w|\geq  \frac{1}{CN}.
\end{equation}
\end{prop}
\begin{proof}[Proof of Proposition \ref{prop:EigDist}]
For $i=1,\dots, M$ let $\zeta_i \in p^{-1}(z_0)$ and 
notice that $M< +\infty$. By \eqref{sp2.0} and the implicit function theorem, 
there exist complex open neighborhoods $U_i$ of $z_0$ and $V_i$ of $\zeta_i$ 
such that $p: V_i \to U_i$ is a diffeomorphism. Setting $\gamma_i = p(S^1\cap V_i)\subset U_i$  when
$S^1\cap V_i \neq \emptyset$ , we have that
\begin{equation}\label{spn1}
	| \widehat{\zeta}^i_n - \widehat{\zeta}^i_m| \asymp |\widehat{z}^i_n -\widehat{z}^i_m |, 
\end{equation}
where $ \widehat{\zeta}^i_n \in\widehat{S}_{\widetilde{N}}\cap V_i $ and 
$\widehat{z}^i_n \in p(\widehat{S}_{\widetilde{N}})\cap \gamma_i $, for $n \in J_i \subset \N$, some 
index set which is non-empty for $N>1$ sufficiently large. Since $M$ is finite, the 
claim follows by \eqref{g27.5} and by taking $U= \bigcap_{i=1}^M U_i$ and by potentially shrinking the 
segements $\gamma_i$. 
\end{proof}
\subsection{Restrictions to intervals}\label{gendisc_s1}
If $K\subset \Z$ is a finite set or an infinite interval, we identify 
\begin{equation}\label{g7}
	\ell^2(K) \simeq \ell^2_K \defeq \{ u\in \ell^2(\Z); ~\supp u \subset K \}.
\end{equation}
We define, 
\begin{equation}\label{g8}
	P_K \defeq 1_K\, p(\tau): ~ \ell^2_K \longrightarrow \ell^2_K, \quad 
	\text{and } \quad P_{\Z} = p(\tau). 
\end{equation}
In the following we assume \eqref{g9}. When $K$ is 
an interval we define the length of $K$ to be $\# K = |K|$. 
\begin{prop}\label{g_prop1}
	Let $K$ be an interval of length $\leq N_+ +N_-$. Any function 
	$u:~K\to \C$ can be extended to a solution $\widetilde{u}:~\Z \to \C$ 
	to $(p(\tau)-z)u=0$. The space of such extensions is affine of 
	dimension $N_++N_- - \# K$. In particular the extension is unique 	
	when $N_++N_- =\# K$. 
\end{prop}
\begin{proof}
  If $\# K < N_+ + N_-$, let $\widetilde{K}\supset K$ be an interval
  with $\# \widetilde{K} = N_++N_-$. The extensions
  $\hat{u}:~\widetilde{K} \to C$ form an affine space of dimension
  $N_++N_--K$, so it suffices to treat the case $\#K =N_+ + N_-$.
	\\
	\par
	Let $K = [M,M+N_++N_-[$ and write $(p(\tau)-z)\widetilde{u}=0$,  i.e. 
	\begin{equation}\label{g10}
		a_{N_+}\widetilde{u}(\nu - N_+) + 
		\dots +(a_{0}-z)\widetilde{u}(\nu) + \dots + 
		a_{-N_-}\widetilde{u}(\nu + N_-) = 0.
	\end{equation}
	\par
	For $\nu = N_+ + M$, $\nu+N_-$ is the first point in $\Z\backslash K$ to 
	the right of $K$ and $\nu +N_--1, \dots ,\nu - N_+$ belong to $K$, so 
	\eqref{g10} defines $\widetilde{u}(\nu+N_-)$ uniquely. 
	Replacing $\nu$ with $\nu+1=M+N_++1$, we get $\widetilde{u}(M+N_++N_-+1)$ 
	and by repeating the procedure we get $\widetilde{u}(M+N_++N_-+\mu)$ for all $\mu\geq 0$. 
	\par
	For $\nu=M+N_+-1$, we have $\nu-N_+\notin K$ while $\nu-N_++1,\dots,\nu+N_-\in K$, 
	and therefore \eqref{g10} determines $\widetilde{u}(M-1)$ uniquely. Iterating the procedure, 
	we get all values of $\widetilde{u}(M-\mu)$, for $\mu > 0$.
\end{proof}
It follows that the space of solutions to $(p(\tau)-z)u=0$ is of dimension $N_++N_-=m_++m_-$, 
cf. \eqref{g6}. Recall \eqref{g6.2}, \eqref{g6.3}, \eqref{g6}, \eqref{g9}, \eqref{g2.0}  and \eqref{g6.1}. The space of exponential solutions, 
spanned by the functions 
\begin{equation}\label{g11}
 \Z \ni \nu \mapsto \nu ^{k} (\zeta _{j_\pm}^\pm)^\nu, \quad  \text{for } 1\le j_{\pm}\le \widetilde{m}_\pm, ~ 
 	0\leq k \leq \mathrm{mult}(\zeta _{j_\pm}^\pm)-1,
\end{equation}
is also of dimension $m_++m_-$, since these functions form a linearly independent system by 
Proposition \ref{exp1}. Hence, assuming \eqref{g9}, \eqref{g2.0}, they form a basis of the space of solutions 
$u: \Z \to \C$ to $(p(\tau)-z)u=0$. We conclude the following
\begin{prop}\label{g_prop2}
	Suppose \eqref{g9} and \eqref{g2.0}. Then, the general solution $u:\Z\to \C$ of $(p(\tau)-z)u=0$ is of the form 
	\begin{equation}\label{g12}
		u(\nu) = \sum_{j=1}^{\widetilde{m}_+}\sum_{k=0}^{\mathrm{mult}(\zeta_j^+)-1}a_{j,k}^+\nu^k(\zeta_j^+)^\nu 
		+ \sum_{j=1}^{\widetilde{m}_-}\sum_{k=0}^{\mathrm{mult}(\zeta_j^-)-1}a_{j,k}^-\nu^k(\zeta_j^-)^\nu, 
		\quad a^{\pm}_{j,k} \in \C. 
	\end{equation}
	The subspace of solutions decaying at $\nu\to \pm\infty$ is given by 
	\begin{equation}\label{g12a}
		a_{j,k}^{\mp}  = 0, \quad \text{for } j = 1,\dots, \widetilde{m}_\mp, ~ k =0, \dots, \mathrm{mult}(\zeta_j^+)-1.
	\end{equation}
\end{prop}
\begin{remark}
	Enumerate all the roots of \eqref{g3} as 
	\begin{equation*}
	w_j = \begin{cases}
		\zeta_j^+, \text{ for } j=1,\dots,\widetilde{m}_+ \\
		\zeta_{j-\widetilde{m}_+}^-, \text{ for } j=\widetilde{m}_++1,\dots,\widetilde{m}_+ +\widetilde{m}_-
	\end{cases}
	\end{equation*}
	so that \eqref{g12} takes the form 
	\begin{equation}\label{g13}
		u(\nu) = \sum_{j=1}^{\widetilde{m}_++\widetilde{m}_-}\sum_{k=0}^{\mathrm{mult}(w_j)-1}b_{j,k}\nu^k(w_j)^\nu, 
		\quad b_{j,k}\in \C.
	\end{equation}
	We then recover the fact that the following Van der Monde type determinant 
	\begin{equation}\label{g13a}
		\det(A_1,\dots,A_{\widetilde{m}_++\widetilde{m}_-})
	\end{equation}
	is non-vanishing. Here, the block matrices $A_{j}$, $j=1,\dots,\widetilde{m}_++\widetilde{m}_-$, are given by 
	\begin{equation}\label{g13b}
		A_j = ( \nu^k w_j^{\nu} )_{\substack{\nu\in K \\ 0 \leq k \leq \mathrm{mult}(w_j)-1}} \in \C^{|K|\times \mathrm{mult}(w_j)}
	\end{equation}
	where $K$ is any interval of length $|K|=m_++m_-$. 
\end{remark}
We next look at $P_K$ where $K$ is the half-axis $[A,+\infty[$ or $]-\infty,A]$. The two 
cases are similar and we may assume by translation invariance that
$K=[0,+\infty [$. 
\par
Let $u:K\to\C$ have its support in $[0,\infty[$ and satisfy 
\begin{equation}\label{g14}
	(p(\tau)-z)u=0 ~~ \text{in} ~~ [0,+\infty[.
\end{equation}
More explicitly, by \eqref{int.0}, 
\begin{equation}\label{g16}
	\left(\sum_{j=-N_-}^{N_+} a_j \tau^j-z\right)u(\nu)=0,\ \ \nu =0,1,\dots
\end{equation}
The left most equation for $\nu =0$ is 
\begin{equation*}
	a_{N_+}u(-N_+) + a_{N_+-1}u(1-N_+) + \dots 
	+ (a_0-z)u(0) + \dots + a_{-N_-}u(N_-) =0.
\end{equation*}
Here, $u(-N_+)=\dots=u(-1)=0$, when $N_+\leq 1$. We know how to extend 
$u|_{[-N_+,+\infty[}$ to a function $\widetilde{u}:\Z\to\C$, by solving \eqref{g16} with 
$u$ replaced by $\widetilde{u}$ for $\nu =-1,-2,\dots$. The equation for $\nu=-1$ 
defines $\widetilde{u}(-N_+-1)$, the next one gives $\widetilde{u}(-N_+-2)$ and so 
on. In this way we get a solution $\widetilde{u}$ on $\Z$ of 
\begin{equation}\label{g14.1}
	(p(\tau)-z)\widetilde{u}=0.
\end{equation}
Consequently $\widetilde{u}$ has the form of the right hand side in \eqref{g12}. Now restrict 
the attention to solutions $u\in\ell^2_{[0,+\infty[}(\Z)$ of \eqref{g14}. The corresponding extension 
$\widetilde{u}$ is of the form \eqref{g12} with $a_{j,k}^-=0$, since it must decay to the right. Hence, 
\begin{equation}\label{g15.1}
	\widetilde{u}(\nu) = 
	\sum_{j=1}^{\widetilde{m}_+}\sum_{k=0}^{\mathrm{mult}(\zeta_j^+)-1}a_{j,k}^+\nu^k(\zeta_j^+)^\nu
\end{equation}
and by construction $\widetilde{u}(\nu)=u(\nu)=0$ for $\nu\in [-N_+,-1]$. More 
explicitly, using \eqref{g6.4}, we have 
\begin{equation}\label{g16.1}
\begin{split}
	&0=A \begin{pmatrix} a_{1,0}^+ \\ \vdots \\ a_{\widetilde{m}_+,\mathrm{mult}(\zeta_{\widetilde{m}_+})-1}^+ \end{pmatrix},  
	\quad A=(A_1^+,\dots,A_{\widetilde{m}_+}^+) \in \C^{N_+\times m_+},
	\\
	& A_j^+ = ( \nu^0(\zeta_j^+)^\nu,\dots, \nu^{\mathrm{mult}(\zeta_j^+)-1}(\zeta_j^+)^\nu)_{-N_+ \leq \nu\leq -1}, \quad 
	\text{for } j=1,\dots, \widetilde{m}_+.
\end{split}
\end{equation}
Notice that $A$ is a rectangular generalized matrix of Van der Monde type, 
of size $N_+\times m_+$. Arguing as at the end of the proof of Proposition 
\ref{exp1} and using \eqref{g6.4}, we see that $A$ is of maximal rank 
$\min (N_+,m_+)$. Thus
\begin{itemize}
	\item if $N_+ \geq m_+$, then 
		\begin{equation*}
		\ker\left(P _{[0,+\infty[}-z\right) =0.
		\end{equation*}
	\item If $N_+ < m_+$, then 
		\begin{equation*}
			\dim \ker\left(P _{[0,+\infty[}-z\right)  = m_+ - N_+.
		\end{equation*}
	\end{itemize}
For $(P_{]-\infty,0]}-z)$ we have the corresponding statements with $N_+,m_+$ 
replaced by $N_-,m_-$.

\begin{lemma}\label{lem:ToOpFred} 
Let $z\notin p(S^1)$, then the operators $(P _{[0,+\infty[}-z): \ell^2([0,+\infty[) \to  \ell^2([0,+\infty[) $ and 
$(P _{]-\infty,0]}-z): \ell^2(]-\infty,0]) \to  \ell^2(]-\infty,0]) $ are Fredholm. 
\end{lemma}
\begin{proof}
	We give the proof for $(P _{[0,+\infty[}-z)$, the one for $(P _{]-\infty,0]}-z) $ 
	is similar. 
	\par
	Recall \eqref{g20}, and define for $z\notin p(S^1)$
	\begin{equation*}
		E(z) = \mathbf{1}_{[0,+\infty[} ( p(\tau) - z )^{-1}\mathbf{1}_{[0,+\infty[}.
	\end{equation*}
	Then,	
	\begin{equation*}
		( P _{[0,+\infty[} - z )E(z) = \mathbf{1}_{[0,+\infty[} + R(z)
	\end{equation*}
	and
	\begin{equation*}
		E(z) ( P _{[0,+\infty[} - z )= \mathbf{1}_{[0,+\infty[} + L(z).
	\end{equation*}
	where $R(z),L(z)$ are compact. Indeed, we have 
	\begin{equation*}
	R(z) 
	= - \mathbf{1}_{[0,+\infty[} ( p(\tau) - z )\mathbf{1}_{]-\infty,0[} 
		( p(\tau) - z )^{-1}\mathbf{1}_{[0,+\infty[}.
	\end{equation*}
	By \eqref{int.0} we see that $R(z) = 1_{[0,N_+[}R(z)$, so 
    $R(z)$ is of finite rank and thus compact. Similarly, we have  
	\begin{equation*}
	L(z) 
	= - \mathbf{1}_{[0,+\infty[} ( p(\tau) - z )^{-1}\mathbf{1}_{]-\infty,0[} 
		( p(\tau) - z )\mathbf{1}_{[0,+\infty[}.
	\end{equation*}
	We notice that $L(z) = L(z)1_{[0,N_+[}$ is of finite rank, hence compact.
\end{proof}
Next, notice that by \eqref{g2}, \eqref{g4}, \eqref{g5}, $p(\tau)^*$ is similar 
to $p(\tau)$ just with the roles of $N_+,m_+$ and $N_-,m_-$ exchanged. 
More explicitly, by \eqref{int.0}, 
\begin{equation*}
	p^*(\tau) = \overline{p}(\tau^{-1}) 
	= \sum_{-N_-}^{N_+} \overline{a}_{j} \tau^{-j}
	= \sum_{-N_+}^{N_-} \overline{a}_{-j} \tau^{j}.
\end{equation*}
The analogue of \eqref{g2} is $\overline{p}(\omega) - \overline{z} =0$, 
or equivalently $p(\overline{\omega}) - z =0$, since 
$\overline{p}(\omega) = \overline{p(\overline{\omega})}$. 
In view of \eqref{g4}, \eqref{g5}, 
we get the roots $\omega_j^{\pm} = 1/\overline{\zeta_j^{\pm}}$.
Remembering \eqref{g8}, we have 
\begin{equation*}
	P_K^* = 1_K \overline{p}(\tau^{-1})1_K
\end{equation*}
Therefore, the above statements remain valid with $(p(\tau)-z)$ 
replaced by $(p^*(\tau)-\overline{z})$  and $N_+,m_+$ exchanged with 
$N_-,m_-$. 
\\
\par
By Lemma \ref{lem:ToOpFred} we get that for $z\notin p(S^1)$ 
\begin{equation*}
	\dim \ker\left(P _{[0,+\infty[}-z\right)^* = \dim \coker\left(P _{[0,+\infty[}-z\right)
\end{equation*}
Hence, using \eqref{g6} we conclude the following 
\begin{prop}\label{g_prop3}
	Assume that $z\notin \{0,+\infty\}\cup p(S^1)$ and that \eqref{g9} holds.
	\begin{itemize}
	\item If $N_+ \geq m_+$, then 
		\begin{equation*}
		\ker\left(P _{[0,+\infty[}-z\right) =0
		\end{equation*}
		and
		\begin{equation*}
		\dim \coker\left(P _{[0,+\infty[}-z\right)  = N_+-m_+.
		\end{equation*}
	\item If $N_+ < m_+$, then 
		\begin{equation*}
			\coker\left(P _{[0,+\infty[}-z\right) =0
		\end{equation*}
		and
		\begin{equation*}
			\dim \ker\left(P _{[0,+\infty[}-z\right)  = m_+ - N_+.
		\end{equation*}
	\end{itemize}
	For $(P_{]-\infty,0]}-z)$ we have the corresponding statements with $N_+,m_+$ replaced by 
	$N_-,m_-$.
\end{prop}
It will be convenient to replace $P_{[0,+\infty[}$ with the unitarily equivalent operator 
$P_{[N_+,+\infty[}$. Moreover, let us recall that the \emph{index} of a Fredholm 
operator $A$ is defined by 
\begin{equation*}
		\mathrm{Ind} \, A \defeq	\dim \ker A -\dim \coker\ A.
\end{equation*}
There is a very nice relation between the index of the Fredholm operator 
$(P_{[N_+,+\infty[}-z) $ and the winding number of the curve $p(S^1)$ around the 
point $z$. 
\begin{prop}\label{g_prop4.1}
	Let $z\notin\{0,+\infty\}\cup p(S^1)$ and suppose that \eqref{g9} holds. 
	Then $(P_{[N_+,+\infty[}-z)$ is Fredholm of index
	\begin{equation}
		\mathrm{Ind}(P_{[N_+,+\infty[}-z) = m_+-N_+ = -\mathrm{ind}_{p(S^1)}(z). 
	\end{equation}
\end{prop}
\begin{proof}
	The first equality follows from Proposition \ref{g_prop3} (see also  
	\eqref{g17.1}, \eqref{g20.1}). To see the second equality, notice that 
	\begin{equation}
		-\mathrm{ind}_{p(S^1)}(z)=
		\frac{1}{2\pi i} \int_{S^1}\frac{d}{d\eta} \log (p(1/\eta) -z)d\eta.
	\end{equation}
	The integral on the right hand side is equal to the number of zeros minus the 
	number of poles of 
	$p(1/\eta)-z$ in $D(0,1)$, where both are counted including multiplicity. This is equal 
	to $m_+-N_+$ by \eqref{g3}, \eqref{g4} and \eqref{g5}.
\end{proof}
\begin{remark}
This result has been obtained by M.G. Krein via a different method. See 
\cite[Chapter 1.5]{BoSi99} for a detailed exposition. 
\end{remark}
\subsection{Zone of zero winding number}
In this section we show that in regions in $\C$, for which the winding number of the curve 
$p(S^1)$ is zero, the norm of the resolvent of $P_N$ is controlled by a constant. Hence, 
we can consider such regions to ``spectrally stable" for $P_N$. 
\begin{prop}\label{resUB}
 Let $\Omega\Subset \C\backslash( \{0\}\cup p(S^1))$ be a compact set and  
suppose that for every $z\in\Omega$ \eqref{g9} holds and 
\begin{equation}\label{st1}
		\mathrm{ind}_{p(S^1)}(z)= 0.
\end{equation}
 Then, there exists a constant $C>0$ such that for $N>0$ sufficiently large and 
 for any $z\in\Omega$
 \begin{equation*}
		\|(P_N-z)^{-1}\| \leq C.
\end{equation*}
\end{prop}
\begin{proof}
By Propositions \ref{g_prop4.1}, \ref{g_prop3} and by \eqref{st1}, we know that 
$(P_{[1,+\infty[}-z)$ and $(P_{]-\infty,N]}-z)$ are bijective on $\ell^2$ with uniformly 
bounded inverses when $z\in\Omega$. By the Combes-Thomas argument the same 
holds after conjugation with a factor $\e^{\varepsilon \varphi}$ if $\varphi$ is Lipschitz 
of modulus $\leq 1$ and $|\varepsilon|$ is small enough. 
\par
Let 
 \begin{equation*}
	Q_N(z) = \mathbf{1}_{[1,N]}( (P_{[1,+\infty[}-z)^{-1}\mathbf{1}_{[1,[N/2]]} 
	+ (P_{]-\infty,N]}-z)^{-1}\mathbf{1}_{[[N/2]+1,N]}  ).
 \end{equation*}
 Then, using the stability under exponential conjugation, it follows that 
  \begin{equation*}
	(P_{[1,N]}-z)Q_N(z) = 1 + R, \quad \| R\| \leq \mO(1) \e^{-N/C}.
 \end{equation*}
 Hence, for $N>1$ large enough, $P_{[1,N]} : \ell^2([1,N]) \to \ell^2([1,N])$ has a uniformly bounded right 
 inverse which is also a left inverse since $P_{[1,N]}$ is a finite square matrix. 
\end{proof}
\section{A Grushin Problem}\label{grush}
We begin by giving a short refresher on Grushin problems. See 
\cite{SjZw07} for a review. The central idea is to set up an auxiliary 
problem of the form 
\begin{equation*}
 \begin{pmatrix}
  P(z) & R_- \\ 
  R_+ & R_{+-} \\
 \end{pmatrix}
 :
 \mathcal{H}_1\oplus \mathcal{H}_- 
 \longrightarrow \mathcal{H}_2\oplus \mathcal{H}_+,
\end{equation*}
where $P(z)$ is the operator under investigation and $R_{\pm}, R_{+-}$ are 
suitably chosen. We say that the Grushin problem is well-posed 
if this matrix of operators is bijective. If 
$\dim\mathcal{H}_-  = \dim\mathcal{H}_+ < \infty$, one typically 
writes 
\begin{equation*}
 \begin{pmatrix}
  P(z) & R_- \\ 
  R_+ & R_{+-}\\
 \end{pmatrix}^{-1}
 =
 \begin{pmatrix}
  E(z) & E_{+}(z) \\ 
  E_{-}(z) & E_{-+}(z) \\
 \end{pmatrix}.
\end{equation*}
The key observation goes back to the Shur complement formula or, 
equivalently, the Lyapunov-Schmidt bifurcation method, 
i.e. the operator $P(z): \mathcal{H}_1 \rightarrow \mathcal{H}_2$ 
is invertible if and only if the finite dimensional matrix 
$E_{-+}(z)$ is invertible and when $E_{-+}(z)$ is invertible, 
we have 
\begin{equation*}
  P^{-1}(z) = E(z) - E_{+}(z) E_{-+}^{-1}(z) E_{-}(z).
\end{equation*}
$E_{-+}(z)$ is sometimes called effective Hamiltonian. 
\subsection{A Grushin problem for the unperturbed operator}\label{grush_upo}
Let $J\subset \Z$ be a fixed interval of length $\# J = N_++N_-$. More precisely, we 
choose
\begin{equation}\label{gu1}
	J=[-N_-,N_+[. 
\end{equation}
If $M > N_+ + N_-$ we view $J$ as a segment of $S_M$, cf. the beginning of Section 
\ref{gendisc_s2}. More precisely we define a segment $[a,b]\subset S_M$, $a,b\in S_M$, 
to be the set of points in $S_M$ that we get by picking first $a$, then $a+1$ and so on until 
we reach $b$ (mod $M\Z$) with the last point $b$ included. Similarly we define $[a,b[$, 
$]a,b[$,  $]a,b]$. Recall that $S_{\infty}=\Z$. 
\par
Suppose that
\begin{equation}\label{gu2.0}
	N\geq N_++N_-+1.
\end{equation}
When $N$ is finite we decompose 
\begin{equation}\label{gu2}
	S_{N+N_++N_-} = J \cup I_N,
\end{equation}
\begin{equation}\label{gu3}
	I_N = [N_+,-N_--1] 	
\end{equation}
where $I_N \simeq [N_+,-N_--1+N_++N_-+N]= [N_+,N_++N-1] \text{ in } \Z$. 
When $N=\infty$, we decompose
\begin{equation}\label{gu4}
	\Z=S_{\infty} = J \cup I_{\infty},
\end{equation}
\begin{equation}\label{gu5}
	I_{\infty} = ]-\infty, -N_--1]\cup [N_+,\infty[.
\end{equation}
Since $\# I_N= N$, we can identify 
\begin{equation}\label{gu5.1}
	P_N\simeq P_{I_N},
\end{equation}
in view of \eqref{int.2}, when $N$ is finite, while $P_{I_{\infty}}$ is the direct sum 
\begin{equation}\label{gu6}
	P_{]-\infty,-N_--1]}\oplus P_{[N_+,\infty[} \simeq P_{]-\infty,0]}\oplus P_{[0,\infty[}.
\end{equation}
In both cases we identify 
\begin{equation*}
	\ell^2(S_{N+N_-+N_-}) \simeq \ell^2(I_N) \oplus \ell^2(J)
\end{equation*}
so that 
\begin{equation}\label{gu7}
	(P_{S_{N+N_-+N_-}}-z) \defeq \mathcal{P}_N(z) = 
	\begin{pmatrix} 
		P_{I_N} -z & R_-^N \\
		R_+^N & R_{+-}^N(z) 
	\end{pmatrix}
	: \ell^2(I_N) \oplus \ell^2(J) \to \ell^2(I_N) \oplus \ell^2(J)
\end{equation}
where 
\begin{equation}\label{gu8}
	\begin{split}
		&P_{I_N}-z = 1_{I_N} (p(\tau) -z)1_{I_N}, \quad R_-^N= 1_{I_N} p(\tau) 1_{J}, \\ 
		&R_+^N = 1_{J} p(\tau) 1_{I_N}, \quad R_{+-}^N(z)= 1_{J} (p(\tau)-z) 1_{J}. \\ 
	\end{split}
\end{equation}
\begin{lemma}\label{gu_lem1}
	$R_+^N$ is surjective and $R_-^N$ is injective. 
\end{lemma}
\begin{proof}
	Suppose that $\supp u \subset [-N_--N_+,-N_-[ \subset I_N$. Then, 
	$\supp R_+^N u \subset  [-N_-,-N_-+N_+[$. By fixing the values of 
	$u(-N_- -1),\dots,u(-N_+-N_-)$ we can arrange so that $R_+^N u $ is equal to 
	any given function with support in 
	$[-N_-,-N_-+N_+[$. 
	\par
	Similarly, if $\supp u \subset [N_+,N_++N_-[$ then 
	$\supp R_+^N u \subset  [-N_+-N_-,N_+[$ and a convenient choice of such 
	a $u$ will produce any given function with support in $[N_+-N_-,N_+[$. Since 
	$J=[-N_-,-N_-+N_+[\cup [N_+-N_-,N_+[$ and $[-N_- -N_+, -N_-[$, $[N_+,N_++N_-[$ 
	are by \eqref{gu2.0} disjoint subsets of $I_N$, we see that $R_+^N$ is surjective.
	\par
	For the same reason $^t(R_-^N)=1_J{^t p(\tau)}1_{I_N}$ is surjective and 
	therefore $R_-^N$ is injective. 
\end{proof}
Recall \eqref{g27.5}. If $z\notin\spec(P_{S_{\widetilde{N}}})$, where $\widetilde{N}=N+N_-+N_-$, then  
$\mathcal{P}_N(z)$ in \eqref{gu7} is bijective and invertible with bounded inverse
\begin{equation}\label{gu9}
	\mathcal{E}_N(z) = 
	\begin{pmatrix} 
		E^N(z) & E_+^N(z) \\
		E_-^N(z) & E_{-+}^N(z) 
	\end{pmatrix}
	: \ell^2(I_N) \oplus \ell^2(J) \to \ell^2(I_N) \oplus \ell^2(J).
\end{equation}
We have 
\begin{equation}\label{gu10}
	\begin{split}
		E^N(z) = 1_{I_N} (P_{S_{\widetilde{N}}}-z)^{-1}1_{I_N}, 
		\quad E_+^N(z)= 1_{I_N} (P_{S_{\widetilde{N}}}-z)^{-1} 1_{J}, \\ 
		E_-^N(z) = 1_{J} (P_{S_{\widetilde{N}}}-z)^{-1} 1_{I_N}, 
		\quad  E_{-+}^N(z) = 1_{J} (P_{S_{\widetilde{N}}}-z)^{-1} 1_{J}. \\ 
	\end{split}
\end{equation}
If $z\notin p(S^1)$, then this also holds for $N=\infty$. We now recall 
Proposition \ref{g_prop4} and \eqref{g25} 
with $N$ replaced by $\widetilde{N}$. On the level of matrices we get with 
$\pi=\pi_{\widetilde{N}}: \Z \to S_{\widetilde{N}}$
\begin{equation}\label{gu11}
	\begin{split}
		&E^N(z;\nu,\mu) = \sum_{\widetilde{\nu}\in\pi^{-1}(\nu)}
				          E^{\infty}(z;\widetilde{\nu},\widetilde{\mu}), 
				          \quad \widetilde{\mu}\in\pi^{-1}(\mu), ~\mu,\nu\in I_N, \\
		&E^N_+(z;\nu,\mu) = \sum_{\widetilde{\nu}\in\pi^{-1}(\nu)}
				          E_+^{\infty}(z;\widetilde{\nu},\mu), 
				          \quad \mu\in J, ~\nu\in I_N, \\
		&E^N_-(z;\nu,\mu) = \sum_{\widetilde{\mu}\in\pi^{-1}(\mu)}
				          E_-^{\infty}(z;\nu,\widetilde{\mu}), 
				          \quad \nu\in J, ~\mu\in I_N,
	\end{split}
\end{equation}
and
\begin{equation}\label{gu12}
		E^N_{-+}(z;\nu,\mu) = \sum_{j\in\Z} (p(\tau)-z)^{-1}(\nu+j\widetilde{N},\mu)
				          \quad \mu,\nu\in J. 
\end{equation}
In these formulas we used that $J$ is naturally defined both as a subset of $S_{\widetilde{N}}$ 
and of $\Z$. We can consider a similar non-canonical identification of $I_N$ with 
$\widetilde{I}_N\subset \Z$ given by $[-M,-N_-[\cup [N_+,-M+\widetilde{N}[$, 
$\widetilde{N} = N + N_- + N_+$, where we choose 
$M$ so that $\Theta N \leq M \leq (1-\Theta)N$ for some $\Theta\in]0,1[$, with $N\gg 1$. Then, 
\eqref{gu11} has a more explicit form: 
\begin{equation}\label{gu13}
	\begin{split}
		&E^N(z;\nu,\mu) = \sum_{j\in\Z}
				          E^{\infty}(z;\nu+j\widetilde{N},\mu), 
				          \quad \mu,\nu\in \widetilde{I}_N, \\
		&E^N_+(z;\nu,\mu) = \sum_{j\in\Z}
				          E_+^{\infty}(z;\nu+j\widetilde{N},\mu), 
				          \quad \nu\in \widetilde{I}_N, ~ \mu\in J, \\
		&E^N_-(z;\nu,\mu) = \sum_{j\in\Z}
				          E_-^{\infty}(z;\nu,\mu+j\widetilde{N}), 
				          \quad \nu\in J, ~\mu\in \widetilde{I}_N, \\
		&E^N_{-+}(z;\nu,\mu) = \sum_{j\in\Z}
				          E_{-+}^{\infty}(z;\nu+j\widetilde{N},\mu), 
				          \quad \nu,\mu\in J.				    
	\end{split}
\end{equation}
In particular, due to the exponential decay, 
\begin{equation}\label{gu14}
	\begin{split}
		&E^N_+(z;\nu,\mu) =E_+^{\infty}(z;\nu,\mu) 
					       + \mO\left( \e^{-\frac{N}{C}}\right), 
				          \quad \nu\in \widetilde{I}_N, ~ \mu\in J, \\
		&E^N_-(z;\nu,\mu) = E_-^{\infty}(z;\nu,\mu) 
					       + \mO\left( \e^{-\frac{N}{C}}\right), 
				          \quad \nu\in J, ~\mu\in \widetilde{I}_N, \\
		&E^N_{-+}(z;\nu,\mu) = E_{-+}^{\infty}(z;\nu,\mu) 
					             + \mO\left( \e^{-\frac{N}{C}}\right), 
				          \quad \nu,\mu \in J.  
	\end{split}
\end{equation}
We next look at some general properties of $E_{-+}^N$. We are mainly interested in the case 
$N=+\infty$, but the discussion holds for all $N$, so we drop the superscript $N$. 
From $(P-z)E_+ + R_-E_{-+} = 0$, conclude that 
\begin{equation}\label{gu15}
	\ker(E_{-+}) \stackrel{E_+}{\longrightarrow} \ker(P-z).
\end{equation}
From $E_{-+}R_+ + E_{-}(P-z) = 0$ we see that 
\begin{equation}\label{gu15.5}
	\ker(P-z) \stackrel{R_+}{\longrightarrow} \ker(E_{-+}).
\end{equation}
Also notice that since $R_+E_+ + R_{+-}E_{-+} = 1$, we have 
\begin{equation}\label{gu16}
	R_+E_+ =1 \text{ on } \ker(E_{-+}).
\end{equation}
Similarly for $E(P-z)+E_+R_+=1$, we have 
\begin{equation}\label{gu17}
	E_+R_+ =1 \text{ on } \ker(P-z),
\end{equation}
so \eqref{gu15}, \eqref{gu15.5} are bijective and inverse to each other. 
\\
\par
Let $N=+\infty$. From Proposition \ref{g_prop3} we know that 
\begin{enumerate}
	\item if $N_+=m_+$, then $N_-=m_-$, by \eqref{g6}, and 
			\begin{equation*}
				\ker(P_{I_{\infty}}-z) =0.
			\end{equation*}
		Since $\ker(E_{-+}^{\infty})=R_+\ker(P_{I_{\infty}}-z)$ we 
		conclude that $E_{-+}^{\infty}$ is injective and hence bijective.
	\item if $N_+ < m_+$, then $N_->m_-$ and 
		\begin{equation*}
			u \in \ker(P_{I_{\infty}}-z) \iff 
			\begin{cases}
				u|_{]-\infty,-N_--1]}=0, \\
				u\in\ker(P_{[N_+,+\infty[}-z).
			\end{cases}
		\end{equation*}
		Moreover, 
		\begin{equation}\label{g17.1}
			\begin{split}
			\dim \ker(E_{-+}^{\infty}(z)) &=\dim \ker(P_{I_{\infty}}-z)\\
			&=\dim \ker(P_{[N_+,+\infty[}-z)=m_+-N_+.
			\end{split}
		\end{equation}
	\item if $N_+>m_+$, then $N_-<n_-$ and 
		\begin{equation*}
			u \in \ker(P_{I_{\infty}}-z) \iff 
			\begin{cases}
				u|_{[N_+,+\infty[}=0, \\
				u\in\ker(P_{]-\infty,-N_--1]}-z).
			\end{cases}
		\end{equation*}
		Moreover,
		\begin{equation*}
			\begin{split}
			\dim \ker(E_{-+}^{\infty}(z)) &=\dim \ker(P_{I_{\infty}}-z)\\
			&=\dim \ker(P_{]-\infty,-N_--1]}-z)=m_--N_-.
			\end{split}
		\end{equation*}
\end{enumerate}
In all cases $\ker(E_{-+}^{\infty}(z))=R_+ \ker(P_{I_{\infty}}-z)$. 
\par
Suppressing again the superscripts, we can describe by duality 
$\mathcal{R}(E_{-+})^{\perp}=\ker(E_{-+}^*)$. In fact by \eqref{gu8}
\begin{equation}\label{gu18}
	\begin{split}
		&R_-^*= 1_{J}\, p(\tau)^*\, 1_{I_N}, \quad R_+^* = 1_{I_N}\, p(\tau)^* \,1_{J}\\ 
		&R_{+-}^*= 1_{J}\, (p(\tau)-z)^* \,1_{J}. \\ 
	\end{split}
\end{equation}
So 
\begin{equation}\label{gu19}
	\mathcal{P}^*(z) = 
	\begin{pmatrix}
		(P_I-z)^* & R_+^* \\ 
		R_-^* & R_{+-}^* \\
	\end{pmatrix}
\end{equation}
is obtained from $(P_{S_{N+N_-+N_+}}-z)^*$ in exactly the same way as $\mathcal{P}(z)$ from 
$(P_{S_{N+N_-+N_+}}-z)$, cf. \eqref{gu7}. The inverse is 
\begin{equation}\label{gu20}
	\mathcal{E}^*(z) = 
	\begin{pmatrix} 
		E(z)^* & E_-(z)^* \\
		E_+(z)^* & E_{-+}(z)^* 
	\end{pmatrix},
\end{equation}
and we get from $N=+\infty$ that 
$\ker(E_{+-}(z)^*)=(R_-)^* \ker((P_{I_{\infty}}-z)^*)$. 
\\
\par
For $u\in\ell^2(\Z)$ let $\Gamma u = \overline{u}$. In view of \eqref{int.0} we see that 
	\begin{equation*}
		(p(\tau)-z)^* = \Gamma (p(\tau^{-1})-z) \Gamma
	\end{equation*}
as operators acting on $\ell^2(\Z)$. By Proposition \ref{g_prop3} we get  
\begin{enumerate}
	\item if $N_+=m_+$, then $N_-=m_-$, by \eqref{g6}, and 
			\begin{equation*}
				\ker\big((P_{I_{\infty}}-z)^*\big)=0.
			\end{equation*}
		Since $\ker((E_{-+}^{\infty})^*)=R_-^*\ker((P_{I_{\infty}}-z)^*)$, we 
		conclude that $(E_{-+}^{\infty})^*$ is injective and hence bijective.
	\item if $N_+ < m_+$, then $N_->m_-$ and 
		\begin{equation*}
			u \in \ker\big((P_{I_{\infty}}-z)^*\big) \iff 
			\begin{cases}
				u|_{[N_+,\infty[}=0, \\
				u\in\ker\big((P_{]-\infty,-N_--1]}-z)^*\big).
			\end{cases}
		\end{equation*}
		Moreover, 
		\begin{equation*}
			\begin{split}
			\dim \ker\big(E_{-+}^{\infty}(z)^*\big)&
			=\dim \ker\big((P_{I_{\infty}}-z)^*\big) \\
			&=\dim \ker\big((P_{]-\infty,-N_--1]}-z)^*\big)
			=m_+-N_+.
			\end{split}
		\end{equation*}
	\item if $N_+>m_+$, then $N_-<m_-$ and 
		\begin{equation}\label{g20.1}
			u \in \ker\big((P_{I_{\infty}}-z)^*\big) \iff 
			\begin{cases}
				u|_{]-\infty,-N_--1]}=0, \\
				u\in\ker\big((P_{[N_+,+\infty[}-z)^*\big).
			\end{cases}
		\end{equation}
		Moreover, 
		\begin{equation*}
			\begin{split}
			\dim \ker\big(E_{-+}^{\infty}(z)^* \big)&=
			\dim \ker\big((P_{I_{\infty}}-z)^*\big) \\
			&=\dim \ker\big((P_{[N_+,+\infty[}-z)^*\big)=m_--N_-.
			\end{split}
		\end{equation*}
\end{enumerate}
\subsection{Estimates on the singular values of $E_{\pm}$}
In this section we will give bounds on the singular values of $E_{\pm}$, 
see \eqref{gu10}. We will treat both the case when $N\geq N_++N_-+1.$ and the 
limiting case when $N=+\infty$. First, notice that 
\begin{equation}\label{gu21.0}
	\mathrm{rank}(E_{\pm}^{N}) \leq |J| = N_+ + N_-.
\end{equation}
When $N\geq N_++N_-+1.$ and $z\notin \spec(P_{S_{\widetilde{N}}})$, 
let 
\begin{equation}\label{gu21.1}
	0 \leq s^{N,\pm}_{|J|} \leq \dots \leq s^{N,\pm}_1 = \|E_{\pm}^{N}\|
\end{equation}
denote the singular values of $E_{\pm}^{N}$. When $N=+\infty$ and $z\notin p(S^1)$, 
let 
\begin{equation}\label{gu21.2}
	0 \leq s^{\infty,\pm}_{|J|} \leq \dots \leq s^{\infty,\pm}_1 = \|E_{\pm}^{\infty}\|
\end{equation}
denote the singular values of $E_{\pm}^{\infty}$. Although we have not denoted it 
explicitly here, the singular values \eqref{gu21.2}, \eqref{gu21.2}, depend 
on $z$. Recall \eqref{gu7} and notice that since 
the operator $p(\tau)$ acting on $\ell^2(S_{\widetilde{N}})$ and on  $\ell^2(\Z)$ is normal, we have 
the trivial upper bounds 
\begin{equation}\label{gu21.3}
	s^{N,\pm}_1\leq \frac{1}{\dist(z,\spec(P_{S_{\widetilde{N}}}))}, \quad
	s^{\infty,\pm}_1  \leq \frac{1}{\dist(z,p(S^1))}
\end{equation}
\begin{lemma}\label{lem:SV} 
Let $N\geq 2(N_++N_-)+1$ and let $\Omega\Subset\C$ be a compact set. Then, 
\begin{enumerate}
\item there exists a constant $C>0$, such that for all $z\in\Omega\backslash \spec(P_{S_{\widetilde{N}}})$
\begin{equation*}
	\frac{1}{C} \leq s_j^{N,\pm} \leq  
	\frac{1}{\dist(z,\spec(P_{S_{\widetilde{N}}}))}, \quad j = 1,\dots, N_++N_-.
\end{equation*}
In particular $E_+^N$ is injective and $E_-^N$ is surjective. 
\item there exists a constant $C>0$, such that for all $z\in\Omega\backslash p(S^1)$
\begin{equation*}
	\frac{1}{C} \leq s_j^{\infty,\pm} \leq  
	\frac{1}{\dist(z,p(S^1))}, \quad j = 1,\dots, N_++N_-.
\end{equation*}
In particular $E_+^{\infty}$ is injective and $E_-^{\infty}$ is surjective. 
\end{enumerate}
\end{lemma}
\begin{remark} 
Notice that in both cases the lower bound on the singular values 
only depends on the compact set $\Omega$ and is independent of $N$. 
This is due to the fact that the only moment 
in the proof of Lemma \ref{lem:SV} where we need that $z\notin \spec(P_{S_{\widetilde{N}}})$ 
(respectively  $z\notin p(S^1)$ when $N=+\infty$) is 
when we use that $\mathcal{E}$ \eqref{gu9}- the inverse of the Grushin problem $\mathcal{P}$ 
\eqref{gu7} -  exists, see \eqref{sv19} below.
\end{remark}
\begin{proof}[Proof of Lemma \ref{lem:SV}]
We begin with the case (1): The upper bounds follow from \eqref{gu21.3}. 
\par
Let us now turn to the lower bounds. We begin by recalling the Grushin problem 
\eqref{gu7}: for $z\in\Omega\backslash \spec(P_{S_{\widetilde{N}}})$, the operator %
\begin{equation*}
	(p(\tau)-z) : \ell^2(S_{\widetilde{N}}) \longrightarrow \ell^2(S_{\widetilde{N}})
\end{equation*}
is bijective with bounded inverse $\mathcal{E}_N(z)$, see \eqref{gu9}. Here, 
$S_{\widetilde{N}}\simeq \Z/\widetilde{N}\Z$, $\widetilde{N} = N + N_+ + N_-$. 
Recall the notation introduced in the 
discussion after \eqref{gu1} where we write segments of $S_{\widetilde{N}}$ as 
intervals modulo $\widetilde{N}\Z$. We write 
\begin{equation*}
	S_{\widetilde{N}} = J\cup I_N
\end{equation*}
where $J=[-N_-,N_+[$ is naturally defined both as a subset of $S_{\widetilde{N}}$ and of $\Z$. 
For $I_N$ we write 
\begin{equation*}
	 I_N  = S_{\widetilde{N}}\backslash J \equiv [N_+,-N_--1]   \subset S_{\widetilde{N}}.
\end{equation*}
Moreover, we will use the notation $a+J = [a - N_-, a+ N_+[$,  
$a\in S_{\widetilde{N}}$.
\\
\par
Next, suppose that $z\in\Omega$ and let 
\begin{equation}\label{svn0}
	(p(\tau) -z ) u = v \quad \text{on }  S_{\widetilde{N}}, \quad 
	\text{with } \supp v \subset J.
\end{equation}
Fix $a_+,a_-\in S_{\widetilde{N}}\backslash J$, so that 
\begin{equation}\label{svn1}
	N_+ + N_-+1 \leq \dist_{S_{\widetilde{N}}}(a_{+},N_{+}-1) 
	=
	\dist_{S_{\widetilde{N}}}(a_{+},J)=	\mO(1)
\end{equation}	
and 
\begin{equation}\label{svn2}
	 N_+ + N_-+1 \leq \dist_{S_{\widetilde{N}}}(a_{-},-N_{-})  =
 	\dist_{S_{\widetilde{N}}}(a_{-},J)=\mO(1).
\end{equation}	
Notice that 
\begin{equation}\label{svn3}
 	(a_{\pm} +[-N_--N_+, N_+ +N_-]) \cap J = \emptyset.
\end{equation}	
By \eqref{svn0} we see that 
\begin{equation}\label{svn4}
	(p(\tau) -z ) \mathbf{1}_{ [a_-,a_+[}\, u =
	\begin{cases}
		0,  \text{ on } S_{\widetilde{N}}\backslash [a_--N_-,a_++N_+[\, , \\
		v, \text{ on }  [a_-+N_+,a_+-N_-[\, , \\
		w_-, \text{ on } a_- +J\, , \\
		w_+, \text{ on }  a_++ J\, ,
	\end{cases}
\end{equation}
where $w_{\pm}\in \ell^2(S_{\widetilde{N}})$ and $\supp w_{\pm} \subset a_{\pm} +J$. 
Since $\supp v \subset J$, we see by 
\eqref{svn1}, \eqref{svn2} and \eqref{svn4}, that  
\begin{equation}\label{svn4.1}
	(p(\tau) -z ) \mathbf{1}_{[a_-,a_+[}\, u = v + w_+ + w_-\,, 
\end{equation}
and 
\begin{equation}\label{svn5}
	\| w_{\pm} \| \leq \mO(1) \| \mathbf{1}_{a_{\pm}+[-N_+-N_-,N_++N_-[}\, u\|.
\end{equation}
Next, write 
\begin{equation}\label{svn6}
	\tau^{-N_+}(p(\tau) -z ) \mathbf{1}_{[a_-,a_+[}\, u = \tau^{-N_+}(v + w_+ + w_-)
\end{equation}
and
\begin{equation}\label{svn7}
	\tau^{N_-}(p(\tau) -z ) \mathbf{1}_{[a_-,a_+[}\, u = \tau^{N_-}(v + w_+ + w_-). 
\end{equation}
We will use these two equations to estimate $\| \mathbf{1}_{[0,N_+[}\, u\|$, when 
$N_+\geq 1$, and $\| \mathbf{1}_{[-N_-,0[}\, u\|$, when $N_-\geq 1$. 
\\
\par
In view of \eqref{int.0}, \eqref{g9}, we see that $\tau^{-N_+}(p(\tau) -z ) $ is upper triangular 
with a non-vanishing constant entry at the diagonal. Since 
$\supp \tau^{-N_+}(v + w_+ + w_-) \subset [a_--N_+-N_-,a_+[$ and 
$\supp \mathbf{1}_{[a_-,a_+[}\, u \subset [a_-,a_+[$, we see that 
\begin{equation}\label{svn8}
	\| \mathbf{1}_{[0,a_+[}\, u \| \leq \mO(1) 
	\| \mathbf{1}_{[0,a_+[}\tau^{-N_+}(v + w_+ + w_-)\|,
\end{equation}
where the constant is uniform in $z\in\Omega$ and independent of $N$. Here, 
\begin{equation*}
	\mathbf{1}_{[0,a_+[}\tau^{-N_+}(v + w_+ + w_-) = \mathbf{1}_{[0,a_+[}\tau^{-N_+} w_+\, ,
\end{equation*}
so, by \eqref{svn8}, \eqref{svn5}, 
\begin{equation*}
	\| \mathbf{1}_{[0,a_+[}\, u \|
	\leq \mO(1)  \| \mathbf{1}_{a_++[-N_+-N_-,N_++N_-[}\, u\| 
\end{equation*}
which, using \eqref{svn1}, implies
\begin{equation}\label{svn9}
	\| \mathbf{1}_{[0,a_++(N_++N_-)[}\, u \| 
	\leq \mO(1) \| \mathbf{1}_{[N_+,a_+ + N_++N_-[}\, u\|.
\end{equation}
Notice that when $N_+=0$ this holds trivially. 
\\
\par
When $N_-\geq 1$, we use that $\tau^{N_-}(p(\tau) -z )$ is lower triangular with 
a non-vanishing constant entry at the diagonal. In \eqref{svn7} we have that 
$\supp \tau^{N_-}(v + w_+ + w_-) \subset [a_-,a_++N_++N_-[$ and 
$\supp \mathbf{1}_{[a_-,a_+[}\, u \subset [a_-,a_+[$. 
We therefore 
deduce that 
\begin{equation}\label{svn10}
	\| \mathbf{1}_{[a_-,0[}\, u \| \leq \mO(1) 
	\| \mathbf{1}_{[a_-,0[}\tau^{N_-}(v + w_+ + w_-)\|,
\end{equation}
where the constant is uniform in $z\in\Omega$ and independent of $N$. 
Since 
\begin{equation*}
\mathbf{1}_{[a_-,0[}\tau^{N_-}(v + w_+ + w_-) = \mathbf{1}_{[a_-,0[}\tau^{N_-} w_-,
\end{equation*}
we obtain by \eqref{svn10}, \eqref{svn5}, \eqref{svn2} that 
\begin{equation}\label{svn11}
	\| \mathbf{1}_{[a_--(N_++N_-),0[}\, u \| \leq \mO(1) 
	 \| \mathbf{1}_{[a_--(N_++N_-),-N_-[}\, u\|,
\end{equation}
which holds trivially when $N_- =0$.  
\\
\par
Combining \eqref{svn9}, \eqref{svn11} gives 
\begin{equation}\label{svn12}
	\| \mathbf{1}_{[a_--(N_++N_-),a_++(N_++N_-)[}\, u \| \leq \mO(1) 
	 \| \mathbf{1}_{[a_--(N_++N_-),a_++(N_++N_-)[\backslash J} \, u\| .
\end{equation}
Since $v=(p(\tau)-z)u$ is supported in $J$, we have that
\begin{equation*}
\begin{split}
\| v\| &\leq \mO(1) 
	\| \mathbf{1}_{[-(N_++N_-),(N_++N_-)[}\, u \| \\
	& \leq \mO(1)
	 \| \mathbf{1}_{[a_--(N_++N_-),a_++(N_++N_-)[} \, u\|, 
\end{split}
\end{equation*}
where the constant in the estimate is uniform in $z\in \Omega$ and independent of $N$. 
Combining this with \eqref{svn12} shows that 
\begin{equation*}
\| v\| \leq \mO(1) 
	 \| \mathbf{1}_{[a_--(N_++N_-),a_++(N_++N_-)[\backslash J} \, u\|.
\end{equation*}
Now suppose that $z\in\Omega\backslash \spec(P_{S_{\widetilde{N}}})$ and recall from \eqref{gu7}, \eqref{gu9}, that when $u\in \ell^2(S_{\widetilde{N}})$, we have that 
\begin{equation}\label{sv19}
u=\mathcal{E}_N(z)v, \quad \text{with } v =  \mathbf{1}_{J}v_+, ~ v_+\in\ell^2(J).
\end{equation}
Hence, by \eqref{gu9}, $u = E_+^N v_+$ on $I_N =S_{\widetilde{N}}\backslash J$. Thus, 
\begin{equation*}
\begin{split}
\| v_+\| & \leq \mO(1) 
	\| \mathbf{1}_{[a_--(N_++N_-),a_++(N_++N_-)[\backslash J} \, E_+^N v_+\| \\
	& \leq \mO(1) \| \mathbf{1}_{I_N} \, E_+^N v_+\|,
\end{split}
\end{equation*}
where the constant in the estimate is uniform in $z\in\Omega$ and independent of $N$. 
This concludes the proof for the singular values of $E_+^N$. The proof of the statement 
for $E_-^N$ follows exactly the same lines using $(E_-^N)^*$ instead of $E_+^N$. 
\par
The proof of the statement in the case (2), when $N=\infty$, is similar, using that 
$S_{\infty}\simeq \Z = ]-\infty,N_-[\cup [N_-,N_+[ \cup [N_+,+\infty[$.
\end{proof}
\section{A Grushin Problem for the perturbed operator}\label{sec:GP}
Our aim is to study the following random perturbation of $P_0=P_\mathrm{I_N}$: 
\begin{equation}\label{gp1}
 P_N^{\delta} \defeq P_N^0 + \delta Q_{\omega}, 
 \quad
 Q_{\omega}=(q_{j,k}(\omega))_{1\leq j,k\leq N}, 
\end{equation}
where $0\leq\delta\ll 1 $ and $q_{j,k}(\omega)$ are independent and 
identically distributed complex Gaussian random variables, 
following the complex Gaussian law $\mathcal{N}_{\C}(0,1)$. Here, 
$1 \ll N < \infty$. Consider the space 
$\mathcal{H}_N \defeq (\C^{N\times N}, \| \cdot \|_{\mathrm{HS}})$ 
of $N\times N$ complex valued matrices equipped with the Hilbert-Schmidt norm. 
We equip $\mathcal{H}_N$ with the probability measure 
\begin{equation}\label{gp1.0}
	\mu_N (dQ) \defeq \pi^{-N^2} \e^{- \| Q\|_{\mathrm{HS}}^2} L(dQ), 
\end{equation} 
where $L(dQ)$ denotes the Lebesgue measure on $\mathcal{H}_N$. 
For $C_1>0$, let $\mathcal{Q}_{C_1N}\subset \mathcal{H}_N$ be the 
subset where 
\begin{equation}\label{gp1.1}
	\| Q\|_{\mathrm{HS}} \leq C_1 N.
\end{equation} 
Markov's inequality \cite[Lemma 3.1]{Kal97} implies that if $C_1>0$ is large enough,
\begin{equation}\label{gp2}
\mathds{P}\left[ 
\Vert Q_{\omega}\Vert_\mathrm{HS}\le C_1N
\right] 
= \mu_N(\mathcal{Q}_{C_1N}) \ge 1-e^{-N^2}.
\end{equation}
\subsection{A general discussion}\label{GPP0}
We begin with a formal discussion of a Grushin problem for the 
perturbed operator $P_{\delta}$. Recall from Section \ref{grush} 
that the Grushin problem for the unperturbed operator is of the form 
\begin{equation*}
 \mathcal{P}_0=\begin{pmatrix}
 		P_0-z & R_-\\
                 R_+ & R_{+-} \\
                \end{pmatrix} 
                :~\ell^2({I_N})\times \ell^2(J) \longrightarrow \ell^2({I_N})\times \ell^2(J),
\end{equation*}
We added a subscript $0$ to indicate that we deal with the 
unperturbed operator. Suppose that $\mathcal{P}_0$ is bijective with inverse
\begin{equation*}
 \mathcal{E}_0=\begin{pmatrix}
                 E^0 & E_+^0 \\
                 E_-^0 & E_{-+}^0\\
                \end{pmatrix}
                :~\ell^2({I_N})\times \ell^2(J)\longrightarrow \ell^2({I_N})\times \ell^2(J),
\end{equation*}
where we added a superscript $0$ for the same reason. 
Supposing that 
\begin{equation}\label{gp3}
	\|\delta Q_{\omega}\| \|E^0\| <1, 
\end{equation}
we see by a Neumann series argument that 
\begin{equation*}
 \mathcal{P}_{\delta} \defeq
	       \begin{pmatrix}
 		P_{\delta} - z & R_-\\
                 R_+ & R_{+-} \\
                \end{pmatrix}
                :~\ell^2({I_N})\times \ell^2(J)\longrightarrow \ell^2({I_N})\times \ell^2(J),
\end{equation*}
is bijective and admits the inverse
\begin{equation*}
 \mathcal{E}_{\delta}=\begin{pmatrix}
                 E^{\delta} & E_+^{\delta} \\
                 E_-^{\delta} & E_{-+}^{\delta} \\
                \end{pmatrix}
                :~\ell^2({I_N})\times \ell^2(J) \longrightarrow \ell^2({I_N})\times \ell^2(J),
\end{equation*}
where 
\begin{equation}\label{gp4}
\begin{split}
 &E_+^{\delta} 
 = (1 +  E^0(\delta Q_{\omega}))^{-1}E_+^0, \\
 &E_-^{\delta} 
 = E_-^0(1 + \delta Q_{\omega}E^0)^{-1},\\
 &E^{\delta} 
 = E^0(1 + \delta Q_{\omega}E^0)^{-1}, \\
 &E_{-+}^{\delta} 
 =  E_{-+}^0 -  E_-^0 \delta Q_{\omega}(1 +  E^0(\delta Q_{\omega}))^{-1}E_+^0. 
 \end{split}
\end{equation}
One obtains the following estimates
\begin{equation}\label{gp5}
\begin{split}
 & \|E^{\delta}\| \leq 
 \frac{\| E^0\|}{1 - \|\delta Q_{\omega}\| \|E^0\|}, ~
 \|E_{\pm}^{\delta}\| \leq 
 \frac{\| E_{\pm}^0\|}{1 - \|\delta Q_{\omega}\| \|E^0\|}, \\
 & \| E_{-+}^{\delta} - E_{-+}^{0}\| \leq  
 \frac{\| E_{+}^0\| \| E_{-}^0\| \|\delta Q_{\omega}\|}{1 -\|\delta Q_{\omega}\| \|E^0\|}.
\end{split}
\end{equation}
Differentiating the equation $\mathcal{E}^{\delta}\mathcal{P}^{\delta}=1$ with 
respect to $\delta$ yields
\begin{equation}\label{gp6}
 \partial_{\delta}\mathcal{E}^{\delta} = 
 - \mathcal{E}^{\delta}(\partial_{\delta}\mathcal{P}^{\delta})\mathcal{E}^{\delta}
 =
 -\begin{pmatrix}
   E^\delta Q_{\omega} E^{\delta} & E^\delta Q_{\omega} E_+^{\delta} \\
    E_-^\delta Q_{\omega} E^{\delta} & E_-^\delta Q_{\omega} E_+^{\delta} \\
  \end{pmatrix}.
\end{equation}
Integrating this relation from $0$ to $\delta$ yields 
\begin{equation}\label{gp7}
 \|E^{\delta} -E^0\| \leq 
 \frac{\|\delta Q_{\omega}\| \| E^0\|^2}{(1 - \|\delta Q_{\omega}\| \|E^0\|)^2}, ~
 \|E_{\pm}^{\delta} -E_{\pm}^0\| \leq 
 \frac{\|\delta Q_{\omega}\| \| E_{\pm}^0\|\|E^0\|}{(1 - \|\delta Q_{\omega}\| \|E^0\|)^2}.
\end{equation}
\par
Since $\mathcal{P}^{\delta}$ is invertible and of finite 
rank, we know that 
\begin{equation*}
 |\partial_{\delta} \ln\det\mathcal{P}^{\delta}| 
 = |\mathrm{tr}(\mathcal{E}^{\delta}\partial_{\delta}\mathcal{P}^{\delta})|.
\end{equation*}
Letting $\|\cdot\|_{\tr}$ denote the trace class norm, we get 
\begin{equation}\label{gp8}
 |\partial_{\delta} \ln\det\mathcal{P}^{\delta}| 
 = |\mathrm{tr}(Q_{\omega}E^{\delta})|
 \leq \|Q_{\omega}\|_{\mathrm{tr}} \|E^{\delta}\|
 \leq \frac{\| E^0\| \|Q_{\omega}\|_{\mathrm{tr}}}{1 - \|\delta Q_{\omega}\| \|E^0\|},
\end{equation}
where $\|Q_{\omega}\|_{\mathrm{tr}} \leq N^{1/2}\|Q_{\omega}\|_{\HS}$. Integration from 
$0$ to $\delta$ yields
\begin{equation}\label{gp9}
 \left| \ln |\det\mathcal{E}^{\delta}| - \ln |\det\mathcal{E}^{0}| \right| = 
 \left| \ln |\det\mathcal{P}^{\delta}| - \ln |\det\mathcal{P}^{0}| \right|
 \leq 
 \frac{ \| E^0\| \|\delta Q_{\omega}\|_{\mathrm{tr}}}{1 - \|\delta Q_{\omega}\| \|E^0\|}.
\end{equation}
\par
Sharpening the assumption \eqref{gp3} to 
\begin{equation}\label{gp10}
  \|\delta Q_{\omega}\| \|E^0\| < \frac{1}{2},
\end{equation}
we get
\begin{equation}\label{gp11}
 \|E^{\delta}\| \leq 2\| E^0\|, ~
 \|E_{\pm}^{\delta}\| \leq 2\| E_{\pm}^0\|,~ 
\| E_{-+}^{\delta} - E_{-+}^{0}\| \leq  
 2 \| E_{+}^0\| \| E_{-}^0\| \|\delta Q_{\omega}\|.
\end{equation}
By \eqref{gp6} we know that 
$\partial_{\delta}E_{-+}^{\delta} = - E_-^{\delta} Q_{\omega} E_+^{\delta}$. 
Therefore, using \eqref{gp5}, \eqref{gp7} and \eqref{gp11} we get 
\begin{equation}\label{gp12}
\begin{split}
 \|\partial_{\delta}E_{-+}^{\delta} +E_-^0 Q_{\omega} E_+^0\| 
  &\leq 
  \|E_-^0 Q_{\omega} \| \|E_+^{\delta}-E_+^0\| + 
  \|Q_{\omega} E_+^{\delta}  \| \|E_-^{\delta}-E_-^0\| \\
  &\leq 12\delta   \|Q_{\omega}\|^2 \|E_-^0\| \|E_+^0\| \| E^0\|.
  \end{split}
\end{equation}
By integration from $0$ to $\delta$, we conclude 
\begin{equation}\label{gp13}
 E_{-+}^{\delta} = E_{-+}^{0} - E_-^0(\delta  Q_{\omega}) E_+^0 
 + \mO(\|\delta Q_{\omega}\|^2 \|E_-^0\| \|E_+^0\| \| E^0\|).
\end{equation}
\subsection{A Grushin problem for the perturbed operator}\label{grupp}
Recall from \eqref{gu7} that 
\begin{equation*}
	\mathcal{P}_N(z)=(P_{S_{\widetilde{N}}}-z), \quad \widetilde{N}=N+N_-+N_+ 
\end{equation*}
and from \eqref{g27.5} that its spectrum is equal to $p(\widehat{S}_{\widetilde{N}})$. 
Suppose that 
$z\notin\spec(P_{S_{\widetilde{N}}})$. As in \eqref{gu9}, $\mathcal{P}_N(z)$ is 
invertible with bounded inverse $\mathcal{E}_N(z)$. 
\\
\par
Suppose that 
\begin{equation}\label{gp14}
	\dist(z,\spec(P_{S_{\widetilde{N}}})) \geq \frac{1}{CN}
\end{equation}
for some fixed sufficiently large constant $C>1$ to be determined later on. 
Since the operator $\mathcal{P}_N(z)$ is normal, it follows that 
\begin{equation}\label{gp15}
	\| \mathcal{E}_N(z)\| = \frac{1}{\dist(z,\spec(P_{S_{\widetilde{N}}}))}.
\end{equation}
In particular
\begin{equation}\label{gp16}
	\| E^N(z)\| ,\| E_-^N(z)\| ,\| E_+^N(z)\| , 
	\| E_{-+}^N(z)\| \leq \frac{1}{\dist(z,\spec(P_{S_{\widetilde{N}}}))}.
\end{equation}
Suppose that 
\begin{equation}\label{gp17}
	0< \delta \ll N^{-2}.
\end{equation}
Then, by \eqref{gp2}, \eqref{gp16}, \eqref{gp14}, with probability 
$\geq 1- \e^{-N^2}$, the assumption \eqref{gp10} is satisfied. Therefore, 
by the discussion in Section \ref{GPP0} we conclude 
\begin{prop}\label{gp:prop1}
With probability $\geq 1 - \exp(-N^2)$ we have: 
Suppose \eqref{gp14}, \eqref{gp17}. Let $\mathcal{P}_N^0(z)= \mathcal{P}_N(z)$ 
be as in \eqref{gu7} and let $\mathcal{E}_N^0(z)= \mathcal{E}_N(z)$ 
be as in \eqref{gu9}. Then, 
\begin{equation*}
	\mathcal{P}_N^{\delta}(z) \defeq
	\begin{pmatrix} 
		P_{I_N}^{\delta} -z & R_-^N \\
		R_+^N & R_{+-}^N(z) 
	\end{pmatrix}
	: \ell^2(I_N) \oplus \ell^2(J) \to \ell^2(I_N) \oplus \ell^2(J)
\end{equation*}
is bijective with bounded inverse 
\begin{equation*}
	\mathcal{E}_N^{\delta}(z) = 
	\begin{pmatrix} 
		E^{N,\delta}(z) & E_+^{N,\delta}(z) \\
		E_-^{N,\delta}(z) & E_{-+}^{N,\delta}(z) 
	\end{pmatrix}
	: \ell^2(I_N) \oplus \ell^2(J) \to \ell^2(I_N) \oplus \ell^2(J).
\end{equation*}
Moreover, 
\begin{equation*}
\| E^{N,\delta}(z)\| ,\| E_-^{N,\delta}(z)\| ,\| E_+^{N,\delta}(z)\|,
 \| E_{-+}^{N,\delta}(z)\|    \leq \frac{2}{\dist(z,\spec(P_{S_{\widetilde{N}}}))}.
\end{equation*}
\end{prop}
\subsection{A lower bound on the determinant of the effective Hamiltonian} \label{sec:LB}
Suppose that $\Omega\Subset\C$ is a compact set. Let $E_{-+}^{N,\delta}$ 
be as in Proposition \ref{gp:prop1}. In this section we are interested in 
estimating the probability that $\log |\det E_{-+}^{N,\delta}(z)| \leq a $ 
for $a\in\R$ and for some $z\in\Omega \backslash p(\widehat{S}_{\widetilde{N}})$ which 
may depend on $N$. To obtain this bound we will adapt the approach 
developed in \cite[Section 9]{HaSj08}.
\par
Set 
\begin{equation}\label{lb1}
	\alpha = \alpha (z;N) \defeq  \dist(z,p(\widehat{S}_{\widetilde{N}})).
\end{equation}
Until further notice we suppose that 
\begin{equation}\label{lb2}
	\alpha \geq \frac{1}{CN^{\kappa}}, \quad \text{ for some } C>1, 
\end{equation}
where $\kappa \geq 1 $ is fixed, and we strengthen assumption \eqref{gp17} to 
\begin{equation}\label{lb2.1}
	0 < \delta \ll N^{-1} \min( \alpha, N^{-1} ).
\end{equation}
Recall Proposition \ref{gp:prop1} and \eqref{gp4}. We 
want to study the map 
\begin{equation}\label{lb3}
	\begin{split}
		\mathcal{Q}_{C_1N} \ni Q \mapsto E_{-+}^{\delta}(z,Q) 
		&= E_{-+}^0(z) - \delta E_-^0(z) \left( Q 
		+ \sum_1^{\infty} (-\delta)^n Q(E^0(z) Q)^n
			\right)E_+^0(z) \\
		& \defeq  E_{-+}^0(z) - \delta E_-^0(z) ( Q +T(z,Q,\delta,N)
			)E_+^0(z)
	\end{split}
\end{equation}
where by \eqref{gp16}, \eqref{gp1.1}, 
\begin{equation}\label{lb4}
	\| T\|_{\mathrm{HS}} \leq \mO\!\left(\frac{\delta(C_1N)^2}{\alpha}\right).
\end{equation}
Next, recall \eqref{gp1.0}, and notice that the measure $\mu_N$ is invariant 
under the left and right action of the group of unitary matrices $\mathcal{U}(N,\C)$ 
on $\mathcal{H}_N$, 
i.e. for any $U,V \in \mathcal{U}(N,\C)$, we have that 
\begin{equation}\label{lb5}
	\mu_N(d(UQV)) = \mu_N(dQ). 
\end{equation}
Furthermore, the left and right action of the group of unitary matrices $\mathcal{U}(N,\C)$ 
leaves $\mathcal{Q}_{C_1N}$ invariant, see \eqref{gp1.1}, and therefore also the probability 
\eqref{gp2}. Thus, we may choose any orthonormal bases (ONB) to represent the matrix 
$Q\in\mathcal{H}_N$. Let $\widetilde{e}_1,\dots,\widetilde{e}_N$ and 
$\widehat{e}_1, \dots, \widehat{e}_N$ be two orthonormal bases of $\C^N$ and write 
\begin{equation}\label{lb6}
	Q = \sum_{i,j=1}^N q_{i,j} \,\widetilde{e}_i \circ \widehat{e}_j ^*, 
	\quad \text{where } q_{i,j} \sim \mathcal{N}_{\C}(0,1) ~ (\mathrm{iid}).
\end{equation}
By Lemma \ref{lem:SV} and \eqref{lb2}, we have for a compact set $\Omega\Subset\C$ and 
for $z\in\Omega\backslash p(\widehat{S}_{\widetilde{N}})$, the following bound on the singular values of $E^N_{\pm}$

\begin{equation}\label{lb6.1}
	\frac{1}{C} \leq s_j^{N,\pm} \leq  
	\frac{1}{\alpha}, \quad j = 1,\dots, |J|=N_++N_-,
\end{equation}
where the constant $C>0$ is uniform in $z\in\Omega$ and independent of $N$.   
\par
By the polar decomposition we write $E_+^0 = S_+D_+$ where $S_+:\C^{|J|}\to \C^N$ 
is an isometry, with $S_+^*S_+ = 1$ and $S_+S_+^*$ is the orthogonal projection 
$\C^N\to \mathcal{R}(E_+^0)$, and $D_+:\C^{|J|}\to\C^{|J|}$ is selfadjoint with eigenvalues 
$s_1^+,\dots,s_{|J|}^+$. Similarly, 
\begin{equation}\label{an1}
	(E_-^0)^* = S_-D_-, \quad  E_-^0 = D_-S_-^*,
\end{equation}
where $S_-:\C^{|J|}\to \C^N$ is an isometry, with $S_-^*S_- = 1$ and $S_-S_-^*$ is the orthogonal projection $\C^N\to \mathcal{R}((E_-^0)^*)$, and $D_-:\C^{|J|}\to\C^{|J|}$ is selfadjoint with eigenvalues $s_1^-,\dots,s_{|J|}^-$. 
\par
From \eqref{lb23}, we get 
\begin{equation}\label{an2}
\begin{split}
E_{-+}^{\delta} &= E_{-+}^0 - \delta D_-S_-^* ( Q  +  T ) S_+D_+\\ 
& = D_-\big(\widehat{E}_{-+}^0 - \delta ( S_-^*Q S_+ + S_-^*T S_+)\big)  D_+,
\end{split}		
\end{equation}
where $\widehat{E}_{-+}^0  = D_-^{-1}E_{-+}^0D_+^{-1}$. Moreover, set 
\begin{equation}\label{an2.1}
	\widehat{T}=S_-^*T S_+.	
\end{equation}

View $\C^{|J|}$ as a subspace of $\C^N$ by 
considering that $J\subset \{1,\dots, N\}$. Let $\Pi_0 : \C^N \to \C^{|J|}$ be the orthogonal 
projection and, whenever convenient, view $\Pi_0$ as the inclusion map 
$\Pi_0 : \C^{|J|} \xhookrightarrow{} \C^N$. Let $\mathcal{S}_+: \C^N \to \C^N $ 
be unitary with $\mathcal{S}_+|_{\C^{|J|}} = S_+ $ and similarly for $\mathcal{S}_-$. 
Then, 
\begin{equation}\label{an2.2}
 S_+ = \mathcal{S}_+ \Pi_0,
\end{equation}
where $\Pi_0$ is viewed as a map $\C^{|J|} \to \C^N$. Similarly, 
\begin{equation}\label{an2.3}
 S_- = \mathcal{S}_- \Pi_0, \quad  S_-^* = \Pi_0 \mathcal{S}_-^* = \Pi_0 S_-^{-1}.
\end{equation}
Then, 
\begin{equation}\label{an3}
E_{-+}^{\delta} 
= D_-\big(\widehat{E}_{-+}^0 - \delta ( \Pi_0\widehat{Q} \Pi_0 
+ \widehat{T})\big)  D_+,	
\quad \widehat{Q} =  \mathcal{S}_-^*Q \mathcal{S}_+.
\end{equation}
Let 
 $\delta_j\in \C^N$, with $\delta_j(i) = 1$ if $i = j$ and $=0$ else, denote 
the standard ONB of $\C^N$. For $k=1,\dots, N$ set  
\begin{equation*}
	\widehat{e}_k \defeq \mathcal{S}_+ \delta_k,\quad \widetilde{e}_k \defeq 
	\mathcal{S}_-^*\delta_k
\end{equation*}
in \eqref{lb6}. Hence,
\begin{equation*}
	\widehat{Q} = \mathcal{S}_-^*Q \mathcal{S}_+ = (q_{j,k})_{1 \leq j,k\leq N} 
\end{equation*}
where 
$q_{j,k} \sim \mathcal{N}_{\C}(0,1)$ are independent and identically 
distributed complex Gaussian random variables. 
\par
By \eqref{an2.1}, \eqref{an2.2} and \eqref{an2.3}, we see that 
$\widehat{T}(Q) =\Pi_0 \widehat{T}(Q) \Pi_0$ and that 
the map $ \mathcal{H}_N \ni Q \mapsto \widehat{T}(Q) \in \mathcal{H}_{|J|}$ 
satisfies 
\begin{equation}\label{lb15.2}
\| \widehat{T}(Q)  \|_{\mathrm{HS}} 
 	\leq  \mO\!\left(\frac{\delta(C_1N)^2}{\alpha}\right).
\end{equation}
where the estimate is uniform in $Q\in\mathcal{Q}_{C_1N}$. 
\\
\par 
By \eqref{an3}
\begin{equation}\label{lb16}
		\mathcal{Q}_{C_1N} \ni Q \mapsto \det E_{-+}^{\delta}(z,Q) 
		=\prod_{k=1}^{|J|} (s_k^+s_k^-) \det \left( \widehat{E}_{-+}^0(z) 
		- \delta ( \Pi_0 \widehat{Q} \Pi_0 +
		\widehat{T}(Q) )\right)
\end{equation}
Recall from \eqref{an1} and from the discussion after 
\eqref{lb6.1} that $s_{k}^{+}$ (resp. $s_{k}^{-}$) denote the singular 
values of $E_{+}^0$ (resp. $(E_-^0)^*$). 
\\
\par
The Cauchy inequalities and \eqref{lb15.2} imply that 
\begin{equation}\label{lb17}
 \| d_Q \widehat{T}\|_{\mathcal{H}_N\to  \mathcal{H}_{|J|} } 
 \leq \mO\!\left(\frac{\delta C_1N}{\alpha}\right),
\end{equation}
uniformly for $Q\in\mathcal{Q}_{C_1N}$. Technically, we can only apply the Cauchy inequalities 
in $\|Q\|_{\mathrm{HS}} \leq \eta\, C_1 N$ for some $\eta \in ]0,1[$. However, we have room for 
that if we start with a slightly large parameter $C_1>0$ to begin with and then restrict to a 
$C_1>0$ such that \eqref{lb17} and \eqref{gp2} hold.
\\
\par
Next, we define the maps 
\begin{equation}\label{lb18}
\begin{split}
\kappa :  \mathcal{H}_N \supset &\mathcal{Q}_{C_1N} \longrightarrow 
\kappa(\mathcal{Q}_{C_1N}) \subset \mathcal{H}_N \\
&Q\longmapsto \kappa(Q) \defeq \widehat{Q}  + \widehat{T} (Q),
\end{split}
\end{equation}
where we identify $\widehat{T} (Q)$ with its image in $\mathcal{H}_N$ under the natural 
inclusion map $\mathcal{H}_{|J|} \xhookrightarrow{}\mathcal{H}_{N}$, which has the left 
inverse 
\begin{equation}\label{lb19}
	\widetilde{\Pi}_0:~ \mathcal{H}_N \to \mathcal{H}_{|J|}: ~ 
	Q\mapsto \widetilde{\Pi}_0(Q)\defeq \Pi_0 Q \Pi_0 
\end{equation}
Moreover, we define the map $\Pi:\mathcal{H}_N \supset \mathcal{Q}_{C_1N}\to  
\mathcal{H}_{|J|}$ by 
\begin{equation}\label{lb20}
	 \Pi \defeq \widetilde{\Pi}_0\circ \kappa.
\end{equation}
In analogy with \eqref{gp1.0} we define the probability measure $\mu_J$ on $\mathcal{H}_{|J|}$ 
by
\begin{equation}\label{lb21}
	\mu_J (dQ) \defeq \pi^{-|J|^2} \e^{- \| Q\|_{\mathrm{HS}}^2} L(dQ).
\end{equation}
We will estimate the probability 
\begin{equation}\label{lb22}
	\mu_N\left(\log |\det E_{-+}^{\delta}(z,Q) |^2 \leq a \text{ and } Q\in\mathcal{Q}_{C_1N}\right).
\end{equation}
To begin, we strengthen \eqref{lb2.1} to 
\begin{equation}\label{lb23}
	0 < \delta \ll \frac{\alpha}{(C_1N)^3}.
\end{equation}
By \eqref{lb17}, \eqref{lb18}, we see that $\kappa$ is injective, since 
for $Q_1,Q_2 \in\mathcal{Q}_{C_1N}$ 
\begin{equation*}
	\begin{split}
		\| \kappa(Q_1) - \kappa(Q_2)\| 
		&\geq 
		\| Q_1- Q_2\|  - \int_0^1 \| d_Q \widehat{T}(tQ_1 + (1-t)Q_2)\| \cdot \| Q_1- Q_2\| dt \\
		& \geq \left(1 -  \mO\!\left(\frac{\delta C_1N}{\alpha}\right)\right)\| Q_1- Q_2\|. 
	\end{split}
\end{equation*}
Define the restricted measure 
\begin{equation}\label{lb23.1}
	(\mathbf{1}_{\mathcal{Q}_{C_1N}}\mu_N)(A)\defeq \mu_N (A\cap\mathcal{Q}_{C_1N} ), 
	\quad \forall A\in \mathcal{B}(\mathcal{H}_N),
\end{equation}
where $\mathcal{B}(\mathcal{H}_N)$ denotes the Borel $\sigma$-algebra of $\mathcal{H}_N$. 
In view of the discussion after \eqref{lb4}, the measure 
$\mathbf{1}_{\mathcal{Q}_{C_1N}}\mu_N$ is invariant under the change of 
orthonormal basis of $\mathcal{Q}_{C_1N}$. Thus, 
by \eqref{lb20}, \eqref{lb16}, the probability in \eqref{lb22} is equal to 
\begin{equation}\label{lb24}
	\begin{split}
		&(\mathbf{1}_{\mathcal{Q}_{C_1N}}\mu_N)\left[ \log \left| \prod_{k=1}^{|J|} (s_k^+s_k^-) \det \left( \widehat{E}_{-+}^0(z) 
		- \delta ( \Pi_0 \widehat{Q} \Pi_0 +
		\widehat{T}(Q) )\right) \right|^2 \leq a\right] \\
		&= 
		(\mathbf{1}_{\mathcal{Q}_{C_1N}}\mu_N)\left[ \log \left|
			\det ( \delta^{-1}\widehat{E}_{-+}^0(z) - 
			\Pi(Q))\right|^2 \leq b\right] \\
		&= 
		\Pi_*(\mathbf{1}_{\mathcal{Q}_{C_1N}}\mu_N)\left[\log|\det (\delta^{-1}\widehat{E}_{-+}^0(z)  - 
			Q') |^2 \leq b \right],
	\end{split}
\end{equation}
where by \eqref{lb16}, \eqref{lb6.1},
\begin{equation}\label{lb24.1}
\begin{split}
   b &= a - 2|J| \log\delta - 2 \sum_{j=1}^{|J|} \log( s_j^+ s_j^- ) \\
       &\leq a - 2|J| \log\delta + 4 |J| \log C. 
 \end{split}
\end{equation}
\par
Continuing, we will estimate the measure $\Pi_*(\mathbf{1}_{\mathcal{Q}_{C_1N}}\mu_N)$. We begin 
by studying the Jacobian of $\kappa$, \eqref{lb18}. By \eqref{lb17} and \eqref{lb23}, we see that the differential 
of $\widehat{T}$ is bounded with norm $\ll 1$. Moreover, since the rank of 
$d_Q\widehat{T}$ 
is bounded by $|J|^2$, it follows that $\|d_Q\widehat{T} \|_{\mathrm{tr}} \leq |J|^2 \|d_Q\widehat{T}\|$. 
Thus, by \eqref{lb17}
\begin{equation}\label{lb25}
	\begin{split}
		\det \frac{\partial \kappa}{\partial Q} 
		&= \det \left ( 1 + d_Q\widehat{T}\right) \\  
		&=  1 + \mO(\|d_Q\widehat {T}\|_{\mathrm{tr}}) \\ 
		& = 1 + \mO\!\left(\frac{\delta C_1N}{\alpha}\right),
	\end{split}
\end{equation}
where in the last line we used as well that $|J|$ is a constant independent of $N$. 
\par
Since $\kappa$ is a holomorphic map, it follows that 
\begin{equation}\label{lb26}
	\begin{split}
	L(d\kappa(Q)) & = \left| \det \frac{\partial \kappa}{\partial Q} \right|^2 L(dQ) \\ 
			       & = \left( 1 + \mO\!\left(\frac{\delta C_1N}{\alpha}\right) \right ) L(dQ). 
	\end{split}
\end{equation}
Next, we see by \eqref{lb18}, \eqref{lb15.2}, that for $Q\in \mathcal{Q}_{C_1N}$ 
\begin{equation*}
	\begin{split}
	\left| \|\kappa(Q)\|_{\mathrm{HS}}^2 -  \|Q\|_{\mathrm{HS}}^2 \right| 
	&= \left| \|\kappa(Q)\|_{\mathrm{HS}} -  \|Q\|_{\mathrm{HS}} \right| 
	( \|\kappa(Q)\|_{\mathrm{HS}} +  \|Q\|_{\mathrm{HS}}) \\
	&\leq \|\kappa(Q) -Q\|_{\mathrm{HS}}( \|\kappa(Q)\|_{\mathrm{HS}} +  \|Q\|_{\mathrm{HS}})  \\ 
	&= \mO\!\left(\frac{\delta (C_1N)^2}{\alpha}\right)\left( C_1N +  
	\mO\!\left(\frac{\delta (C_1N)^2}{\alpha}\right)\right) \\
	&= \mO\!\left(\frac{\delta (C_1N)^3}{\alpha}\right) \ll 1,
	\end{split}
\end{equation*}
which implies that on $\mathcal{Q}_{C_1N}$ 
\begin{equation}\label{lb27}
	\e^{-\|Q\|_{\HS}^2} = 
	\left(1 + \mO\!\left(\frac{\delta (C_1N)^3}{\alpha}\right)\right)\e^{-\|\kappa(Q)\|_{\HS}^2}.
\end{equation}
\eqref{lb26}, \eqref{lb27} imply that for any bounded continuous function 
$\varphi \in \mathcal{C}_b(\mathcal{H}_N;\R_+)$ with values in $\R_+$, 
\begin{equation*}
\begin{split}
	\int \varphi\,\kappa_*(\mathbf{1}_{\mathcal{Q}_{C_1N}}\mu_N) 
	&= \int _{\mathcal{Q}_{C_1N}}\varphi(\kappa(Q)) \mu_N(dQ) \\
	& = \left(1 + \mO\!\left(\frac{\delta (C_1N)^3}{\alpha}\right)\right)
	       \int _{\mathcal{Q}_{C_1N}}\varphi(\kappa(Q)) \e^{-\|\kappa(Q)\|_{\HS}^2}
	      \frac{ L(d\kappa(Q))}{\pi^{N^2}} \\
	&= \left(1 + \mO\!\left(\frac{\delta (C_1N)^3}{\alpha}\right)\right)
	       \int _{\kappa(\mathcal{Q}_{C_1N})}
	       \varphi(\widetilde{Q}) \e^{-\|\widetilde{Q}\|_{\HS}^2}
	       \frac{ L(d\widetilde{Q})}{\pi^{N^2}}.
\end{split}
\end{equation*}
Thus, 
\begin{equation}\label{lb28}
	\kappa_*(\mathbf{1}_{\mathcal{Q}_{C_1N}}\mu_N) =
	\left(1 + \mO\!\left(\frac{\delta (C_1N)^3}{\alpha}\right)\right)
	\mathbf{1}_{\kappa(\mathcal{Q}_{C_1N})}\mu_N.
\end{equation}
This, together with \eqref{lb20}, implies that for any 
$\varphi \in \mathcal{C}_b(\mathcal{H}_{|J|};\R_+)$ 
\begin{equation*}
\begin{split}
	\Pi_*(\mathbf{1}_{\mathcal{Q}_{C_1N}}\mu_N)(\varphi)
	&= \int (\varphi\circ \widetilde{\Pi}_0) \, \kappa_*(\mathbf{1}_{\mathcal{Q}_{C_1N}}\mu_N)  \\
	& = \left(1 + \mO\!\left(\frac{\delta (C_1N)^3}{\alpha}\right)\right)
	       \int \varphi\circ \widetilde{\Pi}_0 \, \mathbf{1}_{\kappa(\mathcal{Q}_{C_1N})}\mu_N \\
	&\leq \left(1 + \mO\!\left(\frac{\delta (C_1N)^3}{\alpha}\right)\right)
	       \int \varphi\circ \widetilde{\Pi}_0\, \mu_N \\
        &\leq \left(1 + \mO\!\left(\frac{\delta (C_1N)^3}{\alpha}\right)\right)
	       \int \varphi(Q')\, \mu_{J}(dQ'),
\end{split}
\end{equation*}
where in the last line we used that $(\widetilde{\Pi}_0)_*\mu_N = \mu_{J}$. Hence, 
by \eqref{lb24} and a density argument, we deduce that  the probability in 
\eqref{lb22} is
\begin{equation}\label{lb29}
\begin{split}
	 \leq \left(1 + \mO\!\left(\frac{\delta (C_1N)^3}{\alpha}\right)\right)
			   \mu_{J}\left[\log|\det(\delta^{-1}\widehat{E}_{-+}^0(z)  - Q')|^2 
			   \leq b\right]. 
\end{split}
\end{equation}
The right hand side can be estimated by \cite[Proposition 7.3]{HaSj08}.
\begin{prop}\label{prop:LB1}
Let $\N \ni M\geq 1$, let $\mu_M$ be the Gaussian measure on $\mathcal{H}_M$ defined 
in \eqref{gp1.0}. Then, there exist constants $\widetilde{C},C'>0$ such that for any 
fixed (deterministic) matrix $D\in\mathcal{H}_M$ 
\begin{equation*}
\begin{split}
	\mu_M(\log|\det(D + Q)|^2 \leq b ) 
	& \leq \mu_M(\log|\det Q)|^2 \leq b ) \\
	&\leq 
	\widetilde{C} \exp \left[ 
	-\frac{1}{2}\left( 
	C' + \left( M-\frac{1}{2}\right) \ln M - 2M -b
	\right)
	\right],
\end{split}
\end{equation*}
when $b \leq C' + \left( M+\frac{1}{2}\right) \ln M - 2M$.
\end{prop}
Combining, \eqref{lb29}, \eqref{lb22}, \eqref{lb24}, \eqref{lb24.1} and \eqref{lb6.1} with 
Proposition \ref{prop:LB1}, we deduce that there exist constants 
$\widetilde{C},C'>0$ such that 
\begin{equation*}
\begin{split}
	&\mu_N(\{\log |\det E_{-+}^{\delta}(z,Q) |^2 \leq a \} \cap \mathcal{Q}_{C_1N})\\
    &\leq 
	\widetilde{C} \exp \left[ 
	-\frac{1}{2}\left( 
	C' + \left( |J|-\frac{1}{2}\right) \ln |J| - 2|J| - b
	\right)
	\right]	\\
	&\leq 
	\widetilde{C} \exp \left[ 
	-\frac{1}{2}\left( 
	C' + \left( |J|-\frac{1}{2}\right) \ln |J| - 2|J| - a + 2|J| \log\delta  - 4|J| \log C
	\right)
	\right]
\end{split}
\end{equation*}
when $b  \leq C' + \left( |J|+\frac{1}{2}\right) \ln |J| - 2|J|$ and thus, by 
\eqref{lb24.1}, when 
\begin{equation*}
a \leq C' + \left( |J|+\frac{1}{2}\right) \ln |J| - 2|J| + 2|J| \log\delta  - 4|J| \log C.
\end{equation*}
Here, the constants $\widetilde{C},C'$ only depend on $J$ and the constant $C$ is 
given by the lower bounds in \eqref{lb6.1} which are uniform in $z\in\Omega$. 
Setting 
\begin{equation*}
\begin{split}
	&C_0 = C' + \left( |J|+\frac{1}{2}\right) \ln |J| - 2|J| - 4|J| \log C, \\
	&a = -t,
\end{split}
\end{equation*}
we conclude, by absorbing the factor $\e^{-\frac{1}{2}(C_0-\log |J|)}$ into the constant $\widetilde{C}$, that 
\begin{equation}\label{lb29.1}
	\mu_N(\{\log |\det E_{-+}^{\delta}(z,Q) |^2 \leq - t  \} \cap \mathcal{Q}_{C_1N})
	\leq 
	\widetilde{C} \exp \left[ -\frac{1}{2}t - |J|\log \delta
	\right]
\end{equation}
when $t \geq C_0 - 2|J|\log \delta$. Finally, since 
\begin{equation*}
\mathds{P}[A^c\cap B] = \mathds{P}[ B] -\mathds{P}[A\cap B],
\end{equation*}
where $A^c$ denotes the complement of the measurable set $A$, we obtain, by combining 
\eqref{lb29.1}and \eqref{gp2},  
\begin{prop}\label{prop:LB}
Let $\kappa \geq 1$, let $\Omega\Subset\C$ be a compact set, let $C>0$ and let 
$C_1>0$ be such that \eqref{gp2} holds. 
Then, there exist constants $C_0\in\R$ and $C_2>0$, such that for any $z\in\Omega$, 
with
\begin{equation*}
	\alpha(z;N)=\dist(z,p(\widehat{S}_{\widetilde{N}})) \geq \frac{1}{CN^{\kappa}}, 
\end{equation*}
we have that 
\begin{equation*}
	\mathds{P} \left [ \log|\det E_{-+}^{\delta}(z,Q) |^2 \geq 
		-t \text{ and } \| Q\|_{\mathrm{HS}} \leq C_1 N
	\right] \geq 1 - \e^{-N^2} - C_2\,\delta ^{-|J|} \e^{ -t/2}, 
\end{equation*}
when
\begin{equation*}
t \geq C_0 - 2|J|\log \delta 
\end{equation*}
and 
\begin{equation*}
	0 < \delta \ll \frac{\alpha}{(C_1N)^3}.
\end{equation*}
\end{prop}
\section{Counting eigenvalues}
In this section we count the eigenvalues of the perturbed operator 
\begin{equation}\label{c1}
	P_N^{\delta} = P_N^0 + \delta Q_{\omega},
\end{equation}
near the curve $p(S^1)$, see also \eqref{gp1}. Recall from \eqref{gu5.1} that $P_N^0= P_{I_N}$,  
see also\eqref{gu7}. Similarly, we have $P_N^{\delta}= P^{\delta}_{I_N}$ as in Proposition \ref{gp:prop1}. 
\par
Until further notice, we will work in the restricted probability space where 
\eqref{gp1.1} holds (see also \eqref{gp2}) and work under the assumptions that 
\begin{equation}\label{c3}
	0 < \delta \ll \frac{\alpha}{N^3}, \quad  \frac{1}{CN} \leq \alpha \leq \mO(1),
\end{equation}
for some sufficiently large constant $C>0$ to be determined later on, see also \eqref{gp14}, \eqref{lb23}. 
Here $\alpha$ is as in \eqref{lb1}.
\\
\par
Counting the number of eigenvalues of $P_{I_N}^{\delta}$ in some 
domain $\Omega \Subset \C$ is equivalent to counting the number of 
zeros of the holomorphic function $u(z;N) = \det (P_{I_N}^{\delta} -z)$ 
in $\Omega$. 
The Shur complement formula 
and Proposition \ref{gp:prop1} imply that, away from $\spec(P_{S_{\widetilde{N}}})$, 
$P_{I_N}^{\delta} -z$ is invertible 
if and only if $E_{-+}^{N,\delta}(z) $ is invertible, and that
\begin{equation}\label{gp18}
	\log|\det (P_{I_N}^{\delta} -z) |= 
	\log|\det \mathcal{P}_N^{\delta}(z) |+
	\log|\det E_{-+}^{N,\delta}(z) |. 
\end{equation}
\subsection{Counting zeros of holomorphic functions of exponential growth}
\label{sec:c1}
We recall  Theorem $1.1$ in 
\cite{Sj09b}, in a form somewhat adapted to our formalism:
\\
\par
\paragraph{\textit{1) Domains with associated Lipschitz weight}}
Let $N\geq1$ be a large parameter, and let $\Omega \Subset \C$ be an open 
simply connected set with Lipschitz boundary $\omega =\partial\Omega$ which may 
depend on $N$. More precisely, we assume that $\partial\Omega$ is Lipschitz 
with an associated Lipschitz weight $r:\omega \to ]0,+\infty[$, which is 
a Lipschitz function of modulus $\leq 1/2$, in the following way : 
\par
There exists a constant $C_0>0$ such that for every $x\in\omega$ 
there exist new affine coordinates $\widetilde{y} = (\widetilde{y}_1,\widetilde{y}_2)$ 
of the form $\widetilde{y} = U(y-x)$, $y\in \C \simeq \R^2$ being the old coordinates, 
where $U=U_x$ is orthogonal, such that the intersection of $\Omega$ and the 
rectangle $R_x := \{ y\in\C; |\widetilde{y}_1| < r(x), |\widetilde{y}_2| < C_0r(x)\}$ 
takes the form 
\begin{equation}\label{eq.LD1}
	\{ y \in R_x; ~\widetilde{y}_2 > f_x(\widetilde{y}_1), |\widetilde{y}_1| < r(x)\},
\end{equation}
where $f_x(\widetilde{y}_1)$ is Lipschitz on $[-r(x),r(x)]$, with Lipschitz modulus 
$\leq C_0$. 
\begin{remark}\label{rem:LD}
 Notice that \eqref{eq.LD1} remains valid if we shrink the weight function $r$.
\end{remark}
\paragraph{\textit{2) Thickening of the boundary and choice of points}}
Define 
\begin{equation*}
	\widetilde{\omega}_r = \bigcup_{x\in \omega}D(x,r(x))
\end{equation*}
and let $z_j^0 \in \omega$, $j\in \Z/ M\Z$, with $M\in \N$ which may depend on $N$, 
be distributed along the boundary in the positively oriented sense such that 
\begin{equation*}
	r(z_j^0)/4 \leq |z_{j+1}^0 - z_j^0| \leq r(z_j^0)/2.
\end{equation*}
\begin{theo}[Theorem 1.1 in \cite{Sj09b}]\label{thm:Count}
	Let $C_0>0$ be as in 1) above. 
	There exists a constant $C_1>0$, depending only on $C_0$, 
	such that if $z_j \in D(z_j^0, r(z_j^0)/(2C_1))$ we have the 
	following :
	\par
	Let $N\geq 1$ and let $\phi $ be a continuous subharmonic function on 
	$\widetilde{\omega}_r$ with a distributional extension to 
	$\Omega\cup \widetilde{\omega}_r$, denoted by the same symbol. Then, 
	there exists a constant $C_2>0$ such that if $u$ is a holomorphic function 
	on $\Omega\cup \widetilde{\omega}_r$ satisfying 
	\begin{equation}\label{th:UB}
		\log |u | \ \leq N \phi \hbox{ on } \widetilde{\omega}_r,
	\end{equation}
	\begin{equation}\label{th:LB}
		\log |u (z_j)| \ \geq N (\phi(z_j) - \varepsilon_j), \hbox{ for } 
		j=1, \dots , M,
	\end{equation}
	where $\varepsilon_j \geq 0$, then the number of zeros of $u$ in $\Omega$
	satisfies
	\begin{equation*}
	\begin{split}
		\bigg|\#(u^{-1}(0)\cap \Omega) &- \frac{N}{2\pi} \mu(\Omega)\bigg| \\
		&\leq C_2 N \left(
		\mu(\widetilde{\omega}_r) + \sum_{j=1}^M
		\left( \varepsilon_j 
		+ \int_{D\!\left(z_j,\frac{r(z_j)}{4C_1}\right)}\left| \log \frac{|w-z_j|}{r(z_j)} \right |
		\mu(dw)
		\right)
		\right).
	\end{split}
	\end{equation*}
	Here $\mu \defeq \Delta\phi \in \mathcal{D}'(\Omega\cup \widetilde{\omega}_r)$ 
	is a positive measure on $\widetilde{\omega}_r$ so that $\mu(\Omega)$ and 
	$\mu(\widetilde{\omega}_r)$ are well-defined. Moreover, the constant 
	$C_2>0$ only depends on $C_0$. 
\end{theo}
\subsection{Upper bound on $\log|\det (P_{I_N}^{\delta} -z) |$}\label{sec:c2}
Recall from \eqref{g27.5}, that 
$\#(p(\widehat{S}_{\widetilde{N}}))=\widetilde{N}$ where 
$\widetilde{N} = N + N_- + N_+$. Then, define the subharmonic function 
\begin{equation}\label{gp19}
	\phi(z)\defeq \phi(z;N)
	\defeq \frac{1}{N}\sum_{\lambda \in p(\widehat{S}_{\widetilde{N}})}
	\log|\lambda - z|. 
\end{equation}
Applying \eqref{gp9}, \eqref{gp16}, \eqref{c3} to \eqref{gp18} we can 
express the contribution from the perturbed Grushin problem in \eqref{gp18} 
by the function $\phi$ and a small error term, i.e.  
\begin{equation}\label{gp20}
\begin{split}
	\log|\det (P_{I_N}^{\delta} -z) | &= \log|\det \mathcal{P}_N^{0}(z) |
	+ \mO(\delta\|Q_{\omega}\|_{\mathrm{tr}} \|E^{N,0}\| ) + 
	\log|\det E_{-+}^{N,\delta}(z) | \\
	& =
	N\left(
		\phi(z) + \frac{ \log|\det E_{-+}^{N,\delta}(z) |}{N} + 
		\mO\!\left(\frac{\delta\|Q_{\omega}\|_{\HS}}{N^{1/2}\alpha}\right)
	\right).
\end{split}
\end{equation}
In the last line we used that 
$\|Q_{\omega}\|_{\mathrm{tr}} \leq N^{1/2}\|Q_{\omega}\|_{\HS}$.
\par
By \eqref{c3}, \eqref{gp2} we have that 
$\alpha^{-1}\delta \|Q_{\omega}\|_{HS}  \ll N^{-2}$. Recall that the dimension of 
the matrix $E_{-+}^{\delta}$ is $|J|=N_++N_-$. Therefore, using \eqref{c3}, \eqref{gp16} and 
Proposition \ref{gp:prop1}, we can bound \eqref{gp20} from above and get 
\begin{equation}\label{gp20.1}
\log|\det (P_{I_N}^{\delta} -z) | 
\leq  
	N\left(
		\phi(z) + 
		\mO(N^{-1}| \log \alpha | )+ 
		\mO(N^{-5/2})
	\right).
\end{equation}
In conclusion, assuming \eqref{c3}, we have that 
\begin{equation}\label{c4}
	\log|\det (P_{I_N}^{\delta} -z) | \leq  N \psi(z;N)
\end{equation}
with probability $\geq 1 - \e^{-N^2}$. Here,  
\begin{equation}\label{c5}
	\psi(z;N) \defeq \phi(z) +  \frac{ C\log N }{N}, 
\end{equation}
for some sufficiently large constant $C>0$. 
\subsection{Lower bound on $\log|\det (P_{I_N}^{\delta} -z) |$}\label{sec:c3}
Fix a $\varepsilon_0 \in ]0,1[$. By \eqref{c3} and Proposition \ref{prop:LB} we have for any 
$z_0$, satisfying 
\begin{equation*}
	\alpha(z_0;N) \geq \frac{1}{CN},
\end{equation*}
that 
\begin{equation}\label{c6}
	\mathds{P} \left [\log |\det E_{-+}^{\delta}(z_0,Q) |^2 \geq 
		-N^{\varepsilon_0} \text{ and } \| Q\|_{\mathrm{HS}} \leq C_1 N
	\right] \geq 1 - \e^{-N^2} - C_2\delta^{-|J|} \e^{-\frac{1}{2}N^{\varepsilon_0}}, 
\end{equation}
for
\begin{equation}\label{c6.1}
	\exp \left[\frac{C_0}{2|J|} - \frac{N^{\varepsilon_0}}{2|J|} \right] \leq \delta 
	\ll \frac{\alpha(z_0;N)}{N^3}.
\end{equation}
Thus, assuming \eqref{c6.1} and combining  \eqref{c6}, \eqref{gp20}, \eqref{c3} and \eqref{c5}, 
we get that $\|Q\|_{\HS}\leq C_1N$ and  
\begin{equation}\label{c7}
	\log|\det (P_{I_N}^{\delta} -z_0)| \geq  
	N\left(\psi(z_0;N) - CN^{\varepsilon_0-1}
	\right)
\end{equation}
hold with probability 
\begin{equation}\label{c8}
	 \geq 1 - \e^{-N^2} - C_2\delta^{-|J|} \e^{-\frac{1}{2}N^{\varepsilon_0}}.
\end{equation}

\subsection{Counting eigenvalues in a fixed smooth domain}\label{sec:c4}
Let $\Omega \Subset \C$ be an open simply connected set with smooth 
boundary $\partial\Omega$ which is independent of $N$. Moreover, 
suppose that ($\Omega$\ref{O1})--($\Omega$\ref{O3}) hold.
\\ 
\par
To estimate the number of zeros of $\det(P_{I_N}-z)$, see \eqref{gp18}, in $\Omega$, 
we will apply Theorem \ref{thm:Count}. 
The boundary $\partial\Omega$ is uniformly Lipschitz at scale 
\begin{equation}\label{fd.4}
	r(x) \defeq 
	\frac{1}{C}\left(\dist(x,p(S^1)) + \frac{1}{N} \right), 
	\quad x\in \partial\Omega
\end{equation}
which is Lipschitz of modulus $\leq 1/2$. Here, $C>0$ is chosen sufficiently 
large, and we will potentially increase it later on.
\par
Due to the singularities of $\psi$ at $p(\widehat{S}_{\widetilde{N}})$, see \eqref{c5}, \eqref{gp19}, 
we cannot in general assure that the weight function $\psi$ \eqref{c5} be continuous in
\begin{equation*}
	\bigcup_{x\in\partial\Omega} D(x,r(x)).
\end{equation*}
To remedy this problem we will consider two $N$-dependent perturbations of 
the boundary $\partial\Omega$: let $z_0\in p(S^1)\cap \partial\Omega$ and 
pass to new affine coordinates $\widetilde{y}\in \R^2\simeq \C$ (as in Section \ref{sec:c1}) 
so that the boundary $\partial\Omega$ 
is given by the graph of the smooth function $f_{z_0}$ near $0$, with derivative bounded 
by $C_0>0$. For $C'>1$ and $N>0$ sufficiently large, the intersection of $\partial\Omega$ 
with the rectangle
\begin{equation}\label{fd.4.0}
	R_{z_0}(N) \defeq \{ y \in\C\simeq \R^2; |\widetilde{y}_1| \leq 1/(C'N), ~ 
	|\widetilde{y}_2| \leq 2C_0/(C'N) \}
\end{equation}
takes the form 
\begin{equation*}
	\{ y \in\C\simeq \R^2; |\widetilde{y}_1| \leq 1/(C'N), ~ 
	\widetilde{y}_2 > f_{z_0}(\widetilde{y}_1) \}.
\end{equation*}
Here, $y\in\C\simeq \R^2$ denote the old coordinates and $\widetilde{y}\in\C\simeq \R^2$ 
denote the new ones. 
\par
Next, define the continuous function $\widetilde{ \chi}$, supported in $[-1,1]$  and of Lipschitz modulus $2$, by 
\begin{equation*}
	\widetilde{ \chi}(x) = 
	\begin{cases}
		2(x+1),~-1 \leq x < - 1/2,\\ 
		1, ~ |x| \leq 1/2,\\
		1 - 2(x-1/2),~ 1/2 < x \leq1, \\
	\end{cases}
\end{equation*}
 and set 
\begin{equation*}
	 \chi(\widetilde{y}_1)\defeq \chi(\widetilde{y}_1;N) 
	 \defeq \frac{C_0}{4C'N}\widetilde{ \chi} (C'N \widetilde{y}_1 ).
\end{equation*}
Moreover, we define for $\eta_{\pm}\in [0,1]$ 
\begin{equation*}
	f_{z_0}^{\eta_\pm} ( \widetilde{y}_1) \defeq f_{z_0}( \widetilde{y}_1) \pm
	\eta_{\pm} \chi( \widetilde{y}_1).
\end{equation*}
Since $f_{z_0}$ has Lipschitz modulus $\leq C_0$, if follows that $f_{z_0}^{\eta_{\pm}}$ 
has Lipschitz modulus $\leq 3C_0/2$, for $N>0$ sufficiently large. 
\par 
By Proposition \ref{prop:EigDist}, it follows that the number of eigenvalues of 
$P_{S_{\widetilde{N}}}$ contained in $R_{z_0}(N)$ is bounded by a constant 
depending only $p$, $C'$ and $C_0$. Since the are only finitely many points to avoid, 
there exist $\eta_{\pm} \in [0,1]$ such that 
\begin{equation}\label{fd.4.5}
\{ y\in \R^2\simeq \C;   | \widetilde{y}_1| \leq 1/C'N,~ \widetilde{y}_2 = f_{z_0}^{\eta_\pm}(\widetilde{y}_1) \} 
 \cap \big(\spec (P_{S_{\widetilde{N}}})\cap R_{z_0}(N)\big) 
    =\emptyset.
\end{equation}
For $C,C', \widetilde{C}>0$ large enough we can arrange that 
\begin{equation}\label{fd.4.1}
\bigg(\bigcup_{\substack{y\in \R^2 \hbox{ s.t. } (\widetilde{y}_1,f_{z_0}^{\eta,\pm}(\widetilde{y}_1)) \\ 
    | \widetilde{y}_1| \leq 1/C'N }} \overline{D(y,r(y))} \bigg)\cap 
    \bigcup_{\lambda \in \spec (P_{S_{\widetilde{N}}})} \overline{D(\lambda, 1/(\widetilde{C}N) )}\cap R_{z_0}(N)
    =\emptyset.
\end{equation}
We perform these 
two deformations of $\partial\Omega$ near every point $z_0 \in p(S^1)\cap \partial\Omega$, 
pick $C>0$ in \eqref{fd.4} at least as large as the maximum over all constants $C$ so that 
\eqref{fd.4.1} holds, and call the resulting deformed sets 
\begin{equation}\label{fd.4.2}
	\Omega_{\pm} \hbox{ with boundary } \partial\Omega_{\pm}.
\end{equation} 
Here, we always take the local deformation $f_{z_0}^{\eta_+}$  for $\Omega_+$, 
and $f_{z_0}^{\eta_-}$ for $\Omega_-$. Notice that since 
\begin{equation*}
	f_{z_0}^{\eta_-} ( \widetilde{y}_1) \leq f_{z_0}( \widetilde{y}_1) \leq 
	f_{z_0}^{\eta_+} ( \widetilde{y}_1), \quad |\widetilde{y}_1| \leq 1/(C'N),
\end{equation*}
we have 
\begin{equation}\label{fd.4.4}
	\Omega_+ \subset \Omega \subset \Omega_-,
\end{equation}
where we do not denote the $N$ dependence explicitly. 
\par
By \eqref{fd.4.1} , ($\Omega$\ref{O1}) and ($\Omega$\ref{O3}),  
there exists a $C>0$ such that 
\begin{equation}\label{fd.4.3}
	\dist \left(\bigcup_{x\in\partial\Omega_{\pm}} \overline{D(x,r(x))}
	,\spec (P_{S_{\widetilde{N}}})\right)\geq \frac{1}{CN},
\end{equation}
which also determines the constant $C>0$ in \eqref{c3}. 
Next, choose points $z_j^{0,\pm} \in \partial \Omega_{\pm}$, $j\in \Z/M\Z$, such that 
\begin{equation}\label{fd.5}
	\partial\Omega_{\pm} \subset \bigcup_{j\in \Z/M\Z} D(z_j^0,r_j^{\pm}/2), 
	\hbox{ and }r_j^{\pm}/4 \leq|z^{0,\pm}_{j+1}-z_j^{0,\pm}| \leq r_j^{\pm}/2, 
\end{equation}
where $r_j^{\pm} = r(z_j^{0,\pm})$. 
\begin{lemma}\label{lem:c4.1} Let $M$ be as in \eqref{fd.5}. Then, 
	\begin{equation*}
		M = \mO(\log N).
	\end{equation*}
\end{lemma}
We will postpone the proof of Lemma \ref{lem:c4.1} to the end of this section 
and carry on with the proof of our main result. 
\\
\par
First, notice that \eqref{c4} holds in $ \bigcup_{j=1}^M D(z_j^0,r_j)$ with probability 
$\geq 1- \e^{-N^2}$. By  \eqref{fd.4.3}, it follows that the weight function 
$\psi(z;N)$ \eqref{c5} is continuous on $\bigcup_{x\in\partial\Omega_{\pm}} \overline{D(x,r(x))}$. 
Moreover, by \eqref{fd.4.3}, we have that for any $z_j\in D(z_j^0,r_j/2)$
\begin{equation}\label{c17.1}
	 \alpha(z_j;N) \geq \frac{1}{CN}, 
\end{equation}
and so it follows that \eqref{c7} holds with probability \eqref{c8}, assuming \eqref{c6.1}. Hence, 
using Lemma \ref{lem:c4.1}, we have that \eqref{c7} holds for $z_1^0,\dots,z_M^0$ with probability 
 \begin{equation}\label{c17}
	 \geq 1 - 
	 \mO(\log N)
	 \left(\e^{-N^2} +C_2\delta^{-|J|} \e^{-\frac{1}{2}N^{\varepsilon_0}}\right).
\end{equation}
In view of \eqref{c7}, we can pick $\varepsilon_j = CN^{\varepsilon_0-1} $ in 
Theorem \ref{thm:Count}, so using Lemma \ref{lem:c4.1}, we get 
 \begin{equation}\label{c18}
 \begin{split}
	&\bigg| 
	\#(\spec(P^{\delta}_N)\cap\Omega_{\pm})
	- \frac{N}{2\pi} \int_{ \Omega_{\pm}} \Delta\phi(z)L(dz) \bigg| \\
	&\leq 
	\mO(N)\!\left(  N^{\varepsilon_0-1}\log N 
	+\mu\!\left(\bigcup_{x\in\partial\Omega_{\pm}} D(x,r(x)) \right) 
	+ \sum_{j=1}^M
	 \int_{D\!\left(z_j^0,\frac{r(z_j^0)}{4C_1}\right)}\left|
	  \log \frac{|w-z_j^0|}{r(z_j^0)} \right |
		\mu(dw)
	\right),
\end{split}
\end{equation}
with probability \eqref{c17}, where we used as well that $\Delta\psi(z;N) = \Delta\phi(z)$, 
see \eqref{c5}. Moreover,  since $\Delta_z \log |z-w| = 2\pi \delta_w$, we have 
\begin{equation}\label{c18.1}
	\Delta\phi  = \mu  = \frac{2\pi}{N} \sum_{\lambda \in p(\widehat{S}_{\widetilde{N}})}
	\delta_{\lambda} \hbox{ in } \mathcal{D}'(\C).
\end{equation}
\par
The integral in the first line is up to an error of order $\mO(1)$ the number of eigenvalues 
of $P_{S_{\widetilde{N}}}$ contained in $\Omega\cap p(S^1)$. Hence,  
by \eqref{gp19} and \eqref{sp6}, 
 \begin{equation}\label{c19}
 \frac{N}{2\pi} \int_{ \Omega_{\pm}} \Delta\phi(z)L(dz) 
 = \frac{N}{2\pi} \int_{p^{-1}(\Omega\cap p(S^1))}L_{S^1}(d\theta) +\mO(1).
 \end{equation}
By \eqref{fd.4.3}
\begin{equation}\label{c20}
	\mu\!\left(\bigcup_{x\in\partial\Omega} D(x,r(x)) \right) 
	=
	0
\end{equation}
Similarly, the discs $\overline{D(z_j^0,r(z_j^0)/2)}$ do not contain any 
eigenvalues of $\mathcal{P}_N^0$. Thus, 
\begin{equation}\label{c21}
	\sum_{j=1}^M
	  \int_{D\!\left(z_j^0,\frac{r(z_j^0)}{4C_1}\right)}\left|
	  \log \frac{|w-z_j^0|}{\widetilde{r}(z_j^0)} \right |
	\mu(dw)
	=
	0
\end{equation}
Finally, from \eqref{fd.4.4}, it follows that 
\begin{equation}\label{c21.1}
	\#(\spec(P^{\delta}_N)\cap\Omega_{+})
	\leq 
	\#(\spec(P^{\delta}_N)\cap\Omega)
	\leq 
	\#(\spec(P^{\delta}_N)\cap\Omega_{-}).
\end{equation}
Combining \eqref{c18}, \eqref{c19}, \eqref{c20}, \eqref{c21} and \eqref{c21.1} we get that 
 \begin{equation}\label{c21.2}
	\bigg| 
	\#(\spec(P^{\delta}_N)\cap \Omega )
	-  \frac{N}{2\pi} \int_{p^{-1}(\Omega\cap p(S^1))}L_{S^1}(d\theta)  \bigg| 
	\leq 
	\mO( N^{\varepsilon_0}\log N ) .
\end{equation}
with probability \eqref{c17}, provided \eqref{c6.1} holds. This completes the proof of Theorem 
\ref{thm:t1}.
\begin{proof}[Proof of Lemma \ref{lem:c4.1}] 
\emph{1.} The perturbed boundaries $\partial\Omega_{\pm}$ \eqref{fd.4.2} coincide with 
$\partial\Omega$ outside the rectangles \eqref{fd.4.0}. Recall from ($\Omega$\ref{O1}) 
that there are only finitely many such rectangles. The number of discs of radius $r^{\pm}_j$ 
\eqref{fd.5} needed to cover $\partial\Omega_{\pm}$, as in \eqref{fd.5}, inside these rectangles 
is by \eqref{fd.4} of order 
\begin{equation}\label{fd.6.6}
	\mO(1).
\end{equation}
It remains to estimate the number of discs needed to cover $\partial\Omega$ outside these 
rectangles, which differs from order of the number of discs needed to cover the unperturbed 
$\partial\Omega$ by $\mO(1)$. Hence, it is sufficient to estimate the number of discs needed 
to cover $\partial\Omega$. 
\\
\par
\emph{2.} Since $\Omega$ is relatively 
compact and intersects with $p(S^1)$ at most finitely many points, we see that 
for any fixed constant $C>1$ the number of discs needed to cover 
$ \partial\Omega\cap \{z\in \C;\dist(z,p(S^1))\geq 1/C\}$, 
is of order
\begin{equation}\label{fd.6.7}
	\mO(1).
\end{equation}
\par
\emph{3.}
It remains to estimate the number of discs needed to cover $\partial\Omega$ 
inside $ \{z\in \C;\dist(z,p(S^1))\leq 1/C\}$. By assumption ($\Omega$\ref{O1}) 
and the fact that $\Omega$ 
is relatively compact 
we see that for any $\varepsilon>0 $ there exists $ \delta >0$ such that for any 
$x\in\partial\Omega$ 
 \begin{equation}\label{fd.6.1}
	\dist(x,p(S^1))  < \delta \quad \Longrightarrow
	 \min\limits_{z_0 \in p(S^1)\cap \partial\Omega}
	 \dist(x,z_0) < \varepsilon. 
\end{equation}
Hence, for any fixed $C'>0$,  we have for $C>0$ sufficiently large
 \begin{equation*}
	\partial\Omega \cap \{z\in \C;\dist(z,p(S^1))\leq 1/C\} \subset 
	\bigcup_{z_0 \in p(S^1)\cap \partial\Omega}D(z_0,1/C').
\end{equation*}
By ($\Omega$\ref{O1}), may restrict our attention to one $z_0\in \partial\Omega\cap p(S^1)$ 
and 
\begin{equation}\label{fd.6.1.2}
\beta =\partial\Omega \cap \{z\in \C;\dist(z,p(S^1))\leq 1/C\}\cap D(z_0,1/C').
\end{equation}
For $x,y\in \beta$ let $\dist_{\beta}(x,y)$ denote the length of the curve in $\beta$ 
with endpoints $x$ and $y$. By the transversality assumption ($\Omega$\ref{O3}), 
we see that for $C>0$ sufficiently large
\begin{equation}\label{fd.11}
	\dist_{\beta}(x,z_0) \asymp \dist(x,p(S^1)), \quad x \in \beta, 
\end{equation}
and 
\begin{equation}\label{fd.11.1}
	\dist_{\beta}(x,y) \asymp |x-y|, \quad x,y \in\beta. 
\end{equation}
\par
\emph{4.} 
Notice that $M_{\beta}$, the number of discs $D(z_i^0,r_i/2)$ needed to cover $\beta$, 
as in \eqref{fd.5}, increases when decreasing the scale $r$ \eqref{fd.4}. Using \eqref{fd.11} 
and by possibly increasing $C>0$ in \eqref{fd.4}, we shrink $r$ to the new scale 
\begin{equation}\label{fd.12}
	r(x) = 
	\frac{1}{C}\left(\dist_{\beta}(x,z_0) + \frac{1}{N} \right), 
	\quad x\in \beta, 
\end{equation}
denoted by the same letter. Set
\begin{equation}\label{fd.13}
	d_j \defeq \dist_{\beta}(z_j^0,z_0), \quad 1\leq j \leq M_{\beta},
\end{equation}
and let $j_1$ be the smallest index so that $d_{j_1} \geq  N^{-1}$. Notice that $j_1=\mO(1)$ 
and that $d_{j_1} \asymp N^{-1}$. By \eqref{fd.13}, \eqref{fd.11.1}, \eqref{fd.12} 
we have for $j> j_1$ 
\begin{equation}\label{fd.14}
\begin{split}
	d_j &=\dist_{\beta}(z_j^0,z_{j-1}^0)
	+\dist_{\beta}(z_{j-1}^0,z_0) \\
	      & \geq \frac{1}{C}|z_j^0-z_{j-1}^0| + d_{j-1} \\
	      &\geq (1+C^{-1})d_{j-1} \\
	      &\geq (1+C^{-1})^{j-j_1}d_{j_1},
\end{split}
\end{equation}
where the constant $C>0$ changes from the second to the third line. Similarly 
\begin{equation}\label{fd.15}
	d_j \leq (1+C)^{j-j_1}d_{j_1} .
\end{equation}
Thus, 
\begin{equation}\label{fd.16}
	(1+\widehat{C}^{-1})^{M_{\beta}-j_1}d_{j_1} \leq d_{M_{\beta}}  
	\leq (1+C)^{M_{\beta}-j_1}d_{j_1}.
\end{equation}
Using that the length of $\beta$ is $\asymp 1$, we get that $M_{\beta}\asymp \log N$ 
and therefore, by \eqref{fd.6.6}, \eqref{fd.6.7}, that 
\begin{equation*}
	M = \mO(\log N). \qedhere
\end{equation*}
\end{proof}
\subsection{Counting eigenvalues in thin $N$-dependent domains}\label{sec:c5}
In Section \ref{sec:c4} we saw that most eigenvalues of $P^{\delta}_N$ 
lie ``near'' the curve $p(S^1)$. Now we want to give a quantitative estimate 
on \emph{how close} these eigenvalues  are to the $p(S^1)$. For this purpose 
let $\Omega\Subset\C $ be an open simply connected set with smooth boundary 
$\partial \Omega$ which is independent of $N$ and satisfies 
($\Omega$\ref{O1})--($\Omega$\ref{O3}), as in Section \ref{sec:2.R2}. 
\par 
We consider an open simply connected $N$-dependent set $\Omega_N$, with a  
unifromly Lipschitz boundary $\partial \Omega_N$, which coincides with $\Omega$ 
in small tube around $p(S^1)$. More precisely, let
\begin{equation}\label{eq:s6.4.1}
		\frac{C}{N} \leq \tau \leq \mO(1), \quad C>1, 
\end{equation}
and suppose that 
\begin{equation}\label{eq:s6.4.2}
		\Omega_N \cap \{ z\in \C; \dist(z, p(S^1))< \tau\} = 
		\Omega \cap \{ z\in \C; \dist(z, p(S^1))< \tau\} ,
\end{equation}
and that $\partial\Omega_N$ is uniformly Lipschitz, as in Section \ref{sec:c1}, with 
weight function 
\begin{equation}\label{eq:s6.4.3}
		r(x) \defeq \frac{1}{C}\left(\dist(x,p(S^1)) + \frac{1}{N}\right), \quad 
		x\in \partial\Omega_N \cap \{ z\in \C; \dist(z, p(S^1))< \tau\},
\end{equation}
inside $\{ z\in \C; \dist(z, p(S^1))<\tau\}$ and with constant weight function 
\begin{equation}\label{eq:s6.4.4}
		r(x) \defeq \tau, \quad 
		x\in \partial\Omega_N \cap \{ z\in \C; \dist(z, p(S^1))\geq \tau\}
\end{equation}
outside. Let 
\begin{equation}\label{eq:s6.4.5}
	 	\ell(N) > 0
\end{equation}
be the length of $\partial\Omega_N \cap \{ z\in \C; \dist(z, p(S^1))\geq \tau\}$.
To prove Theorem \ref{thm:t2}, we can follow the proof of Theorem \ref{thm:t1} 
in Section \ref{sec:c4} with some modifications:
\\
\par
By \eqref{eq:s6.4.1} and \eqref{eq:s6.4.2}, we may perform the same perturbations of $\partial\Omega_N$ 
as for $\partial \Omega$ in \eqref{fd.4.0}--\eqref{fd.4.5} so that \eqref{fd.4.4} and 
\eqref{fd.4.3} hold for the perturbed sets 
\begin{equation}\label{eq:s6.4.6}
	 	\Omega_N^{\pm} \hbox{ with boundary } \partial\Omega_N^{\pm}. 
\end{equation}
Next, choose points $z_j^{0,\pm} \in \partial \Omega_N^{\pm}$, $j\in \Z/M\Z$, such that 
\begin{equation}\label{eq:s6.4.7}
	\partial\Omega_N^{\pm} \subset \bigcup_{j\in \Z/M\Z} D(z_j^{0,\pm},r_j^{\pm}/2), 
	\hbox{ and }r_j^{\pm}/4 \leq|z^{0,\pm}_{j+1}-z_j^{0,\pm}| \leq r_j^{\pm}/2, 
\end{equation}
where $r_j^{\pm} = r(z_j^{0,\pm})$. 
\begin{lemma}\label{lem:c5.1} Let $M$ be as in \eqref{eq:s6.4.7}. Then, 

	\begin{equation*}
		M =  \mO(\ell(N) \tau^{-1})+ \mO(\log (\tau N)).
	\end{equation*}
\end{lemma}
\begin{proof}
Following the exact same lines of \emph{Step 1, 3} and \emph{4} of the proof 
of Lemma \ref{lem:c4.1}, while keeping in mind \eqref{eq:s6.4.3} and that by \eqref{eq:s6.4.1}, \eqref{eq:s6.4.2} 
the length of $\partial\Omega_N \cap \{ z\in \C; \dist(z, p(S^1))\leq \tau\}$ is of order 
$\asymp \tau$, we see that the number of discs needed to cover 
$\partial\Omega_N \cap \{ z\in \C; \dist(z, p(S^1))\leq \tau\}$ is of order 
 \begin{equation}
 	\mO(\log (\tau N)).
 \end{equation}
By \eqref{eq:s6.4.5}, \eqref{eq:s6.4.4} we have that we have that the number of discs needed to 
cover $\partial\Omega_N \cap \{ z\in \C; \dist(z, p(S^1))\geq \tau\}$ is of order 
 \begin{equation}
 	\mO(\ell(N)\tau^{-1}). \qedhere 
 \end{equation}
\end{proof}
Since \eqref{fd.4.3} holds for $\partial\Omega_N^{\pm}$, the 
weight function $\psi(z;N)$ \eqref{c5} is continuous on 
	\begin{equation*}
		 \bigcup_{x\in\partial\Omega_N^{\pm}} \overline{D(x,r(x))},
	\end{equation*}
and that \eqref{c4} holds in $\bigcup_{j=1}^M D(z_j^0,r_j)$ \eqref{eq:s6.4.7} 
with probability $\geq 1 - \e^{-N^2}$. Moreover, since \eqref{fd.4.3} holds for 
$\partial\Omega_N^{\pm}$, we have that for any $z_j\in D(z_j^0,r_j/2)$
\begin{equation*}
	 \alpha(z_j;N) \geq \frac{1}{CN}, 
\end{equation*}
and it follows that \eqref{c7} holds with probability \eqref{c8}, assuming \eqref{c6.1}. Hence, 
using Lemma \ref{lem:c5.1}, we have that \eqref{c7} holds for $z_1^0,\dots,z_M^0$ with probability 
 \begin{equation}\label{eq:s6.4.9}
	 \geq 1 - 
	 \mO(M)
	 \left(\e^{-N^2} +C_2\delta^{-|J|} \e^{-\frac{1}{2}N^{\varepsilon_0}}\right).
\end{equation}
\par
In view of \eqref{c7}, we may set $\varepsilon_j = CN^{\varepsilon_0-1} $ in Theorem  \ref{thm:Count} and, by following the exact same arguments as above, from \eqref{c18} to \eqref{c21.1}, while keeping in mind Lemma \ref{lem:c5.1}, we obtain 
\begin{theo}\label{thm:t2}
Let $p$ be as in \eqref{int.0}, set $M = N_+ + N_-$ and let $P_N^{\delta}$ be as in 
\eqref{i0.1}. Let $\tau$ be as in \eqref{eq:s6.4.1} and let $\Omega_N\Subset \C$ be a 
relatively compact open simply connected set satisfying \eqref{eq:s6.4.2}--\eqref{eq:s6.4.5}.
 Pick a $\varepsilon_0\in]0,1[$. 
\par
There exists a constant $C>0$ such that for $N>1$ sufficiently large, if  \eqref{int.3} holds, 
\begin{equation*}
	C\e^{- N^{\varepsilon_0}/(2M)}\leq \delta \leq \frac{ N^{-4}}{C},
\end{equation*}
then,  
 \begin{equation}\label{t2.1}
	\bigg| 
	\#(\spec(P^{\delta}_N)\cap \Omega )
	-  \frac{N}{2\pi} \int_{p^{-1}(\Omega\cap p(S^1))}L_{S^1}(d\theta)  \bigg| 
	\leq 
	\mO(N^{\varepsilon_0}\ell(N) \tau^{-1}+N^{\varepsilon_0}\log (\tau N)).
\end{equation}
with probability 
 \begin{equation}\label{t2.2}
	 \geq 1 - 
	\mO(\ell(N) \tau^{-1}+\log (\tau N))
	 \left(\e^{-N^2} +C_2\delta^{-|J|} \e^{-\frac{1}{2}N^{\varepsilon_0}}\right).
\end{equation} 
\end{theo}
\begin{remark}
In the assumption \ref{eq:s6.4.2} on $\Omega_N$ we assumed that it coincides with 
an $\Omega$ with smooth boundary, which is independent of $N$, inside a tube of radius $\tau$ around $p(S^1)$. Therefore, Assumption \ref{eq:s6.4.2} implies that $\ell(N) \geq 1/C > 0$. However, the proof of Theorems \ref{thm:t1} and  \ref{thm:t2} 
shows that we can allow $\Omega$ to be $N$ dependent as long as its boundary 
$\partial\Omega$ remains uniformly 
Lipschitz in the sense discussed at the beginning of Section \ref{sec:c1} and satisfies 
($\Omega$\ref{O1})--($\Omega$\ref{O3}). Hence, Theorem \ref{thm:t2} holds as well for 
sets $\Omega_N$, satisfying \eqref{eq:s6.4.1}-\eqref{eq:s6.4.4} with 
 \begin{equation}
	 \frac{C}{N} \leq \ell(N).
\end{equation} 
\end{remark}
\section{Convergence of the empirical measure}\label{sec:LP}
In this section we present the two proofs of Corollary \ref{thm:t0}. The first one, 
in Section \ref{Sec:C1}, shows that it is a consequence of Theorem \ref{thm:t1}. The 
second (alternative) proof in Sections Section \ref{Sec:C2}, Section \ref{Sec:C3}, 
shows how one can obtain the result from our methods via analysing the convergence 
of the associated logarithmic potentials, in perhaps a more direct way. 
\subsection{Proof of Corollary \ref{thm:t0}}\label{Sec:C1} 
Let $\Omega $ be a fixed
domain as in Theorem \ref{thm:t1} and choose a sequence $\delta =\delta _N$
satisfying \eqref{int.3}. By the Borel-Cantelli lemma, we know that
a.s.\ (almost surely)
\begin{equation}\label{pfcor.1}
\frac{1}{N}\# (\sigma (P_\delta ^N)\cap \Omega )\to
\mathrm{vol\,}(p^{-1}(\Omega )\cap S^1),\ N\to \infty .
\end{equation} 
Let now $\Omega $ be a square of the form $a_1\le \mathrm{Re}\, z <a_2$, $b_1\le
\mathrm{Im}\,z <b_2$, $a_2-a_1=b_2-b_1>0$. Assume that the corners $a_j+ib_k$ do
not belong to $p(S^1)$. Then the conditions ($\Omega $1)--($\Omega $3)
make sense. If they are fulfilled, then (\ref{pfcor.1}) holds
a.s.. Indeed, let $\Omega _\mathrm{int}$, $\Omega _\mathrm{ext}$ be
sets with smooth boundary such that $\Omega _\mathrm{int}\subset
\Omega \subset \Omega _\mathrm{ext}$ and coinciding with $\Omega $
away from a small neighborhood of the union of the corners of $\Omega
$. Then (\ref{pfcor.1}) holds a.s.\ for $\Omega _\mathrm{int}$ and
$\Omega _\mathrm{ext}$, and the common limit in the right hand side is
$(2\pi )^{-1}\mathrm{vol\,}(p^{-1}(\Omega )\cap S^1)$. Since 
$$
\frac{1}{N}\# (\sigma (P_\delta ^N)\cap \Omega_\mathrm{int} )\le \frac{1}{N}\#
(\sigma (P_\delta ^N)\cap \Omega )\le \frac{1}{N}\# (\sigma (P_\delta ^N)\cap \Omega_\mathrm{ext} ),
$$
we conclude that (\ref{pfcor.1}) holds a.s.\  for $\Omega $.
\par 
Write $p(\zeta )=p_1(\zeta )+ip_2(\zeta )$ so that
${{p_j}_\vert}_{S^1}$ are real analytic. Then for $j=1,2$:
\begin{itemize}
\item[1)] The set $C_j$ of critical values of ${{p_j}_\vert}_{S^1}$ is
  finite.
\item[2)] For $j=1,2$ and for every $a\in\R$ the equation
  $p_j(\zeta )=a$ has at most finitely many solutions in $S^1$. 
\end{itemize}
\par 
Let $\epsilon >0$. Then we can choose $a,b\in {\R}$ (depending
on $\epsilon $) such that $a+{\Z}\epsilon \cap C_1=\emptyset $,
$b+{\Z}\epsilon \cap C_2=\emptyset $. After a slight shift of $b$
we can arrange so that we also have 
$$
(a+{\Z}\epsilon )+i(b+{\Z}\epsilon )\cap p(S^1)=\emptyset .
$$
Then for each $\epsilon >0$ we have a.s.\  that (\ref{pfcor.1}) holds for
$\Omega =\Omega _{\epsilon ,j,k}$ for all $j,k\in {\Z}$. Here, we
put $\Omega _{\epsilon ,j,k}=(a+[j,j+1[\epsilon )+i(b+[k,k+1[\epsilon
[)$. Let $\epsilon _\nu >0$, $\nu \in {\N}$ be a decreasing sequence
tending to zero. Then a.s., (\ref{pfcor.1}) holds for all the $\Omega
_{\epsilon _\nu ,j,k}$. 
\par 
Let $G$ be the set of all step functions of the form, 
\begin{equation}\label{pfcor.2}
\psi =\sum _{j,k}g_{j,k}1_{\Omega _{\epsilon _\nu ,j,k}},\ g_{j,k}\in
{\mathds{Q}},
\end{equation}
Then a.s.\  we have for every $\psi \in G$, that
\begin{equation}\label{pfcor.3}
\int \psi \,\xi _N(dz)\to \int \psi\, p_*\left(\frac{1}{2\pi }L_{S^1}\right)(dz),\
N\to \infty .
\end{equation}
\par
Let $\phi \in C_c ({\C};{\R})$. For every $\epsilon >0$,
we can find $\psi =\psi _\epsilon \in G$, such that  $|\phi -\psi |\le
\epsilon $. $\xi _N$ and $p_*((2\pi )^{-1}L_{S^1})$ are probability
measures, so 
$$
\left|\int \phi\, \xi _N(dz)-\int\psi\, \xi _N (dz)\right|\le \epsilon ,
$$
$$
\left|\int \phi\, p_*\left(\frac{1}{2\pi }L_{S^1}\right)(dz)-\int\psi \,p_*\left(\frac{1}{2\pi }L_{S^1}\right)(dz)\right|\le \epsilon .
$$
It follows that a.s., we have for all $\phi \in C_0({\C})$,
\begin{multline*}
\limsup_{N\to \infty } \left|
\int \phi\, \xi _N(dz)-\int \phi \,p_*\left(\frac{1}{2\pi }L_{S^1}\right)(dz)
\right|\le 
\\ 2\epsilon +\limsup_{N\to \infty }\left|\int \psi \,\xi _N(dz)-\int \psi \,p_*\left(\frac{1}{2\pi
}L_{S^1}\right)(dz)\right| .
\end{multline*}
A.s.\ the last limit is $0$ for all $\psi \in G$, hence a.s.\  we have
that for all $\epsilon >0$ and all $\phi \in C_c({\C})$,
$$
\limsup_{N\to \infty } \left|
\int \phi \,\xi _N(dz)-\int \phi \,p_*\left(\frac{1}{2\pi }L_{S^1}\right)(dz)
\right|\le 2\epsilon .
$$
In other words, a.s.\  we have  
$$
\lim_{N\to \infty }\int \phi\, \xi _N (dz)=\int \phi\, p_*\left(\frac{1}{2\pi }L_{S^1}\right)(dz),
$$
for all $\phi \in C_c({\C})$, so a.s.:
$$
\xi _N(dz)\rightharpoonup p_*\left(\frac{1}{2\pi }L_{S^1}\right),\ N\to \infty .
$$
Notice that almost surely, $\supp \xi_N$ is contained in a fixed compact set. 
\subsection{Logarithmic potential and weak convergence of measure}\label{Sec:C2}
We begin by recalling some basic facts concerning the weak convergence of measures. 
Let $\mathcal{P}(\C)$ denote the space of probability measures $\mu$ on 
$\C$, integrating the logarithm at infinity
\begin{equation}\label{eq:em1}
	\int \log (1 + |x|) \mu(dx) < + \infty.
\end{equation}
We define the \emph{logarithmic potential} of $\mu$ by 
\begin{equation}\label{eq:em1.2}
	U_{\mu}(z) \defeq - \int \log | z- x | \mu(dx).
\end{equation}
Since $U_{\mu} \in L^1_{\mathrm{loc}}(\C,L(dz))$, it follows that $U_{\mu}(z) < +\infty$ for Lebesgue 
almost every (a.e.) $z\in\C$. 
\par
One property of the logarithmic potential is that for a given sequence of probability measures 
$\{\mu_n\}_n \in \mathcal{P}(\C)$, satisfying some suitable uniform integrability assumption, 
one has that almost sure convergence of the associated logarithmic potentials 
$U_{\mu_n}(z)\to U_{\mu}(z)$, for some $\mu \in \mathcal{P}(\C)$, implies the weak  convergence $\mu_n\rightharpoonup \mu$. 
\par
There are various versions of the above observation known in the case of random measures, 
see for instance \cite[Theorem 2.8.3]{Ta12} or \cite{BoCh13}. In the following we describe   
a slightly modified version of \cite[Theorem 2.8.3]{Ta12} for the reader's convenience. 
\begin{theo}\label{thm:h0}  Let $K,K'\Subset \C$ be open relatively compact sets with 
$\overline{K}\subset K'$, 
and let $\{\mu_n\}_{n\in\N} \in \mathcal{P}(\C)$ be as sequence of random measures so that 
almost surely 
\begin{equation}\label{eq:em1.2.4}
	\supp \mu_n \subset K \hbox{ for } n \hbox{ sufficiently large}.
\end{equation}
Suppose that for a.e. $z\in K'$ almost surely
\begin{equation}\label{eq:em1.2.5}
	U_{\mu_n}(z)\to U_{\mu}(z), \quad n \to \infty, 
\end{equation}
where $\mu \in \mathcal{P}(\C)$ is some probability measure with 
$\supp \mu \subset K$. Then, almost surely, 
\begin{equation}\label{eq:em1.2.6}
	\mu_n\rightharpoonup \mu, \quad n \to \infty, \quad \hbox{weakly.}
\end{equation}
\end{theo}
\begin{proof}
\emph{1.} 
Notice that the assumption that for a.e. $z\in K'$  \eqref{eq:em1.2.5} holds almost surely is equivalent to the statement that almost surely \eqref{eq:em1.2.5} holds for a.e. $z\in K'$. To see this, consider the set 
$E = \{ (z,\omega)\in K'\times \Omega; U_{\mu_n}(z)\to U_{\mu}(z), \hbox{ as }n \to \infty\}\subset K'\times\Omega$, where $\Omega$ denotes the underlying probability space. Applying the Tonelli theorem to 
$\mathbf{1}_{E^c}$ lets us conclude the claim. 
\\
\\ 
\emph{2.}
Since $\log|\cdot -w| \in L^2(K')$ uniformly for $w\in K'$, it follows by the Minkowski 
integral inequalities that, almost surely, $U_{\mu_n},U_{\mu}\in L^2(K')$ uniformly. Let 
us remark here that although $\mu_n$ depends on the random parameter $\omega$, 
we do not denote that explicitly. 
\par 
Combining this with \eqref{eq:em1.2.4} and \emph{step 1.} above, we see that 
there exists an $\Omega'\subset \Omega$ with 
$\mathds{P}(\Omega') = 1$, so that for each $\omega\in \Omega'$ we have that 
\begin{itemize}
	\item \eqref{eq:em1.2.5} holds for a.e. $z\in K'$, 
	\item there exists an $n_0\geq 1$ such that $\supp \mu_n \subset K$ for all $n\geq n_0$, 
    \item there exists a $C_{K',\Omega'}>0$, depending only on $K'$ and $\Omega'$, 
    such that $\|U_{\mu_n}\|_{L^2(K')},\|U_{\mu}\|_{L^2(K')}\leq C_{K',\Omega'}$ 
    for any $n\geq 1$ .
\end{itemize}
To show \eqref{eq:em1.2.6} for any $\omega\in\Omega'$, it is 
enough to show that for any real-valued smooth function $\phi \in \mathcal{C}_c^{\infty}(K';\R)$ 
with support contained in $K'$, 
\begin{equation}\label{eq:em1.5}
	\mu_n(\phi) \to \mu(\phi), \quad n\to \infty.
\end{equation}
\emph{3.} Let $\omega\in\Omega'$, and set $g^M_n(z)=\min (|U_{\mu_n}(z) - U_{\mu}(z)|,M)$, 
$z\in K'$, for $M>0$. The dominated convergence theorem shows that $g_n^M\to 0$, as $n\to\infty$, in 
$L^1(K')$ for any $M>0$. Using the $L^2(K)$ bound of $U_{\mu_n}$ and $U_{\mu}$, we see that 
\begin{equation*}
\begin{split}
	\| g^M_n - |U_{\mu_n} - U_{\mu}| \|_{L^1(K')} &\leq 
	\int_{\substack{|U_{\mu_n}- U_{\mu}|\geq M\\ z\in K}}|U_{\mu_n}(z)- U_{\mu}(z)|L(dz) \\ 
	& \leq \sqrt{2}C_{K'} \left( \int_{\substack{|U_{\mu_n}- U_{\mu}|\geq M\\ z\in K}}L(dz)\right)^{1/2} \\ 
	&\leq \frac{\sqrt{2}C_{K'}L(K')^{1/2} }{M}.
\end{split}
\end{equation*}
Hence, for any $w\in\Omega'$ we have that $U_{\mu_n}\to U_{\mu}$ in $L^1(K')$ as $n\to\infty$. 
Thus, almost surely $U_{\mu_n}\rightharpoonup  U_{\mu}$ in $\mathcal{D}'(K')$, and so \eqref{eq:em1.5} holds almost surely, since 
$\Delta_z U_{\mu_n} =-2\pi \mu_n$, $\Delta_z U_{\mu} =-2\pi \mu$ in $\mathcal{D}'(\C)$.
\end{proof}
\subsection{Proof of Corollary \ref{thm:t0}}\label{Sec:C3}
Recall the definition of the empirical measure $\xi_N$ \eqref{eq:cm1} 
and \eqref{i0.1}. By \eqref{int.0.1}, \eqref{int.2} and the Fourier transform $\mathcal{F}$ 
\eqref{g18} we see that the operator norm of the unperturbed operator $P^0_N$ is 
satisfies 
\begin{equation}\label{eq:em2.1}
	\| P^0_N\| \leq \| p \|_{L^{\infty}(S^1)}.
\end{equation}
Suppose \eqref{gp17}, then by \eqref{gp2}, \eqref{eq:em2.1} it follows that 
\begin{equation}\label{eq:em2.2}
	\| P^{\delta}_N\| \leq \| p \|_{L^{\infty}(S^1)}+1
\end{equation}
for $N>1$ sufficiently large, with probability $\geq 1 - \e^{-N^2}$. We deduce by 
a Borel-Cantelli argument that almost surely 
\begin{equation}\label{eq:em2.3}
	\supp \xi_N \subset  \overline{D(0,\| p \|_{L^{\infty}(S^1)}+1)} \defeq K\subset 
	D(0,\| p \|_{L^{\infty}(S^1)}+2)\defeq K'
\end{equation}
for $N$ sufficiently large. For $p$ as in \eqref{int.0.1}, define the 
probability measure 
\begin{equation}\label{eq:em2.4}
	 \xi =p_*\left(\frac{1}{2\pi} L_{S^1}\right)
\end{equation}
which has compact support, 
\begin{equation}\label{eq:em2.5}
	\supp \xi = p(S^1) \subset K.
\end{equation}
Here, $\frac{1}{2\pi} L_{S^1}$ denotes the normalized Lebesgue measure on $S^1$. 
\par
To conclude Corollary \ref{thm:t0} from Theorem \ref{thm:h0} it remains to show that 
for almost every $z\in K'$ we have that $U_{\xi_N}(z) \to U_{\xi}(z)$ almost surely. 
\par
By \eqref{eq:em1.2} we see that for $z\notin \spec( P_N^{\delta})$ 
\begin{equation}\label{eq:em2.6}
	U_{\xi_N}(z) = -\frac{1}{N} \log| \det (P_N^{\delta} -z)|.
\end{equation}
For any $z\in \C$ the set $\Sigma_z = \{ Q \in \C^{N\times N}; \det(P_0 +\delta Q -z) =0\}$ 
has Lebesgue measure $0$, since $\C^{N\times N} \ni Q \mapsto \det (P_N^{\delta} -z)$ is 
analytic and 
not constantly $0$. Thus $\mu_N(\Sigma_z)=0$, where $\mu_N$ is the Gaussian 
measure given in \eqref{gp1.0}, and for every $z\in \C$ \eqref{eq:em2.6} holds 
almost surely (a.s.). 
\par
Next, define the set 
\begin{equation}\label{eq:em2.7}
	E_N \defeq \{z\in \C; \dist(z,p(S^1)) \leq 1/(CN) \} 
\end{equation}
which has Lebesgue measure $L(E_N) = \mO(N^{-1})$. By \eqref{gp19}, \eqref{gp20}, \eqref{gp2} as well as Proposition \ref{gp:prop1} and \eqref{gp17} we have that for 
every $z\in K'\backslash E_N$ 
\begin{equation}\label{eq:em2.8}
	\left| 
	\frac{1}{N} \log| \det (P_N^{\delta} -z)| - \phi(z)
	\right|
	\leq \mO(\delta N^{3/2}) + N^{-1} \big| \log |\det E_{-+}^{\delta}(z)|\big|.
\end{equation}
with probability $\geq 1 - \e^{-N^2}$. 
Using Proposition \ref{gp:prop1}, we see that for every $z\in K'\backslash E_N$
\begin{equation}\label{eq:em2.9}
	\log |\det E_{-+}^{\delta}(z)| \leq \mO(\log N).
\end{equation}
with probability $\geq 1 - \e^{-N^2}$. Let $\varepsilon_0\in ]0,1[$ be as in Corollary \ref{thm:t0} 
and let $\varepsilon_1\in ]0,1[$ with $ \varepsilon_0< \varepsilon_1$. Then, by replacing $\varepsilon_0$ in  \eqref{c6} with $\varepsilon_1$, we have  that 
\begin{equation}\label{eq:em2.10}
	\log |\det E_{-+}^{\delta}(z)| \geq -N^{\varepsilon_1}
\end{equation}
with probability $\geq 1 - \e^{-N^2} - C_2 \delta^{-|J|} \e^{-\frac{1}{2}N^{\varepsilon_1}}$, when 
\begin{equation*}
	\exp \left[\frac{C_0}{2|J|} - \frac{N^{\varepsilon_1}}{2|J|} \right] \leq \delta .
	\ll N^{-4}
\end{equation*}
For $z\notin p(S^1)$ the function $S^1 \ni\zeta \mapsto \log |z - p(\zeta)|$ is continuous. 
Hence, by \eqref{gp19}, \eqref{eq:em2.4}, \eqref{eq:em1.2}, 
and a Riemann sum argument, we see that  for 
\begin{equation}\label{eq:em2.12}
	\left| 
	\phi(z) + U_{\xi}(z) 
	\right| \longrightarrow 0, \quad \hbox{as } N\to \infty.
\end{equation}
For any $z\in K'\backslash p(S^1)$ we have that $z\in K'\backslash E_N$  for $N>1$ 
suffciently large. Thus, by \eqref{eq:em2.6}, \eqref{eq:em2.8}, \eqref{eq:em2.9}, \eqref{eq:em2.10}, 
and \eqref{eq:em2.12} we have for any $z\in K'\backslash p(S^1)$ and $N>1$ sufficiently 
large that 
\begin{equation}
	\left| U_{\xi_N}(z) - U_{\xi}(z) 	\right| = o(1)
\end{equation}
with probability $\geq 1 - \mO(1) \e^{-\frac{1}{2}N^{\varepsilon_1}(1 - |J| N^{\varepsilon_0-\varepsilon_1})}$. 
Here we also used \eqref{int.3}. Since $\varepsilon_0<\varepsilon_1$, we conclude by the Borel-Cantelli theorem that for almost every $z\in K'$ 
\begin{equation}
	U_{\xi_N}(z)   \longrightarrow  U_{\xi}(z), \quad \hbox{as } N\to \infty, \hbox{ almost surely},
\end{equation}
which by Theorem \ref{thm:h0} concludes the proof of Corollary \ref{thm:t0}.
\end{document}